\documentclass[10pt]{amsart}
\usepackage[latin1]{inputenc}
\usepackage{latexsym} 
\usepackage{amssymb} 
\usepackage{amsbsy} 
\usepackage{eqnarray}
\usepackage{amstext} 
\usepackage{amsmath} 
\usepackage{amsthm} 
\usepackage{amscd} 
\usepackage{graphics} 
\usepackage{color} 
\usepackage{epsfig} 
\usepackage[all]{xy} 

\newtheorem{thm}{Theorem}[section]
\newtheorem{prop}[thm]{Proposition}

\newtheorem{lem}[thm]{Lemma}

\theoremstyle{definition}
\newtheorem{defn}[thm]{Definition}
\newtheorem{ex}[thm]{Example}
\newtheorem{rem}[thm]{Remark}

\newcommand{\End}{\operatorname{End}}
\newcommand{\Hom}{\operatorname{Hom}}

\newcommand{\ind}{\operatorname{ind}}

\newcommand{\add}{\operatorname{add}}

\def\SS{\scriptstyle}
\def\mc{\mathcal}

\begin{document}

\title[$m$-cluster tilted algebras of type $A_n$]{Derived equivalence classification of $m$-cluster tilted algebras of type $A_n$}
\author{Graham J. Murphy}
\address{Dept. Pure Mathematics\\
        School of Mathematics\\
         University of Leeds\\
         Leeds\\LS2 9JT\\U.K.}
\email{graham@maths.leeds.ac.uk}
\thanks{}

\date{22nd August 2007}
\begin{abstract} We use the maximal faces of the $m$-cluster complex of type $A_n$ introduced in \cite{FR} to describe the quivers of $m$-cluster tilted algebras of type $A_n$. We then classify connected components of $m$-cluster tilted algebras of type $A_n$ up to derived equivalence using tilting complexes directly related to the combinatorics of the  $m$-cluster complex of type $A_n$. This generalizes a result of Buan and Vatne \cite{BV}.
\end{abstract}
\maketitle

\section{Introduction}
Cluster categories were introduced in \cite{BMRRT} as a representation theoretic framework for cluster algebras of Fomin and Zelevinsky \cite{FZ}. Independently, in \cite{I} Iyama defines maximal $m$-orthogonal modules for Artin algebras, $m\geqslant 1$, and presents a classification of maximal 1-orthogonal modules for finite type self-injective algebras of tree class $A_n,B_n,C_n$ and $D_n$. Motivated by a generalized version of the combinatorics related to cluster algebras, $m$-cluster categories were defined as a generalization of cluster categories \cite{T}, \cite{Z}. The cluster categories and their generalizations have been the subject of much recent research.
\\The $m$-cluster category is defined to be $\mathcal{C}_m(H):=\mathcal{D}^b(H)/\tau^{-1}[m]$ where $H$ is a finite dimensional hereditary algebra over a field $k$ (which we take to be algebraically closed) whose corresponding quiver has no loops or oriented 2-cycles; $\tau$ is the Auslander-Reiten (AR) translate in the bounded derived category $\mathcal{D}^b(H)$; $[1]$ denotes the suspension functor present in the triangulated structure of $\mathcal{D}^b(H)$ and $[m]$ denotes its $m$-th power, that is $[1]^m=[m]$.
\\In $m$-cluster categories one can define $m$-cluster tilting objects as follows \cite{T}.
\begin{defn} An $m$-cluster tilting object $T$ in $\mc{C}_m(H)$ is an object $T$ satisfying the following conditions.
\begin{enumerate}
\item $\Hom_{\mc{C}_m(H)}(T,X[i])=0$, $1\leqslant i\leqslant m$, if and only if $X\in\add(T)$.
\item $\Hom_{\mc{C}_m(H)}(X,T[i])=0$, $1\leqslant i\leqslant m$, if and only if $X\in\add(T)$.
\end{enumerate}
\end{defn}
It is the endomorphism algebras of these $m$-cluster tilting objects, $\End_{\mathcal{C}_m(H)}(T)$, called $m$-cluster tilted algebras, that we study in this paper. We focus on the case where $H$ is an hereditary algebra of Dynkin type $A_n$.
\\Our aim is to generalize a result of Buan and Vatne, \cite{BV}, on the derived equivalence classification of 1-cluster tilted algebras of type $A_n$.
\\In section \ref{sec:m_clust} we give a local description of $m$-cluster tilted algebras of type $A_n$ as quivers with relations using the combinatorics of $m$-cluster complexes \cite[section 4]{FR} which can be related to $m$-cluster tilting objects. More precisely, it can be proven that the maximal faces, or facets, of the generalized cluster complex correspond bijectively to the $m$-cluster tilting objects in $\mc{C}_m(kA_n)$ \cite{T}, \cite{M}.  Then, the combinatorial model for the generalized cluster complex in type $A_n$ can be used to determine the $m$-cluster tilted algebra corresponding to a given $m$-cluster tilting object.
\\This description of the $m$-cluster tilted algebras based on the combinatorial model for the generalized cluster complex enables us to prove that $m$-cluster tilted algebras of type $A_n$ are gentle and that if the quiver of an $m$-cluster tilted algebra of type $A_n$ contains cycles they must be of length $m+2$ and must have full relations, that is the composition of any two consecutive arrows in the cycle must be a relation. These $m+2$-cycles play a crucial r\^{o}le when considering derived equivalences of $m$-cluster tilted algebras. Details of these statements can be found in section \ref{sec:m_clust}.
\\Our main theorem can be stated as follows.
\begin{thm}\label{main} Two connected components of an $m$-cluster tilted algebra of type $A_n$ are derived equivalent if and only if their quivers have the same number of oriented $m+2$-cycles with full relations.
\end{thm}
The statement of theorem \ref{main} must allow for the possibility that for $m\geqslant 2$ the $m$-cluster tilted algebras of type $A_n$ may be disconnected. That these algebras can be disconnected follows from the combinatorial description in section \ref{sec:m_clust}. This is one of the noticeable differences between the $m=1$ and $m\geqslant 2$ cases, since 1-cluster tilted algebras of type $A_n$ are connected \cite{BV}.
\\After section \ref{sec:m_clust} the paper is structured as follows. In section \ref{sec:tilt} we provide tilting complexes which prove that there exist derived equivalences between $m$-cluster tilted algebras which preserve $m+2$-cycles with full relations. In this section we also describe how these tilting complexes are related to the combinatorial model of $m$-cluster tilted algebras in type $A_n$. Then we use \cite[thm 4.1,cor 4.3]{BH}, which states that the Cartan matrix of a gentle algebra is unimodular equivalent to a diagonal matrix determined by the cycles with full relations which appear in the description of that algebra as a quiver with relations, to identify the necessary conditions for derived equivalence of two connected components of an $m$-cluster tilted algebra of type $A_n$. The form of the diagonal matrix is restated in section \ref{sec:tilt}.
\\Finally, in section \ref{sec:algo}, we give the algorithm which uses the results of section \ref{sec:tilt} to mutate any connected component of an $m$-cluster tilted algebra to a particular $m$-cluster tilted algebra, referred to as the \textit{normal form} (to be defined in section \ref{sec:tilt}), in such a way that a derived equivalence is achieved.
\section{$m$-cluster tilted algebras}\label{sec:m_clust}
We begin with a brief summary of the combinatorial background from \cite[section 4]{FR},\cite[section 4]{I}, \cite{R}.
\\The $m$-co-ordinate system, on the translation quiver $\mathbb{Z}A_n$ is defined as shown in the diagram below.
\begin{center}
\includegraphics[scale=0.55]{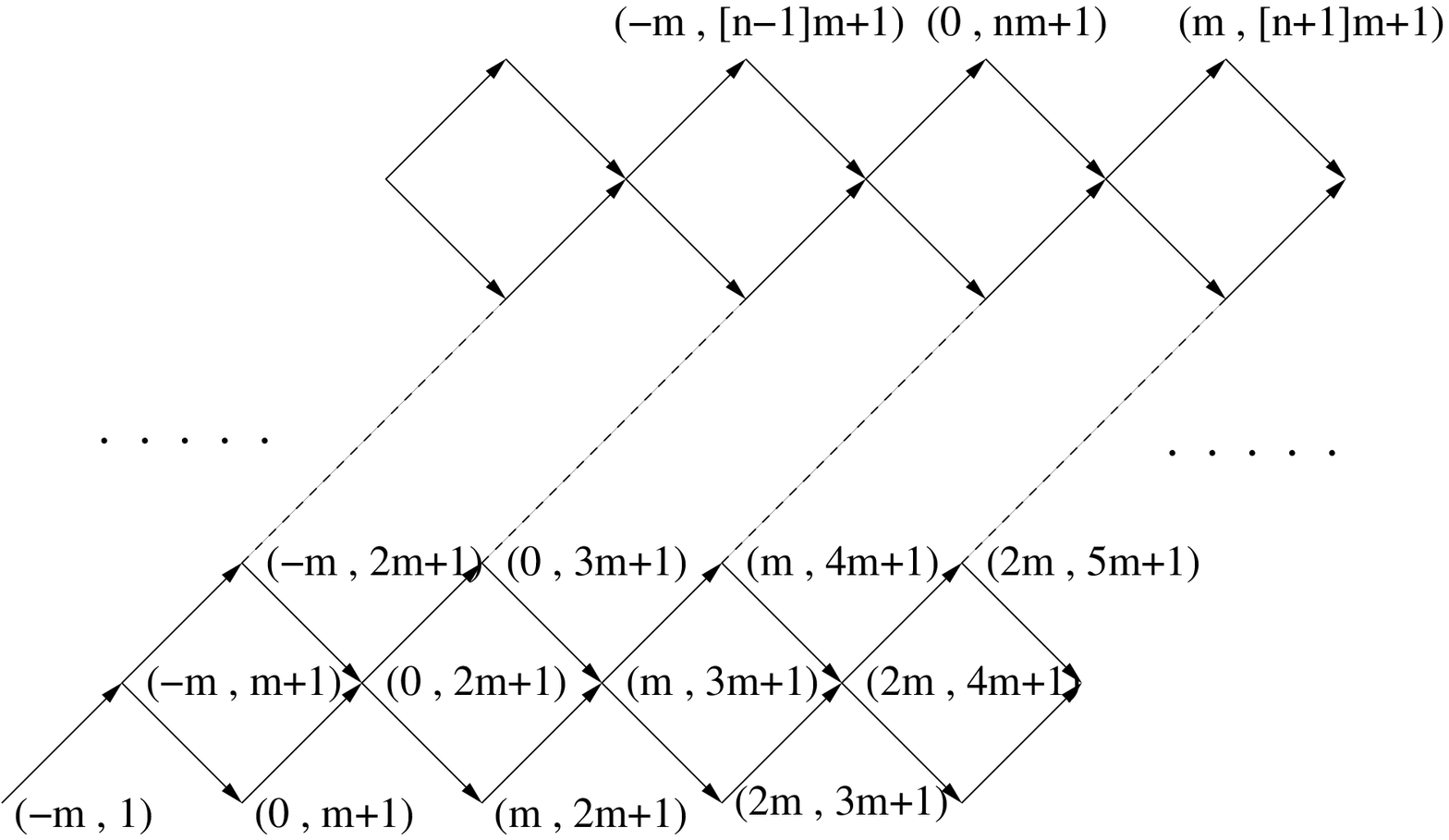}
\end{center}
\begin{rem}This co-ordinate system is similar to the idea of the $m$-th power of the translation quiver defined by Baur and Marsh in \cite{BMa}.
\end{rem}
We now recall the definition of the mesh category associated to the translation quiver $\mathbb{Z}A_n$. The objects of this category are the vertices of $\mathbb{Z}A_n$ and morphisms are the arrows of $\mathbb{Z}A_n$, subject to the \textit{mesh relations}. For each arrow $\alpha:x\rightarrow y$ denote by $\sigma(\alpha)$ the unique arrow $\sigma(\alpha):\tau(y)\rightarrow x$. Notice that the uniqueness of $\sigma(\alpha)$ follows from the bijection between the arrows starting at $x$ and the arrows ending at $x$.
The mesh relations are given by,
$$
\sum_{\alpha:x\rightarrow y}\alpha\sigma(\alpha)=0
$$
for each vertex $y\in(\mathbb{Z}A_n)_0$.
\\Next for a vertex $x=[rm,sm+1]\in\mathbb{Z}A_n$, where $r\in\mathbb{Z}$ and $r+1\leqslant s\leqslant r+n$, $H^-(x)$ and $H^+(x)$ as the vertices of $\mathbb{Z}A_n$ enclosed by and on the boundary of the regions specified in the following figure.
\\[10pt]
\begin{center}
\includegraphics[scale=0.5]{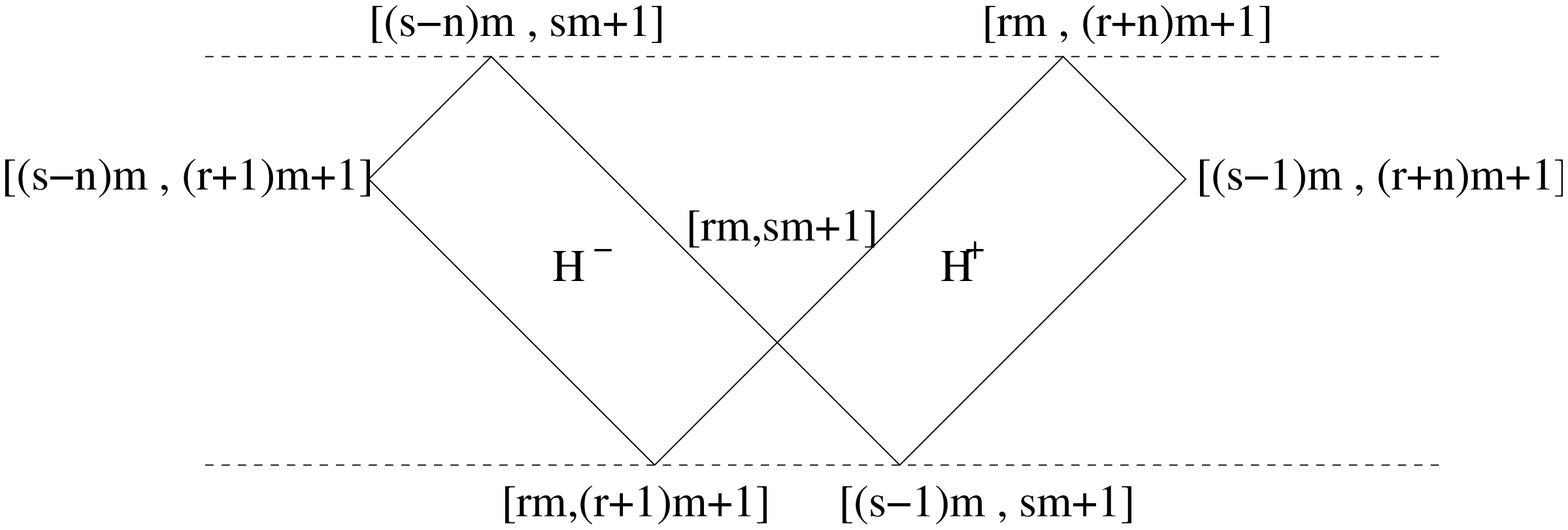}
\end{center}
\begin{rem}
For any $x\in\mathbb{Z}A_n$ by the \textit{left boundary} of $H^-(x)$ we mean the set of all vertices on the boundary of $H^-(x)$ to the left of the vertices  $[rm,(r+1)m+1]$ and $[(s-n)m,sm+1]$. The \textit{right boundary} of $H^-(x)$ is the set of all vertices on the boundary of $H^-(x)$ to the right of $[rm,(r+1)m+1]$ and $[(s-n)m,sm+1]$. We define left and right boundaries of $H^+(x)$ similarly. Notice also that the regions $H^-(x)$ and $H^+(x)$ on $\mathbb{Z}A_n$ describe the set of vertices from which (respectively, to which) there is a non-zero morphism in the mesh category associated to $\mathbb{Z}A_n$.
\end{rem}
Define the following automorphisms on $\mathbb{Z}A_n$:
\begin{eqnarray}
\omega : \mathbb{Z}A_n\rightarrow\mathbb{Z}A_n & , & \tau : \mathbb{Z}A_n\rightarrow\mathbb{Z}A_n \nonumber \\
(i,j)\mapsto (j-(n+1)m,i+1)& {} &(i,j)\mapsto(i-m,j-m). \nonumber
\end{eqnarray}
Further, $\tau_{i}:=\tau\omega^{i-1}$.
\begin{rem}These automorphisms were considered by Iyama in \cite{I} for $m=1$.
\end{rem}
Results by D. Happel show that for a quiver $Q$ of Dynkin type $\Delta$ the AR-quiver of the bounded derived category $\mc{D}^b(kQ)$ is isomorphic as a stable translation quiver to $\mathbb{Z}\Delta$ \cite[I,5.6]{Ha}. From this result it follows that $\omega$ as defined above in type $A_n$ is the action of the \textit{inverse} of the suspension $[1]$ in $\mc{D}^b(kA_n)$ expressed in terms of the $m$-co-ordinate system. Also the automorphism $\tau$ on $\mathbb{Z}A_n$ is the action of the AR-translate in $\mc{D}^b(kA_n)$ expressed in terms of the $m$-co-ordinate system. It is also proven in \cite[I,5.6]{Ha} that for $Q$ a quiver of Dynkin type the mesh category $k(\mathbb{Z}\Delta)$ is equivalent to $\ind \mc{D}^b(kQ)$, where $Q$ is of type $\Delta$. Here $\ind \mc{D}^b(kQ)$ denotes the subcategory of $\mc{D}^b(kQ)$ whose objects are the indecomposable objects of $\mc{D}^b(kQ)$ and whose morphisms are the irreducible morphisms between indecomposable objects in $\mc{D}^b(kQ)$.
\\We can construct the quotient translation quiver $\mathbb{Z}A_n/\langle\tau_{m+1}\rangle$ by identifying the vertices of $\mathbb{Z}A_n$ with their $\tau_{m+1}$-shifts.
\\
\\[10pt]Now we present some polygonal combinatorics which appears in \cite{FR} in connection with the generalized cluster complexes of type $A_n$. We denote by $P$ a (regular convex) $N$-gon where $N\in\mathbb{N}$.
\begin{defn}An $m$-allowable diagonal in $P$ is a chord joining two non-adjacent
boundary vertices which divides $P$ into two smaller polygons $P_1$ and $P_2$ which can themselves be subdivided into
$m+2$-gons by non-crossing chords.
\end{defn}
We label the vertices of $P$ from 0 to $N-1$ in an anti-clockwise direction, then we denote diagonals in
$P$ by $d(i,j)$ where $d(i,j)$ meets the $P$ at the vertices labeled $i$ and $j$.
\begin{rem}In type $A_n$ the regular polygon $P$ will be an $(n+1)m+2$-gon. Notice however that the definition of $m$-allowable diagonals and also the next lemma \ref{div} do not require us to choose a specific value for $N$.
\end{rem}
The following lemma is easy to prove.
\begin{lem}\label{div}A regular convex $N$-gon, $P$, can be divided into convex $m+2$-gons by non-crossing diagonals if, and only
if, $N\equiv 2\,\, mod\,\, m$.
\end{lem}
Lemma \ref{div} implies that an $(n+1)m+2$-gon, $P$, can be divided into $m+2$-gons by non-crossing $m$-allowable diagonals.
\begin{rem}
We can define a simplicial complex with vertex set the $m$-allowable diagonals of $P$. Two $m$-allowable diagonals are
\textit{compatible} if they do not cross in the interior of $P$. The $m$-cluster complex of type $A_n$ is isomorphic as a simplicial complex to the clique complex with respect to the above notion of compatibility
for the set of $m$-allowable diagonals of
an $(n+1)m+2$-gon \cite[sec. 4.1]{FR}. The original definition of $m$-cluster complex of type $A_n$ given in \cite{FR} by Fomin and Reading in terms of ``generalized'' root
systems.
\end{rem}
We will need the next proposition which appears in \cite{M} in the proof of proposition \ref{endoquiv}.
\begin{prop}\label{bij}There is a bijection
$$\theta_{A_n}:\mathcal{D}\longrightarrow \left(\mathbb{Z}A_n /\langle\tau_{m+1}\rangle\right)_0$$
given by $\theta_{A_n}(d(i,j))=(i,j)$ where $\left(\mathbb{Z}A_n /\langle\tau_{m+1}\rangle\right)_0$ denotes the set of vertices of the translation quiver $\left(\mathbb{Z}A_n /\langle\tau_{m+1}\rangle\right)$ and $\mathcal{D}$ denotes the set of all $m$-allowable diagonals of a regular convex $(n+1)m+2$-gon.
\\Also, we have that the subset,
$$(\theta_{A_n}^{-1}\circ\pi)\Big(\bigcup_{1\leqslant c\leqslant m}H^+\big(\tau_c^{-1}(i,j)\big)\big)\Big)\subseteq\mathcal{D}
$$
(where $\pi:\mathbb{Z}A_n\rightarrow \mathbb{Z}A_n /\langle\tau_{m+1}\rangle$ is the canonical projection) is the set of $m$-allowable diagonals crossing $\theta^{-1}_{A_n}((i,j))=d(i,j)$ in $P$.
\end{prop}
\begin{rem}Notice that if any collection of $m$-allowable diagonals in a division of $P$ share a given vertex then when considered as vertices of $\mathbb{Z}A_n/\langle\tau_{m+1}\rangle$ under $\theta_{A_n}^{-1}$ they must have a common co-ordinate. This will be used in the proof of \ref{isoclust}.
\end{rem}
The next result appears in \cite[sec. 1]{BMRRT} (for $m=1$, but which holds for any value of $m\geqslant 1$) relates the mesh category associated to the quotient translation quiver $\mathbb{Z}A_n/\langle\tau_{m+1}\rangle$ to the full subcategory category $\ind(\mc{C}_m(kA_n))$ of $\mc{C}_m$.
\begin{prop}\label{meshcor}Let $Q$ be any quiver of Dynkin type $\Delta$. Then $\ind(\mc{C}_m(kQ))$ is equivalent to the mesh category associated with $\mathbb{Z}\Delta/\langle\tau_{m+1}\rangle$.
\end{prop}
\begin{rem}\label{dimrem}
Notice also that it follows from the mesh relations and the equivalence in proposition \ref{meshcor} that the $\Hom$ spaces in the $m$-cluster category of type $A_n$ are at most one dimensional.
\end{rem}
From now on let $P$ denote a regular convex polygon with $(n+1)m+2$ vertices. Any division of $P$ by $m$-allowable diagonals will be assumed to be maximal and all $m$-cluster-tilted algebras will be of type $A_n$. We will denote by $\mathcal{C}_m$ the $m$-cluster category $\mathcal{D}^b(kA_n)/\tau^{-1}[m]$. Given a division, $\mc{T}$, of a regular convex $(n+1)m+2$-gon into $m+2$-gons by non-crossing $m$-allowable diagonals we describe a finite quiver, $Q_{\mathcal{T}}$. We also describe a set of relations, $\mc{I}_{\mc{T}}$, in the path algebra, $kQ_{\mathcal{T}}$, and show the algebra $kQ_{\mathcal{T}}/\mc{I}_{\mc{T}}$ is isomorphic to an $m$-cluster-tilted algebra. We will show that all $m$-cluster-tilted algebras can be realized in this way.
\\
\\[10pt]Let $\mathcal{T}$ denote the set of $m$-allowable diagonals in a division of $P$. The vertices of the quiver $Q_{\mathcal{T}}$
 are in one-to-one correspondence with the elements of $\mathcal{T}$. For any two vertices $i$ and $j$ in $Q_{\mathcal{T}}$, there is an arrow $i\rightarrow j$ if and only if:
 \begin{enumerate}
 \item the corresponding $m$-allowable diagonals, $d_i$ and $d_j$ share a vertex in the $(n+1)m+2$-gon.
 \item $d_i$ and $d_j$ are edges of the same $m+2$-gon in the division.
 \item $d_j$ follows $d_i$ in a clockwise direction.
 \end{enumerate}
\begin{ex} Here  $Q_{\mathcal{T}}$ is shown for a particular division of a (5+1)+2=8-gon.
\begin{center}
\includegraphics[scale=0.5]{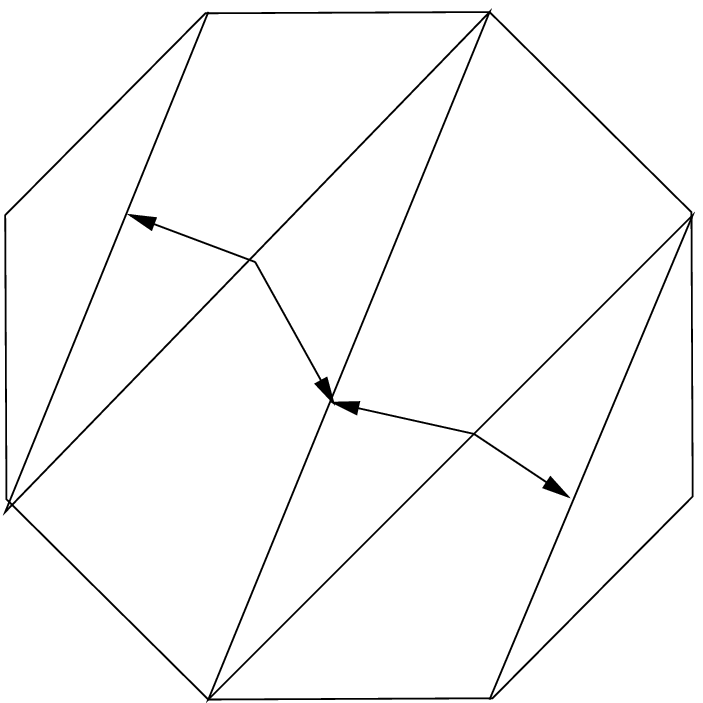}
\end{center}
\end{ex}
In \cite{T} the following statement is proven.
\begin{prop}\label{bij1}
There exists a bijection between the maximal divisions of $P$ by $m$-allowable diagonals and the $m$-cluster-tilting objects in the cluster category $\mc{C}_m$.
\end{prop}
The result also follows from \cite[thm 1.01]{M}.
\\For a given division $\mc{T}$ of $P$ let $T$ denote the corresponding $m$-cluster-tilting object under the bijection of proposition \ref{bij1} and let $\End_{\mc{C}_m}(T)\cong kQ/I$ denote the $m$-cluster-tilted algebra.
\begin{prop}\label{endoquiv} The quiver $Q_{\mathcal{T}}$ is the quiver, $Q$, where $\End_{\mc{C}_m}(T)\cong kQ/I$.
\end{prop}
\begin{proof}For an $m$-allowable diagonal $d_i$ in $\mc{T}$ let $T_i$ denote the indecomposable object of $\mc{C}_m$ which corresponds to the vertex $\theta_{A_n}^{-1}(d_i)$ of the mesh category associated with $\mathbb{Z}A_n/\langle\tau_{m+1}\rangle$.
\\Since by proposition \ref{bij} we have a bijection $\theta_{A_n}:\mathcal{D}\longrightarrow \left(\mathbb{Z}A_n /\langle\tau_{m+1}\rangle\right)_0$ and, by proposition \ref{meshcor}, the category $\ind(\mathcal{C}_m)$ is equivalent to the mesh category of $\mathbb{Z}A_n/\langle\tau_{m+1}\rangle$ it follows that if two $m$-allowable diagonals $d_1$ and $d_2$ share a vertex of $P$ then there must be a non-zero map, $\nu$ say, between the corresponding indecomposable summands $T_1$ and $T_2$ of the $m$-cluster-tilting object, $T$. Indeed, if two $m$-allowable diagonals, $d_1$ and $d_2$, share a vertex of $P$ then they correspond under $\theta^{-1}_{A_n}$ to vertices, $x_1$ and $x_2$ say, of $\mathbb{Z}A_n/\langle\tau_{m+1}\rangle$ which have a common co-ordinate. We assume without loss that  $\nu:T_1\rightarrow T_2$.
\\Then $\nu$ will be irreducible if and only if $d_1$ and $d_2$ are part of the same $m+2$-gon in the division of $P$. Notice that if $d_1$ and $d_2$ lie in different $m+2$-gons but, as we assume, share a vertex of $P$ then a third $m$-allowable diagonal, $d_3$, must exist which also meets $P$ at the vertex shared by $d_1$ and $d_2$. There would then exist non-zero maps from $d_1$ to $d_3$ and from $d_3$ to $d_1$, but since there is a common co-ordinate we must have that the composition of the non-zero maps $T_1\rightarrow T_3$ and $T_3\rightarrow T_2$ is equal to $\nu$.
The definition of $\theta_{A_n}$ ensures that the non-zero map between $T_i\rightarrow T_j$ exists if $d_j$ follows $d_i$ in a clockwise direction around the $m+2$-gon.
\\Notice that it is not possible for other arrows to exist since they would have to correspond to irreducible morphisms on the mesh category.
\end{proof}
Next we define $\mc{I}_{\mc{T}}$. Given consecutive arrows $i\stackrel{\alpha}{\rightarrow} j\stackrel{\beta}{\rightarrow} k$ in $Q_{\mathcal{T}}$ we define the path to be zero if $d_i$, $d_j$ and $d_k$ are in the same $m+2$-gon in the division of $P$. Let $\mathcal{I}_{\mathcal{T}}$ denote the ideal in the path algebra $kQ_{\mc{T}}$ generated by such relations.
\begin{ex}\label{relex}Here are two examples showing the quiver $Q_{\mathcal{T}}$ and relations $\mathcal{I}_{\mathcal{T}}$ for two different 2-allowable divisions of a $16=(6+1)\cdot 2+2$-gon.
\begin{center}
\includegraphics[scale=0.5]{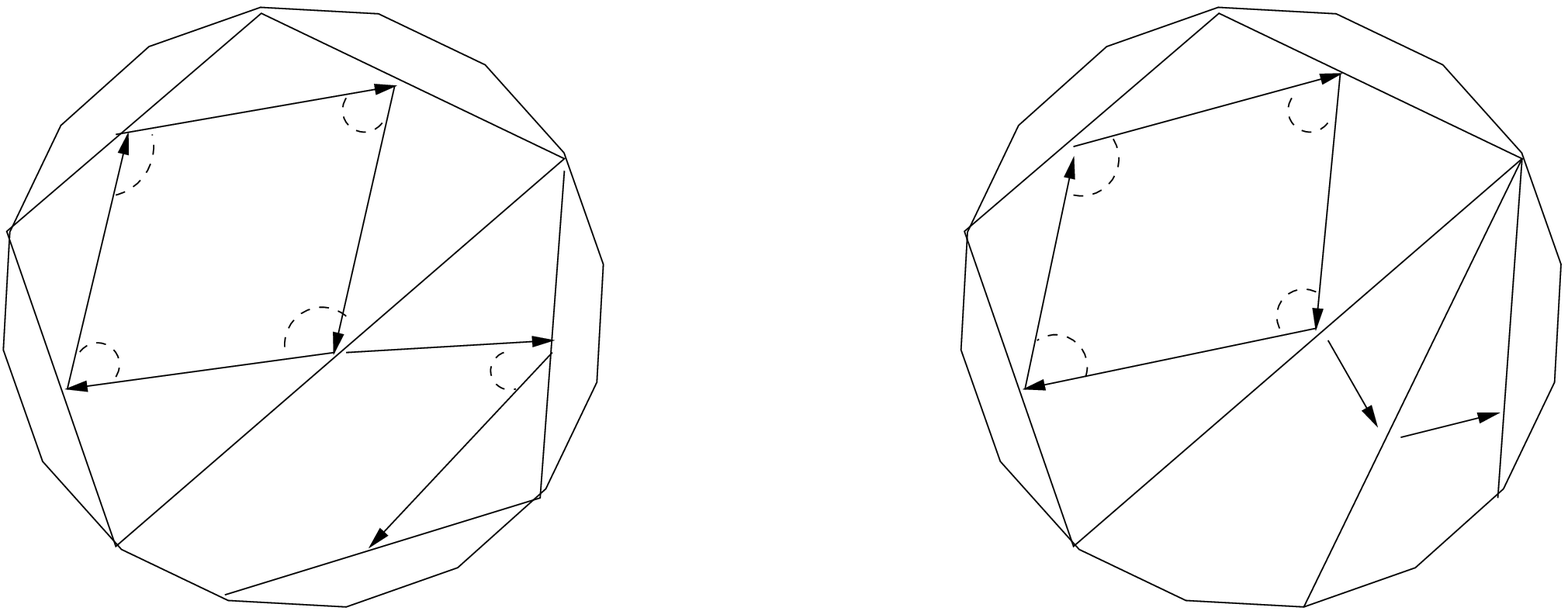}
\end{center}
\end{ex}
The following proposition was known for the $m=1$ case \cite{CCS}, here we state the result for $m\geqslant 1$.
\begin{prop}\label{isoclust} The algebra k$Q_{\mathcal{T}}$/$\mathcal{I}_{\mathcal{T}}$ is isomorphic to the $m$-cluster-tilted algebra $End_{\mathcal{C}_m}(T)$.
\end{prop}
\begin{rem}In the proof of the proposition whenever we refer to $H^+(x)$ we mean $\pi(H^+(x))$, that is we regard the $H^+(x)$ region as a subset of $\mathbb{Z}A_n/\langle\tau_{m+1}\rangle$. Also, when we refer to an $m$-allowable diagonal $d$ we will use the same notation for the $m$-allowable diagonal as the vertex of $\mathbb{Z}A_n/\langle\tau_{m+1}\rangle$
\end{rem}
\begin{proof}
Suppose that $d_i$, $d_j$ and $d_k$ are $m$-allowable diagonals which are part of the same $m+2$-gon arranged as follows,
\begin{center}
\includegraphics[scale=0.5]{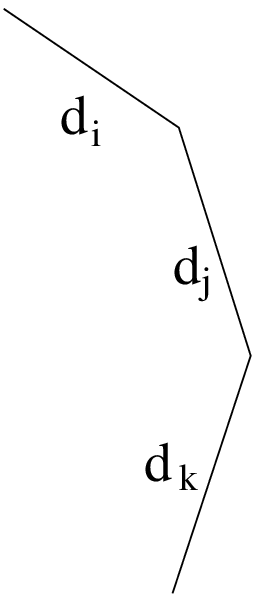}
\end{center}
so that $d_i$ and $d_j$ share a vertex of $P$, $d_j$ and $d_k$ share a vertex of $P$ but $d_i$, $d_j$ and $d_k$ have no common vertex.
\\We know by proposition \ref{endoquiv} that there are arrows $i\stackrel{\alpha}{\rightarrow}j$ and $j\stackrel{\beta}{\rightarrow} k$ in $Q_{\mathcal{T}}$ corresponding to irreducible maps between the indecomposable summands $T_i$ and $T_j$, and $T_j$ and $T_k$ respectively.
\\We claim that the composition of these irreducible maps between $T_i\stackrel{\alpha}{\rightarrow}T_j\stackrel{\beta}{\rightarrow}T_k$ is zero in the $m$-cluster category. It is enough to show that $d_k\notin H^+(d_i)$, the claim then follows from proposition \ref{meshcor}.
\\Since $d_i$ and $d_j$ share a co-ordinate we have that $d_j$ lies on the left boundary of $H^+(d_i)$ and similarly, $d_k$ lies on the left boundary of $H^+(d_j)$. Notice that since $d_i$, $d_j$ and $d_k$ have no common vertex it follows that $d_k$ cannot lie on the part of the left boundary of $H^+(d_j)$ which overlaps with the left boundary of $H^+(d_i)$. Now since $d_k$ does not cross $d_i$ we have that $d_k\notin H^+(\tau^{-1}d_i)$ so that $d_k$ lies on the left boundary of $H^+(d_j)$ but cannot lie in any part of $H^+(d_i)$. Therefore the composition $\beta\circ\alpha$ must be equal to zero in $mc{C}_m$ and so in $\End_{\mathcal{C}_m}(T)\cong kQ/I$ we have $\beta\alpha\in I$.
\\It follows from the above that $\mc{I}_{\mc{T}}\subseteq I$. Finally, we claim that $\mc{I}_{\mc{T}}=I$. To see that this is the case it is  enough to show that no commutativity relations exist in $I$ and that the composition of two arrows in $Q_{\mc{T}}$ not in the same  $m+2$-gon is not a relation.
\\First we consider commutativity relations. Suppose that an $m$-cluster-tilting object $T$ has two indecomposable summands $T_i$ and $T_j$ so that there exist two distinct non-zero paths between the corresponding vertices, $i$ and $j$ of $\End_{\mathcal{C}_m}(T)$. Then the mesh relations dictate that $j\in H^+(\tau^{-1}(i))$ or vice versa, which is impossible since $d_i$ and $d_j$ do not cross.
\\Now suppose that $d_1\stackrel{\gamma}{\rightarrow}d_2$ and $d_2\stackrel{\delta}{\rightarrow}d_1$ are two arrows in $Q_{\mc{T}}$ such that the path $\delta\gamma$ is defined and $\gamma$ and $\delta$ lie in different $m+2$-gons in some division of $P$. Then we must have that $d_1,d_2$ and $d_3$ have a common vertex so that $\theta_{A_n}^{-1}(d_2)$ and $\theta_{A_n}^{-1}(d_3)$ lie on the left boundary of $H^+(\theta_{A_n}^{-1}(d_1))$. Hence the composition $\delta\gamma$ is not a relation in $\End_{\mc{C}_m}(T)$.
\end{proof}
We have now proven that $kQ_{\mathcal{T}}$/$\mathcal{I}_{\mathcal{T}}$ is isomorphic to the $m$-cluster-tilted algebra corresponding to $T$. Since every $m$-cluster-tilting object can be described by a division of $P$ \cite{T} every $m$-cluster-tilted algebra can be realized in this way.
\begin{rem}\label{difference}
Notice that in the proof of proposition \ref{isoclust} if $m=1$ then the $m$-allowable diameters $d_i$, $d_j$ and $d_k$ must form a triangle so that there exists an irreducible map $T_k\rightarrow T_i$ making a three cycle. We collect the following easy consequences of \ref{isoclust},
\begin{enumerate}
\item the only possible cycles which can occur in $kQ_{\mathcal{T}}$/$\mathcal{I}_{\mathcal{T}}$ are $m+2$-cycles with full relations, that is all paths of length two in any $m+2$-cycle are relations.
\item for $m$-cluster-tilted algebras of type $A_n$, $m\neq 1$, relations can occur outwith cycles.
\item there can exist at most $m-1$ consecutive relations outwith a cycle.
\end{enumerate}
It is known for 1-cluster-tilted algebras of Dynkin type that the relations are determined by the quiver \cite{BMR3}. The next example shows that this is not the case for $m$-cluster-tilted algebras of type $A_n$, $m\neq 1$.
\end{rem}
\begin{ex}This example shows that the quiver $\vec{A_2}$ can arise as the quiver of two different 2-cluster-tilted algebras of type $A_3$.
\begin{center}
\includegraphics[scale=0.5]{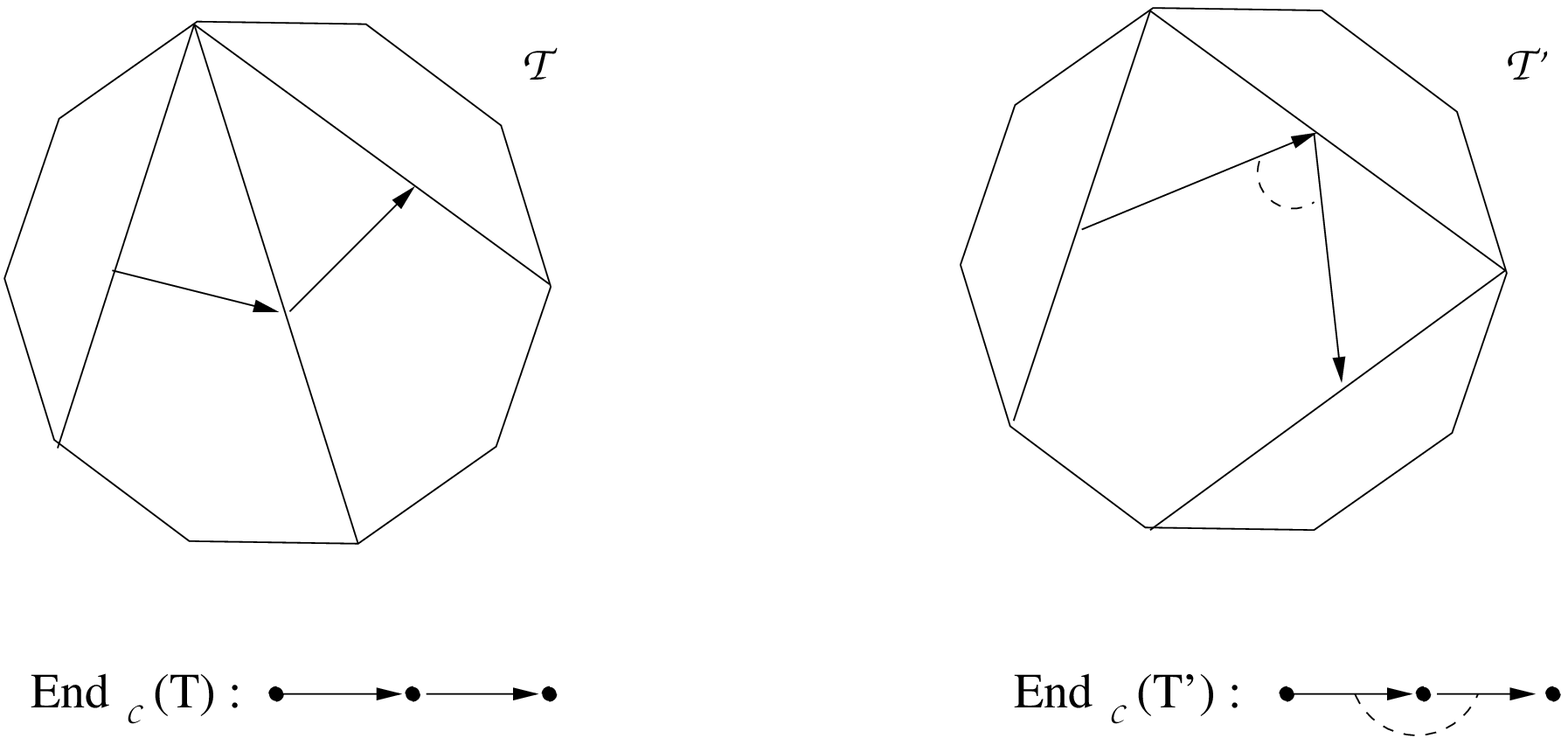}
\end{center}
Note that the figure in example \ref{relex} provides two more examples.
\end{ex}
\section{Tilting complexes and elementary moves}\label{sec:tilt}
The remainder of this paper will demonstrate theorem \ref{main}.
\\We begin with the following observation.
\begin{prop} $m$-cluster-tilted algebras of type $A_n$ are gentle for any $m,n\geqslant 1$.
\end{prop}

\begin{proof} The result follows from considerations of the possible divisions of $P$, the regular $(n+1)m+2$-gon. Let $\mc{T}$ be a division. It is clear that there can be at most two arrows starting or ending at a given vertex of $Q_{\mc{T}}$.
The following figures make the other required properties clear.
$$
\begin{array}{ccc}
\includegraphics[scale=0.5]{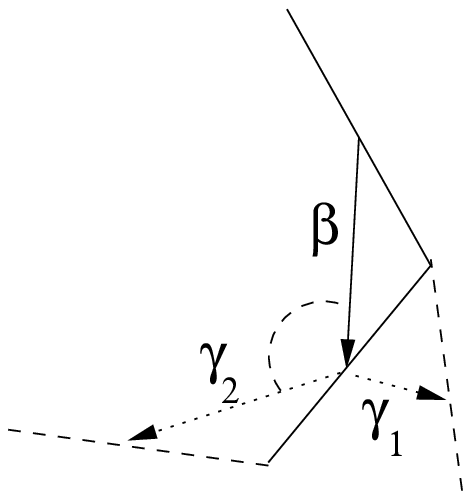}\qquad\qquad &
\includegraphics[scale=0.5]{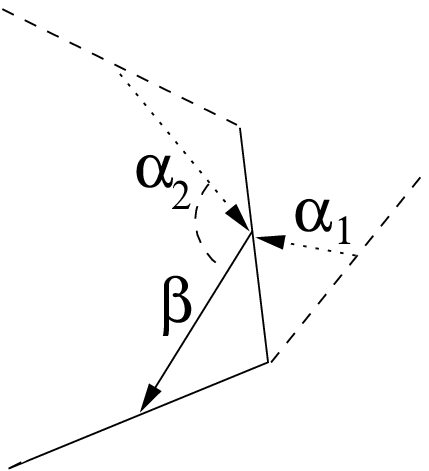}\qquad\qquad &
\includegraphics[scale=0.5]{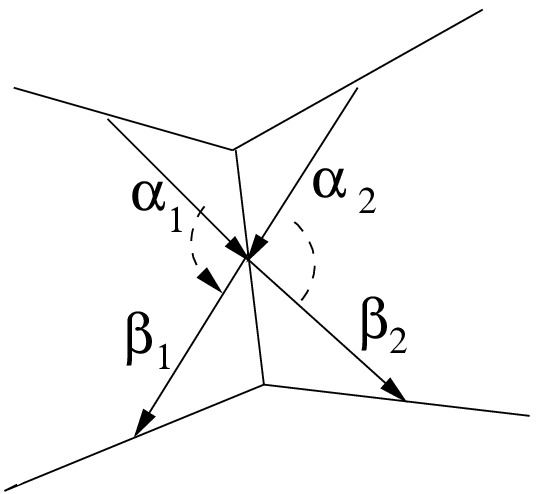}
\end{array}
$$
\end{proof}
\begin{rem}See \cite{AsSk}, for example, for the definition of a gentle algebra.
\end{rem}
Suppose that we have some division of $P$. Notice that each $m$-allowable diagonal, $d$, is a common edge of two $m+2$-gons and that for a given $d_0$ (as in figure the next example) there are $m$-allowable diagonals $d_1,d_2,\ldots,d_m$ of $P$ which we view as anti-clockwise ``rotations" of $d_0$ inside the two $m+2$-gons of which it is an edge. This is a generalization of the well know concept of quadrilateral exchange in triangulations of polygons, indeed quadrilateral exchange is the $m=1$ case.
\\It is clear that for any given $m$-allowable diagonal, $d_0\in\mathcal{T}$, there exist $m$ distinct $m$-allowable diagonals $d_1,d_2,\ldots,d_m$ which could replace $d_0$ and provide a distinct division of $P$.
\begin{ex}The next figure shows an example of the possible replacement 2-allowable diagonals for $d_0$ in a regular hexagon.
\begin{center}
\includegraphics[scale=0.5]{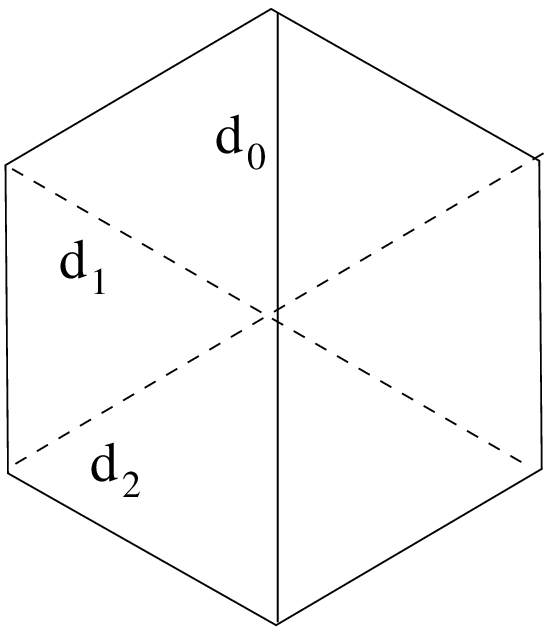}
\end{center}
\end{ex}
\begin{defn}\label{elmove} Suppose that $\mc{T}$ is a maximal division of $P$ by $m$-allowable diagonals. Fix a presentation of $P$ in the plane. For any chosen $m$-allowable diagonal, $d_0$, in we define an operation $\mu_m$. This operation ``rotates'' $d_0$ in an anti-clockwise direction inside the two $m+2$-gons in the division of $P$ of which it is a common edge. Thus we have,
$$\mu_m(d_i)=d_{i+1}\quad, \textnormal{for all $0\leqslant i\leqslant m-1$}$$
and
$$\mu_m(d_m)=d_{0}$$
where $d_1,\ldots,d_m$ are the $m$-allowable diagonals which are the ``rotations'' of $d_0$. Every $\mu_m^k(d_0)$, $1\leqslant k\leqslant m$, can replace $d_0$ and achieve a distinct division of $P$ as described above.
Further we define, $\mu_m^{-1}=\mu_m^m$. We call $\mu_{m}$ and $\mu_m^{-1}$ \textit{elementary polygonal moves}.
\end{defn}
For the most part we will be interested in applying $\mu_{m}$ or $\mu_m^{-1}$ which can be thought of as ``rotations''of an $m$-allowable diagonal in the anti-clockwise and clockwise directions respectively. Given a division of $P$ we can apply the elementary polygonal moves to the $m$-allowable diagonals in the division to create distinct divisions. As in section \ref{sec:m_clust} an $m$-cluster-tilted algebra can be associated with the original division and with the new division and we say that we have \textit{mutated} the algebra corresponding to the original division to the algebra corresponding to the new division.
\begin{rem} $\mu_1$ can be thought of as the quiver mutation operation in a 1-cluster-tilted algebra of type $A_n$ via the approach of \cite{CCS}. Since for 1-cluster-tilted algebras of Dynkin type the quiver determines the relations \cite{BMR3}, it is possible to define $\mu_1$ without reference to relations for 1-cluster-tilted algebras. For $m$-cluster-tilted algebras of Dynkin type this is no longer possible since the quiver does not determine the relations (as we have seen in example \ref{relex}).
\end{rem}
\begin{ex}\label{mu3ex}We illustrate the algebra mutation induced by the elementary polygonal moves with an example. The algebra mutation shown is achieved by applying $\mu_3^{-1}$ at the vertex circled in the first quiver.
\begin{center}
\includegraphics[scale=0.5]{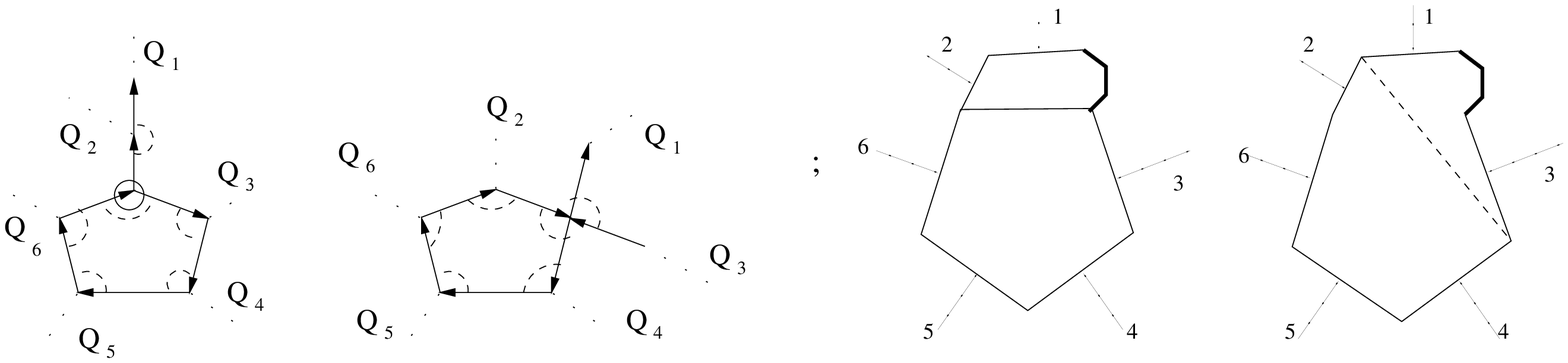}
\end{center}
In the above figure the numbers 1-6 denote where the division continues arbitrarily, the labels $Q_i$ denote the corresponding parts of the quiver with relations. The edges in the polygonal configurations shown in bold typeface are edges of $P$.
\end{ex}
We now wish to show that using the elementary polygonal moves we can induce algebra mutations such that the original and new $m$-cluster-tilted algebras are derived equivalent. The following theorem, due to Bessenrodt and Holm \cite{BH} plays a crucial role in determining when this is possible.
\begin{thm}\cite[thm 4.1,cor 4.3]{BH}\label{BH} Let $\Lambda=kQ/I$ be a gentle algebra with $ec(\Lambda)$ and $oc(\Lambda)$ being the number of cycles of even length with full relations and the number of cycles of odd length with full relations respectively. Then the Cartan matrix, $C_{\Lambda}$, is unimodular equivalent to a diagonal matrix of the form:
$$
\left(
\begin{array}{ccccccccc}
2 & {} & {} & {} & {} & {} & {} & {} & {}\\
{} & \ddots & {} & {} & {} & {} & {} & {} & {}\\
{} & {} & 2 & {} & {} & {} & {} & {} & {}\\
{} & {} & {} & 0 & {} & {} & {} & {} & {}\\
{} & {} & {} & {} & \ddots & {} & {} & {} & {}\\
{} & {} & {} & {} & {} & 0 & {} & {} & {}\\
{} & {} & {} & {} & {} & {} & 1 & {} & {}\\
{} & {} & {} & {} & {} & {} & {} & \ddots & {}\\
{} & {} & {} & {} & {} & {} & {} & {} & 1\\
\end{array}
\right)
 $$
where there are $oc(\Lambda)$ 2's and $ec(\Lambda)$ 0's on the diagonal
\end{thm}
Recall that the Cartan matrix for a finite dimensional algebra $\Lambda \cong kQ/I$ is given by $C_\Lambda=[c_{i,j}]_{1\leqslant i,j,\leqslant |Q_0|}$ where $c_{i,j}=dim_k \Hom_\Lambda(P_i,P_j)=dim_k\Hom_\Lambda(e_i\Lambda,e_j\Lambda)=dim_k\,\,e_j\Lambda e_i$. Here $e_i$ is a principle idempotent and $e_i\Lambda\cong P_i$ are the indecomposable projectives.
\\With the above conventions the columns of the Cartan matrix are the dimension vectors of the indecomposable projective right $\Lambda$ modules. For an $m$-cluster tilted algebra of type $A_n$ we have the following basic properties of the Cartan matrix which follow from the combinatorial descriptions in section \ref{sec:m_clust}.
\begin{itemize}
\item The entries of $C_\Lambda$ are either 0 or 1.
\item There are at most $n$ entries which are 1 in any given row or column. In the case where there are $n$ entries which are 1 in some column then the $m$-cluster tilted algebra in question is isomorphic to the path algebra of the linear orientation of the Dynkin diagram $A_n$.
\item For a vertex $i$ in a cycle there are at least two entries which are 1 in the $i$-th row.
\end{itemize}
It is well known that the unimodular equivalence class of the Cartan matrix of a finite dimensional algebra is a derived invariant. Therefore since the $m$-cluster-tilted algebras of type $A_n$ are gentle we have, by theorem \ref{BH}, that any two $m$-cluster-tilted algebras which are derived equivalent have the same number of $m+2$-cycles. Hence in order to mutate an $m$-cluster-tilted algebra and induce a derived equivalence we must ensure that we preserve the number of $m+2$-cycles.
\begin{rem}It is not always possible to induce a derived equivalence by applying the elementary polygonal moves, for example applying $\mu_3$ at the circled vertex in figure in example \ref{mu3ex} would not preserve the number of cycles and so the mutated algebra would not be derived equivalent to the original algebra.
\end{rem}
It is our aim to show that each connected component of an $m$-cluster-tilted algebra with a given number, $r\in\mathbb{N}_0$ of $m+2$-cycles is derived equivalent to the following \textit{normal form} having the same number of $m+2$-cycles.
\begin{defn}\label{nformdef}
We define the normal form to be the following finite dimensional algebra given as a quiver with relations:
\begin{center}
\includegraphics[scale=0.5]{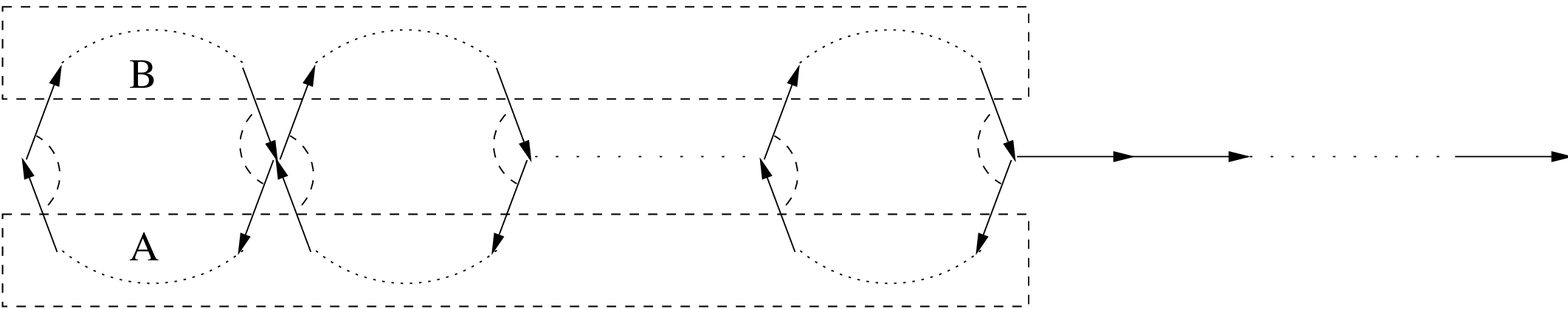}
\end{center}
In the normal form shown the vertices in the region labeled A are called \textit{vertices of type A} while those in the region labeled B are called \textit{vertices of type B}. The vertices which connect cycles to each other or to the linear part of the quiver are not included in either region A or region B (the same holds for the extreme left vertex).
\\The cycles are $m+2$-cycles and have full relations (see section \ref{sec:m_clust}). The number of $m+2$-cycles is $r$, for some $r\geqslant 0$. The total number of vertices is $n$.
\\If $m$ is even then there are $\frac{m}{2}$ vertices of types A and B per cycle. If $m$ is odd there are $\left(\frac{m-1}{2}\right)$ vertices of type B and $\left(\frac{m+1}{2}\right)$ vertices of type A per $m+2$-cycle.
\end{defn}
We remark that there are other possible choices for the normal form in which the number of vertices of types A and B would be different.
\\
\\[10pt]To prove our claim that all connected components of $m$-cluster-tilted algebras are derived equivalent to the quiver with relations specified by the normal form we employ mutations induced by the elementary polygonal moves which result in derived equivalent original and mutated algebras.
\\Shortly we give a list containing local descriptions of some algebra mutations induced by the elementary polygonal moves which preserve the number of $m+2$-cycles occurring in the quiver of an $m$-cluster-tilted algebra. We will use these extensively in section \ref{sec:algo}. Before presenting the list some explanation of the information contained in the figures below is required.
\\Denote by $\Lambda$ the original $m$-cluster-tilted algebra and by $\Lambda^\prime$ the mutated algebra. Each figure contains four diagrams. Two of these are parts of quivers with relations, the first of which is some local situation in $\Lambda$. The circled vertex in this first diagram denotes the vertex of the quiver at which the algebra mutation corresponding to the elementary move takes place. The second of the quiver diagrams shows the resultant local situation in $\Lambda^\prime$. The curved dotted line appearing in both quiver diagrams indicates an $m+2$-cycle.
\\The third and fourth diagrams are the local polygonal configurations corresponding to $\Lambda$ and $\Lambda^\prime$ respectively. The dashed line in the fourth diagram illustrates the new position of the $m$-allowable diagonal corresponding to the mutation vertex after the application of an elementary polygonal move.
\\The edges which appear in bold typeface correspond to boundary edges of $P$. It should be noted that one edge in bold typeface may correspond to as many as $m+1$ edges of $P$. Each figure is general in the sense that it is independent of the values of $m$ and $n$, though of course we assume $n$ is large enough to allow the configurations we show to exist.
\\In some figures we may make further assumptions on the quiver of $\Lambda$, if this is the case then these assumptions will be stated in a short remark.
\\The $(D)$ which sometimes appears below the numbering denotes a mutation which has arisen from $\mu_m^{-1}$ at the mutation vertex.
\begin{center}
\includegraphics[scale=0.5]{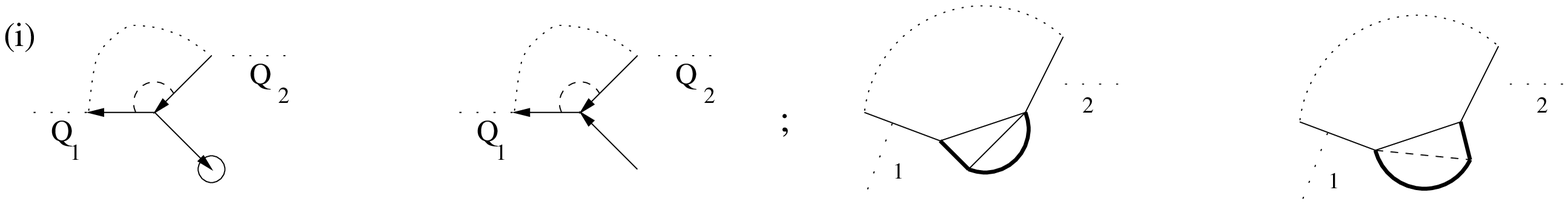}
\end{center}
\begin{center}
\includegraphics[scale=0.5]{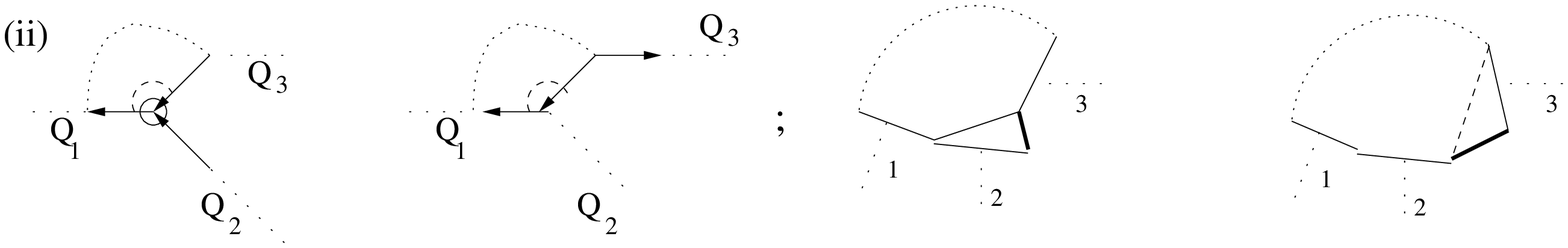}
\end{center}
There is no relation involving an arrow in $Q_2$ and the arrow going into the mutation vertex (and not in the cycle) which is shown.
\begin{center}
\includegraphics[scale=0.5]{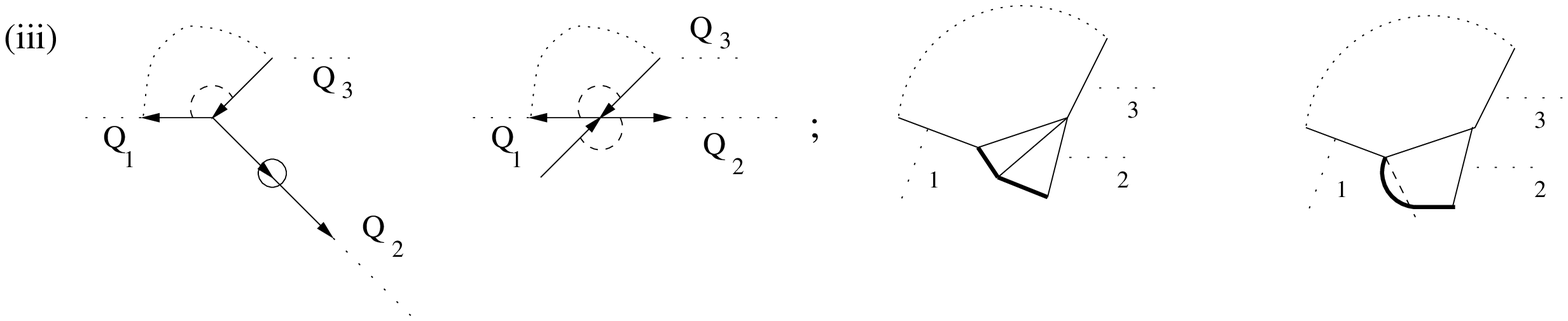}
\end{center}
\begin{center}
\includegraphics[scale=0.5]{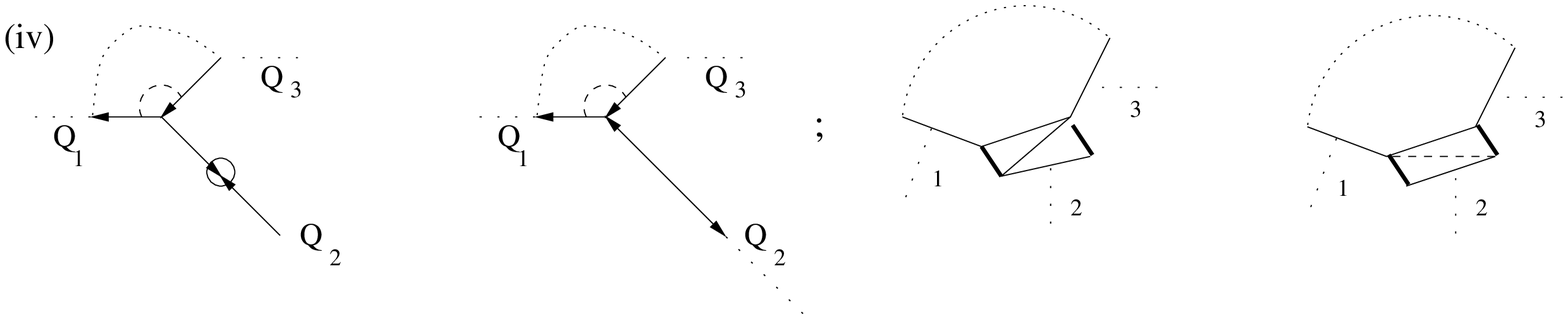}
\end{center}
\begin{center}
\includegraphics[scale=0.5]{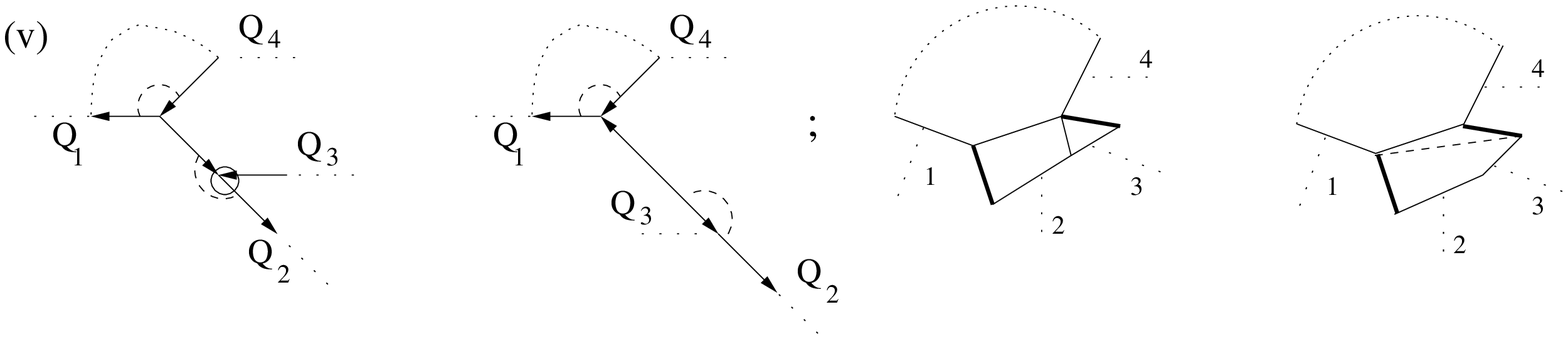}
\end{center}
\begin{center}
\includegraphics[scale=0.5]{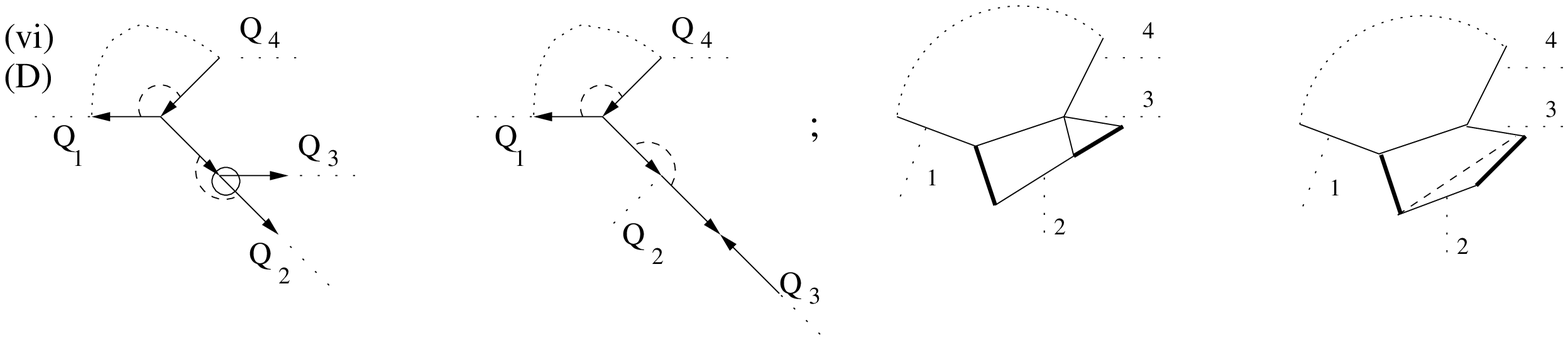}
\end{center}
In the two figures immediately above we assume that there is no relation involving an arrow in $Q_3$ and the arrow into, or our of, the mutation vertex.
\begin{center}
\includegraphics[scale=0.5]{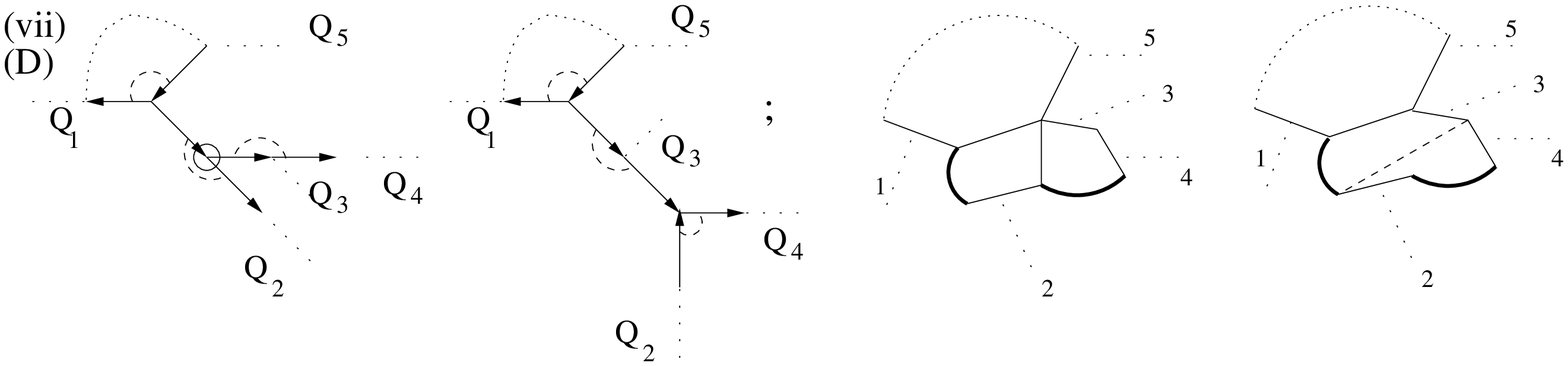}
\end{center}
\begin{center}
\includegraphics[scale=0.5]{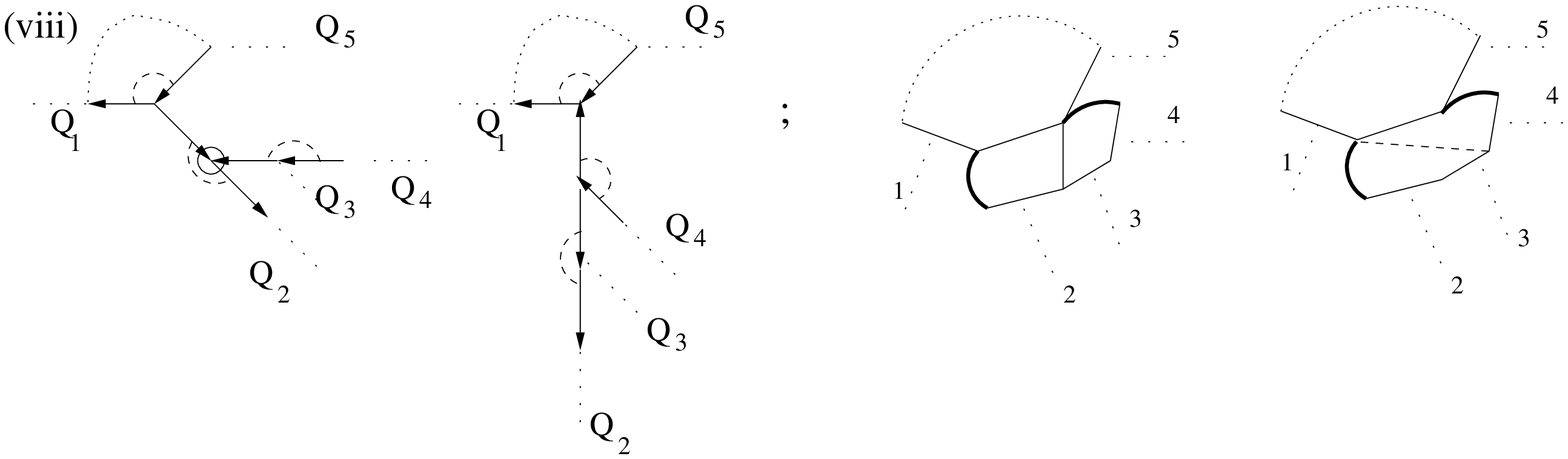}
\end{center}
\begin{center}
\includegraphics[scale=0.5]{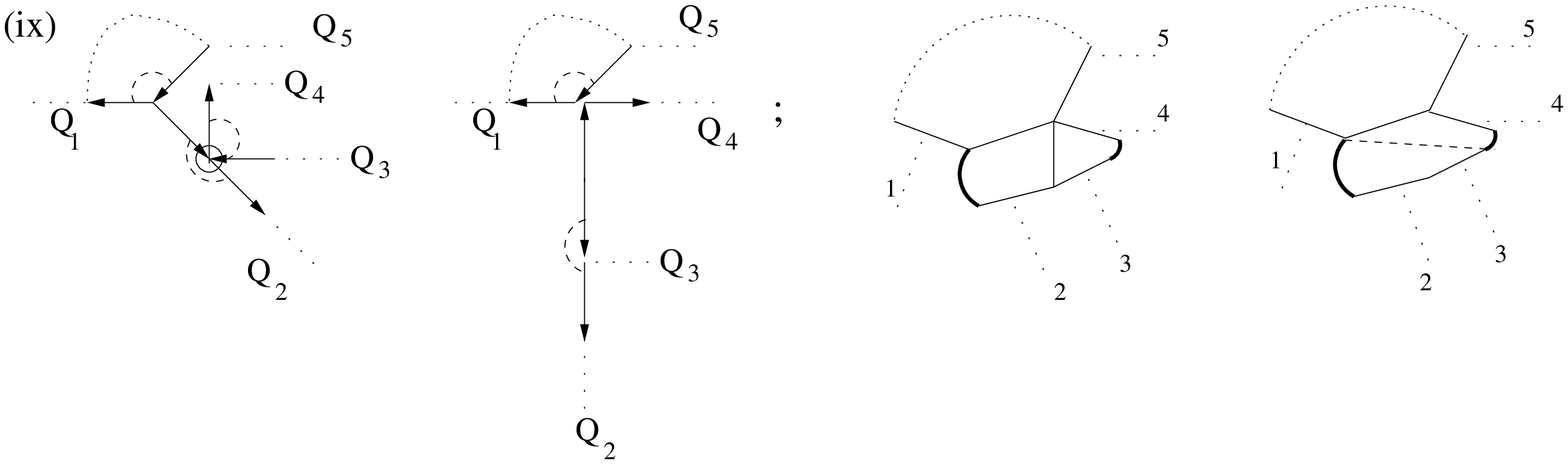}
\end{center}
\begin{center}
\includegraphics[scale=0.5]{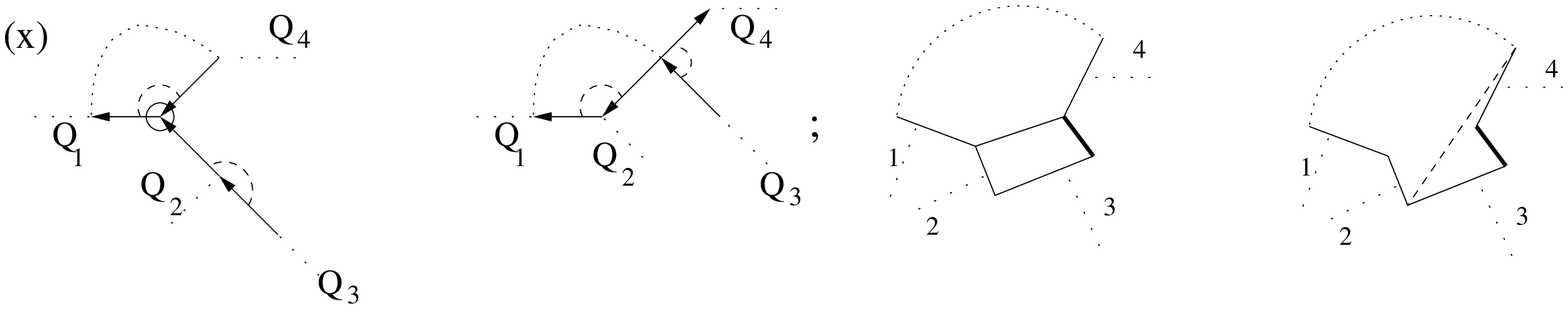}
\end{center}
\begin{center}
\includegraphics[scale=0.5]{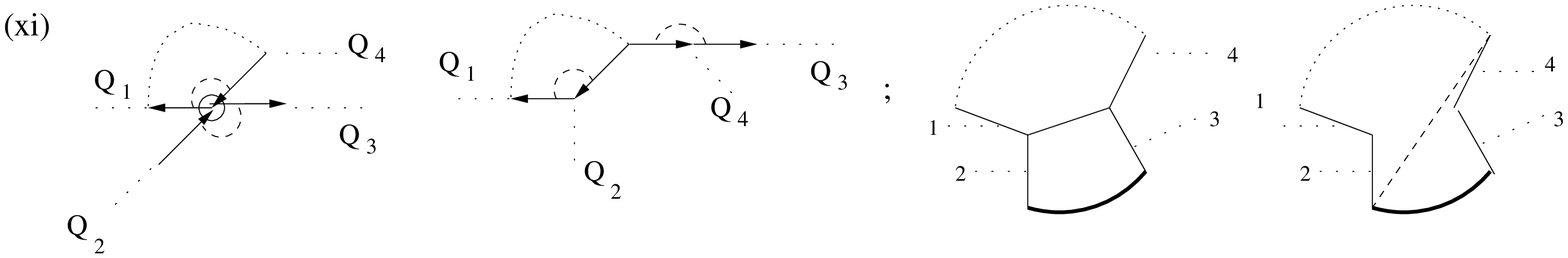}
\end{center}
Here we assume there is no relation between an arrow from $Q_2$ and the arrow going into the mutation vertex (which is not in the cycle) and that there is no relation between the arrow going out of the mutation vertex (which is not in the cycle) and any arrow in $Q_3$.
\begin{center}
\includegraphics[scale=0.5]{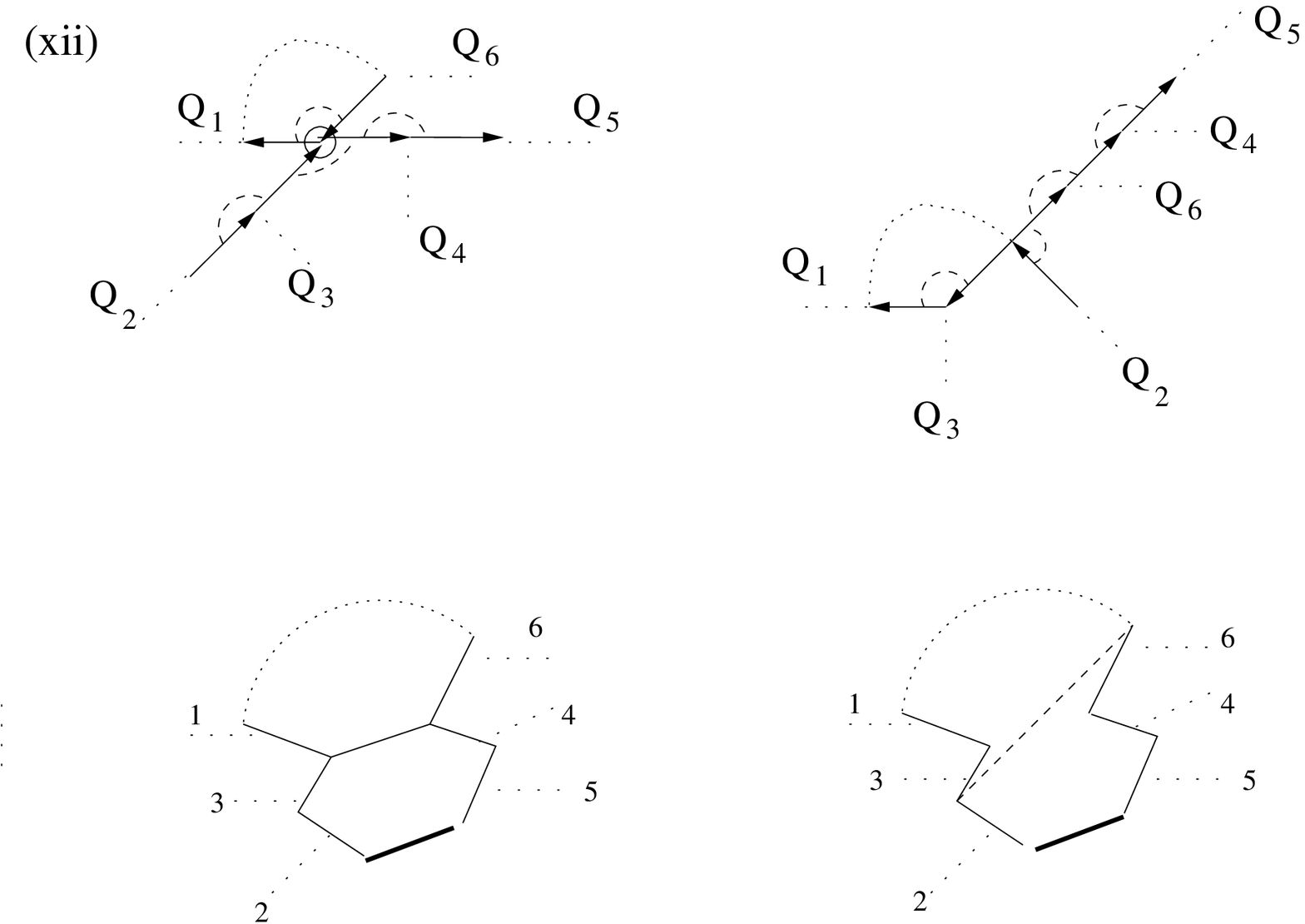}
\end{center}
\begin{center}
\includegraphics[scale=0.5]{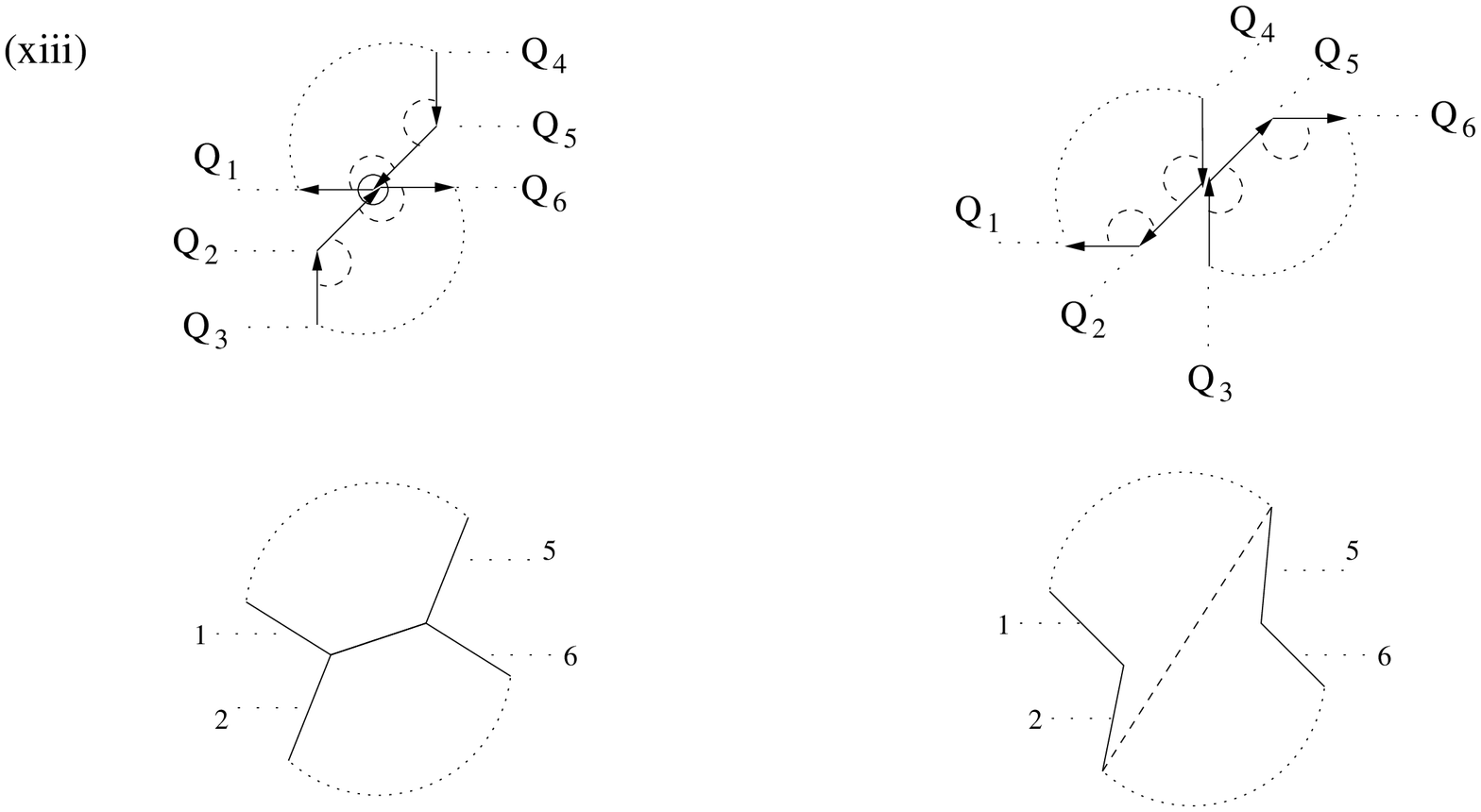}
\end{center}
\begin{center}
\includegraphics[scale=0.5]{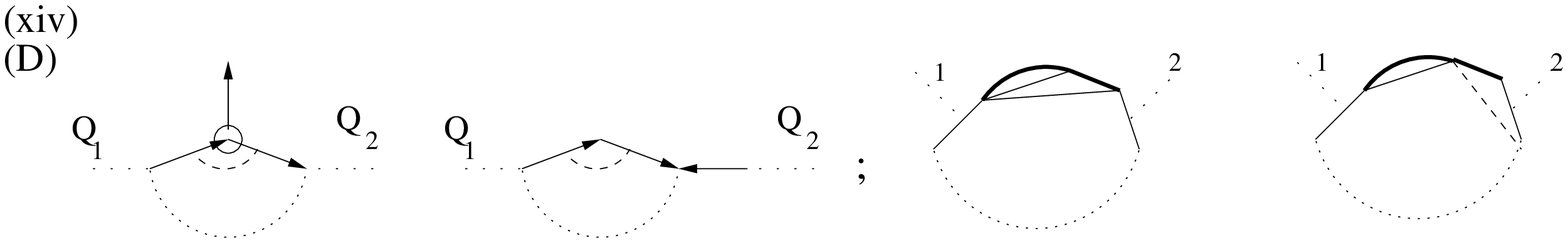}
\end{center}
\begin{center}
\includegraphics[scale=0.5]{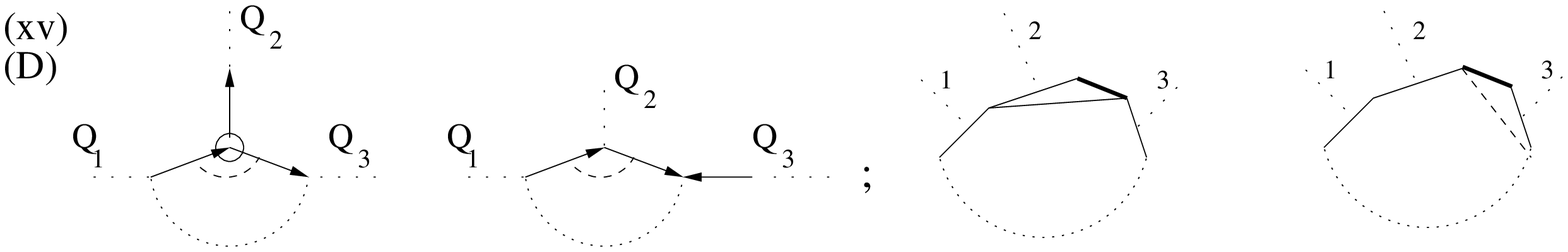}
\end{center}
There is no relation involving an arrow in $Q_2$ and the arrow going out of the mutation vertex (and not in the cycle) which is shown.
\begin{center}
\includegraphics[scale=0.5]{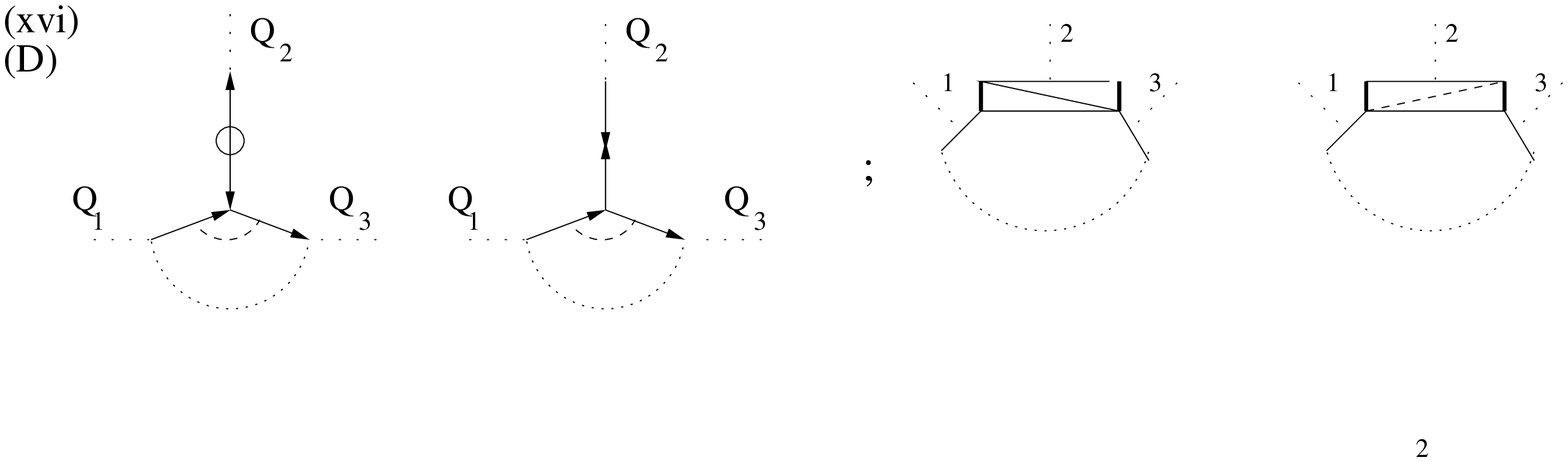}
\end{center}
\begin{center}
\includegraphics[scale=0.5]{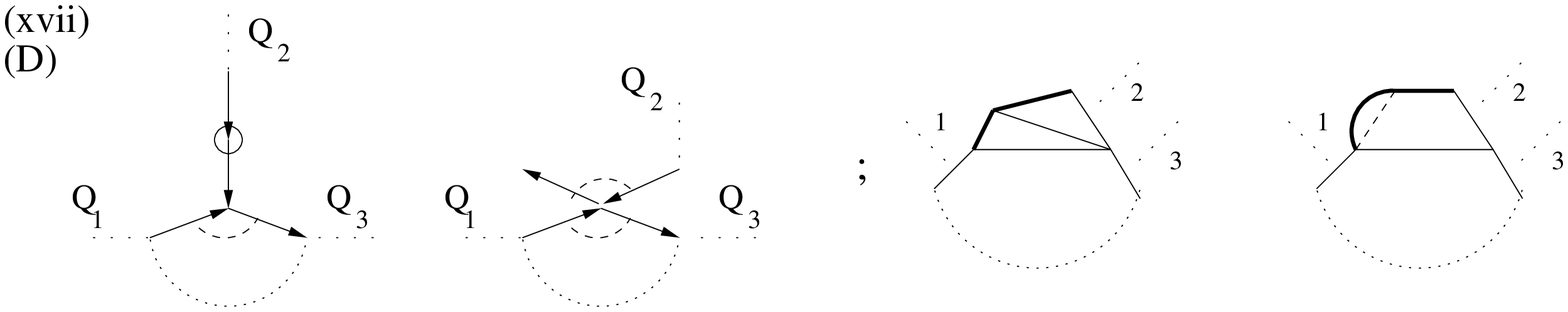}
\end{center}
\begin{center}
\includegraphics[scale=0.5]{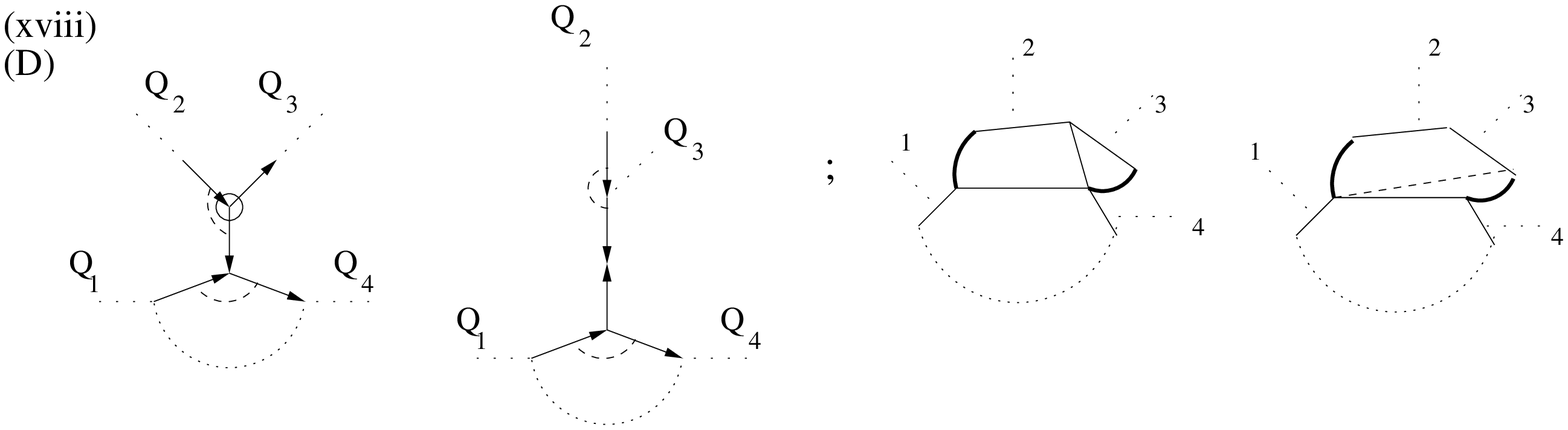}
\end{center}
\begin{center}
\includegraphics[scale=0.5]{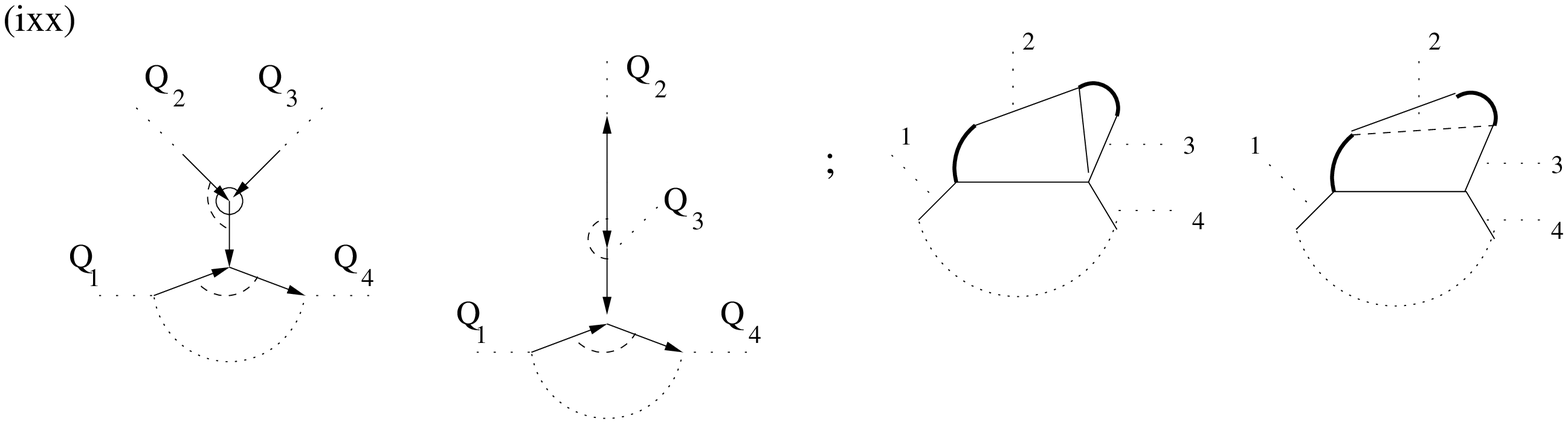}
\end{center}
In the two figures immediately above we assume that there is no relation involving an arrow in $Q_3$ and the arrow into, or our of, the mutation vertex.
\begin{center}
\includegraphics[scale=0.5]{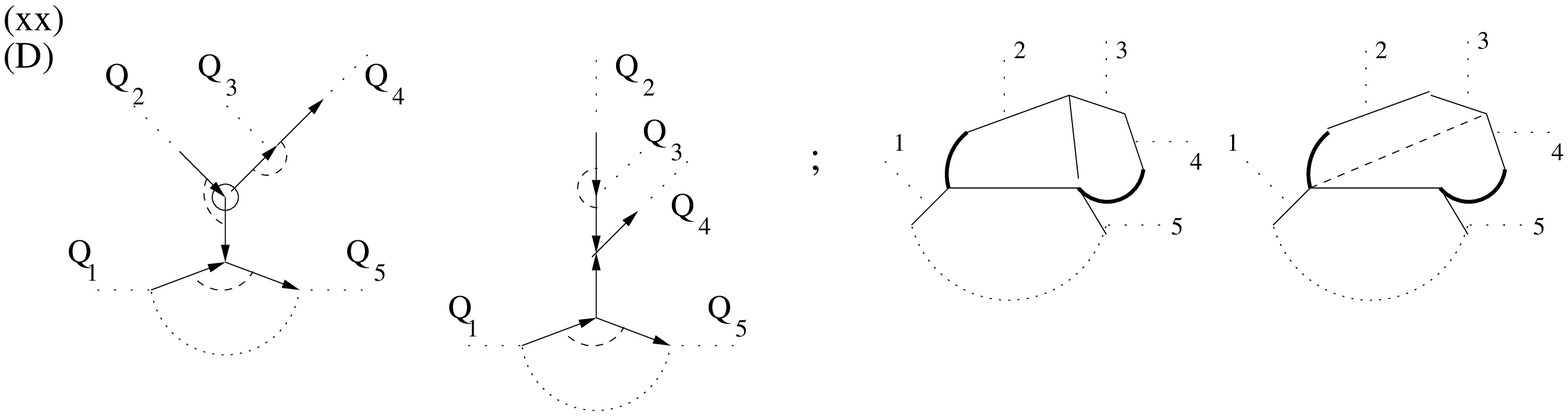}
\end{center}
\begin{center}
\includegraphics[scale=0.5]{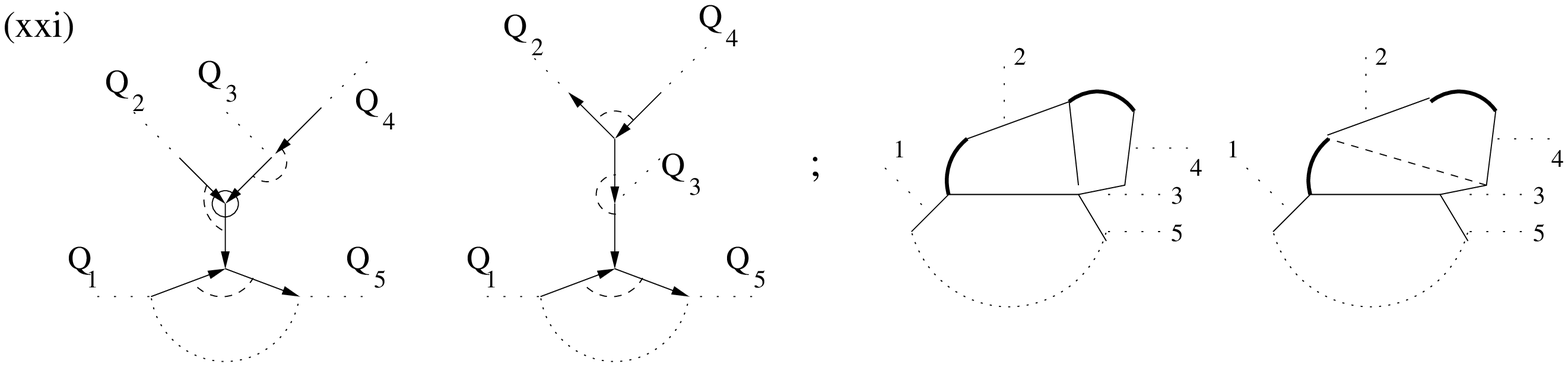}
\end{center}
\begin{center}
\includegraphics[scale=0.5]{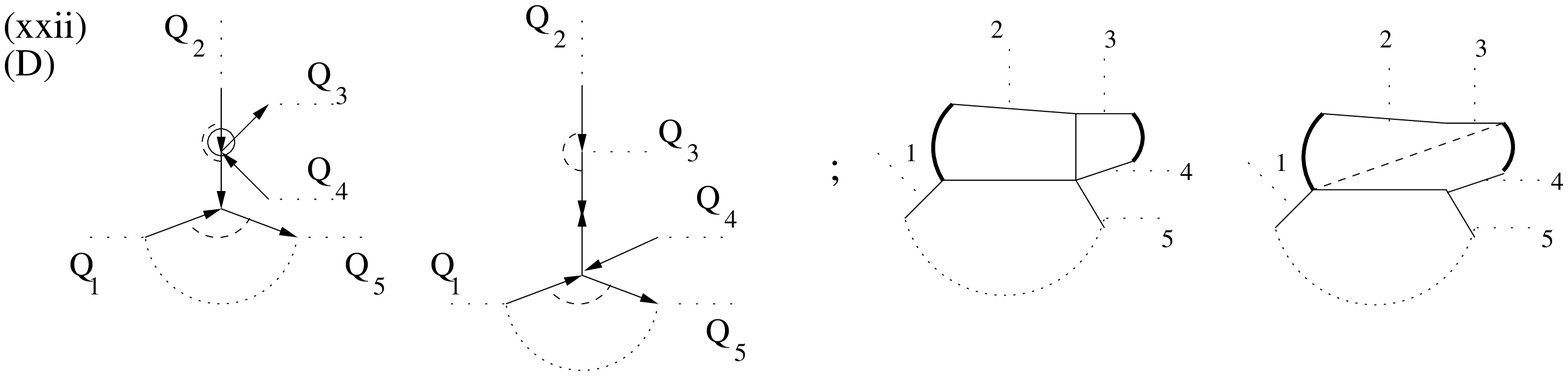}
\end{center}
\begin{center}
\includegraphics[scale=0.5]{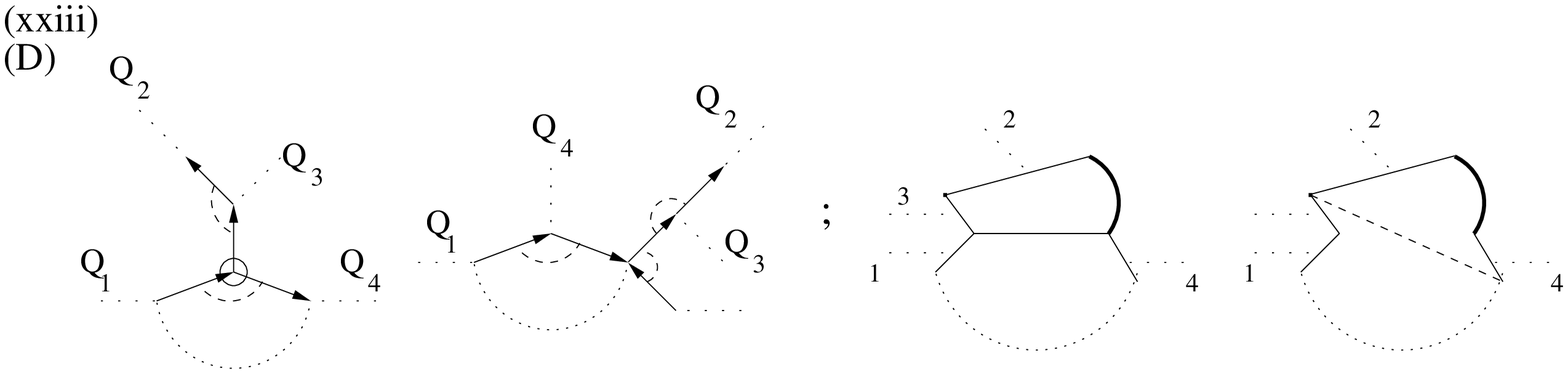}
\end{center}
\begin{center}
\includegraphics[scale=0.5]{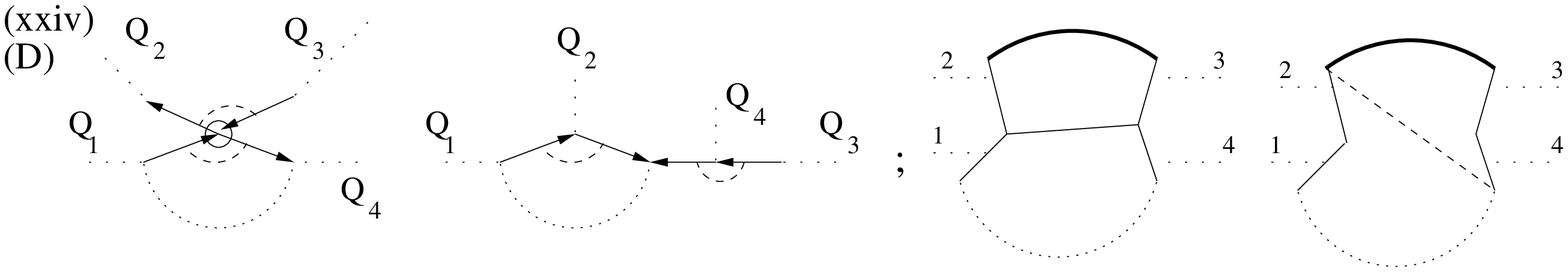}
\end{center}
Here we assume there is no relation between an arrow from $Q_2$ and the arrow going into the mutation vertex (which is not in the cycle) and that there is no relation between the arrow going out of the mutation vertex (which is not in the cycle) and any arrow in $Q_3$.
\begin{center}
\includegraphics[scale=0.5]{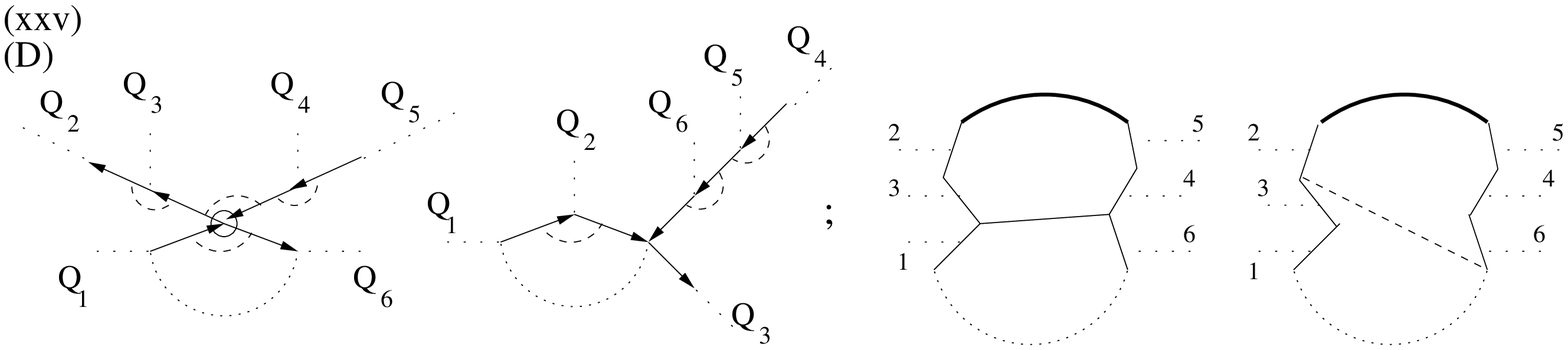}
\end{center}
\begin{center}
\includegraphics[scale=0.5]{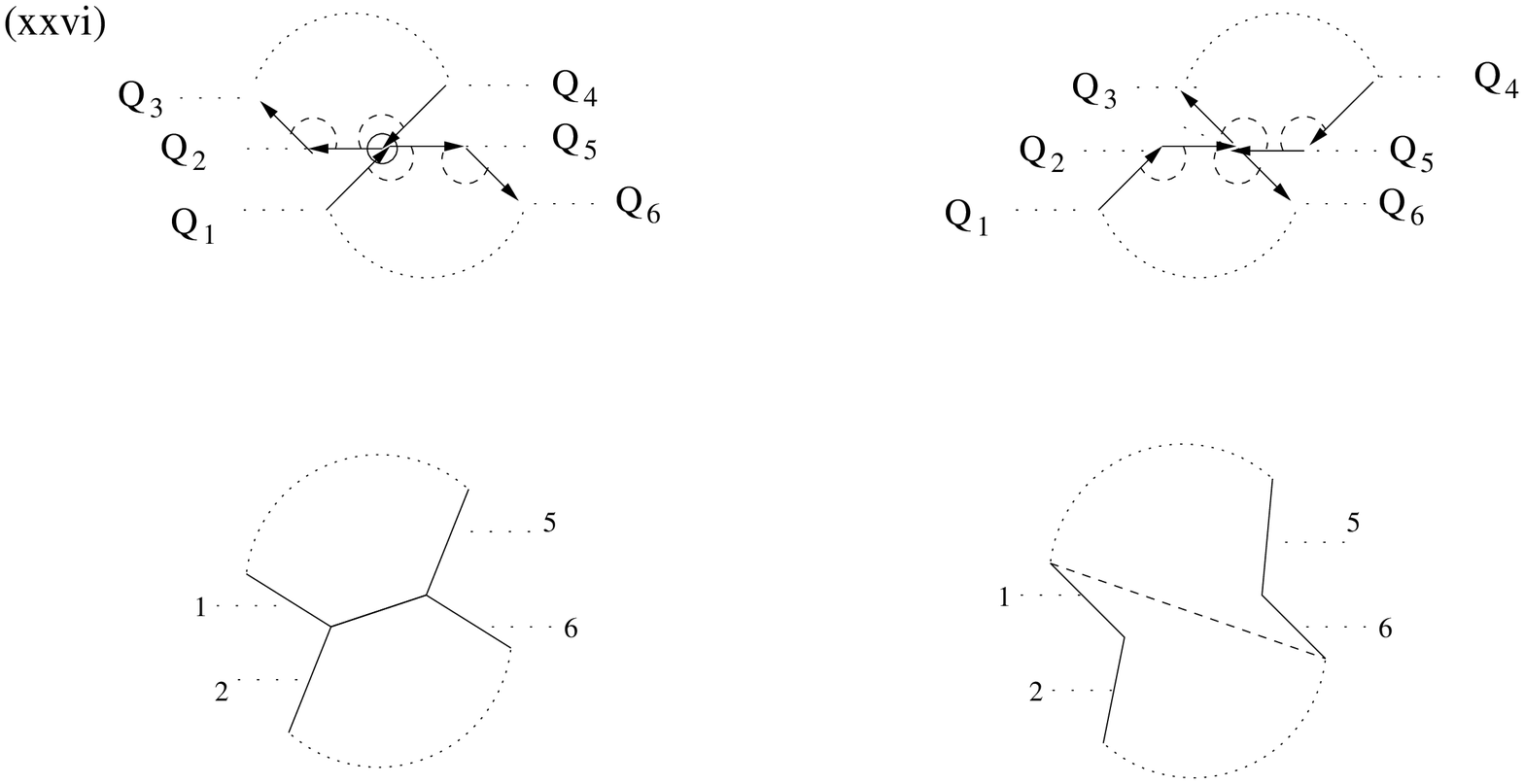}
\end{center}
There are two groups of diagrams. The first group are those numbered (i) to (xiii). The second consists of the remaining figures (xiv) to (xxvi). The two groups will be used at different stages in the algorithm described in section \ref{sec:algo} which reduces any connected component of an $m$-cluster-tilted algebra of type $A_n$ to our chosen normal form in definition \ref{nformdef}.
\begin{rem}Notice that in the above list some of the figures may be viewed as simplified versions of another figure. For example figure (i) is a simplified version of figure (iii). To realize figure (i) from figure (iii) we delete the parts of the quiver which are not needed, namely $Q_2$ and the arrow into $Q_2$. In more complex examples such simplifications may not be immediately apparent. We include the longer list to facilitate checking the algorithm in section \ref{sec:algo}.
\end{rem}
In order to demonstrate that in each of the local algebra mutations above the original and mutated algebras are derived equivalent we now provide explicit tilting complexes which achieve this. There are two such complexes, one corresponding to each of the elementary polygonal moves in definition \ref{elmove}.
\\Before stating the next theorem we must define the notations used. Let $\Lambda\cong kQ_{\mathcal{T}}/\mathcal{I}_{\mathcal{T}}$ be an $m$-cluster-tilted algebra and let $P_i$ denote the finitely generated projective right $\Lambda$-module, $e_i\Lambda$ where $e_i$ is the primitive idempotent corresponding to $i\in (Q_{\mathcal{T}})_0$ (the set of vertices of $Q_{\mathcal{T}}$).
\\Suppose we mutate $\Lambda$ to $\Lambda^\prime$ via an elementary polygonal move, then we distinguish certain vertices in $({Q_{\mathcal{T}}})_0$. We denote by $mut$ the mutation vertex. If there exist arrows in $Q_{\mathcal{T}}$ whose target is $mut$ then we label the source $in_t$ and the arrow $i_t$, where $1\leqslant t\leqslant 2$ (since $\Lambda$ is gentle there can be at most two ingoing and two outgoing arrows). We denote by $out_u$, $1\leqslant u\leqslant 2$, the targets of any arrows out of $mut$ and label such arrows $o_u$. Since the $m$-cluster-tilted algebras in this case are gentle we may adopt the convention that should the paths $o_2i_1$ and $o_1i_2$ exist they are zero.
\\Further, should there exist an arrow $\gamma$ such that $i_t\gamma=0$ then we label the source of $\gamma$ by $pre_t$. Similarly, if there exists an arrow $\delta$ such that $\delta o_u=0$ then we label the target of $\delta$ by $post_u$.
\\Next, $P(\Lambda)$ denotes the full subcategory of finitely generated projective right $\Lambda$-modules; $K^b(P(\Lambda))$ denotes the homotopy category of bounded complexes of finitely generated projective right $\Lambda$-modules; and $\mathcal{D}^b(\Lambda)$ is the bounded derived category of finitely generated $\Lambda$-modules. Our last notational convention is that given an arrow $\alpha:i\rightarrow j$ between any two vertices $i,j\in(Q_{\mc{T}})_0$, the map $P_i\rightarrow P_j$ between the projective indecomposable $\Lambda$-modules induced by left multiplication by $\alpha$ will also be denoted $\alpha$.
\\[10pt]
\\Part of the proof of theorem \ref{tilting_reg} (and part of the proof of theorem \ref{tilting_dual}) involves determining the existence of maps between certain indecomposable summands of a tilting complex. In order to do this we use Happel's alternating sum formula \cite[chapter III, 1.3, 1.4]{Ha}.

\begin{thm}\label{altsum}\cite{Ha} For a finite dimensional $k$-algebra $\Lambda$ let $C=(C^r)_{r\in\mathbb{Z}}$ and $D=(D^s)_{s\in\mathbb{Z}}$ be complexes in $K^b(P(\Lambda))$. Then,
$$\sum_{i}(-1)^i\,dim_k \Hom_{K^b(P(\Lambda))}(C,D[i])=\sum_{r,s}(-1)^{r-s}\,dim_k \Hom_{\Lambda}(C^r,D^s).
$$.
In particular if $C$ and $D$ are direct summands of a tilting complex then,
$$dim_k \Hom_{K^b(P(\Lambda))}(C,D)=\sum_{r,s}(-1)^{r-s}\,dim_k \Hom_{\Lambda}(C^r,D^s).
$$
\end{thm}
The next two results prove that when the elementary polygonal move $\mu_m$ defined in definition \ref{elmove} is applied in the figures labeled (i)-(v),(viii)-(xiii),(ixx) and (xxi) above there is a derived equivalence between the initial algebra $\Lambda$ and the new algebra $\Lambda^\prime$.

\begin{thm}\label{tilting_reg}Suppose $\Lambda$ is a connected $m$-cluster-tilted algebra of type $A_n$ in which we can perform $\mu_m$ at a given vertex, $mut$, and preserve the number of $m+2$-cycles and connectedness. Then then following complex, $T$, in $K^b(P(\Lambda))$, where $P_{mut}$ is in degree 1, is a tilting complex.
$$T:\cdots\rightarrow 0\rightarrow P_{in_1}\oplus P_{in_2}\oplus\left(\bigoplus_{i\neq mut}P_i\right)\stackrel{[i_1,i_2,0]}{\longrightarrow} P_{mut}\rightarrow 0\rightarrow\cdots$$
Moreover, the endomorphism algebra, $\End_{\mathcal{D}^b(\Lambda)}(T)\cong kQ^\prime/I^\prime=\Lambda^\prime$, of $T$ can be obtained from $\Lambda$ using the following algorithm:
\begin{enumerate}
\item reverse all arrows going into the mutation vertex.
\item if there exist arrows $in_t\rightarrow mut\rightarrow out_t$ (which is non-zero by definition), there must exist an arrow $in_t\rightarrow out_t$ in $Q^\prime$ and the arrow $mut\rightarrow out_t$ does not exist in $Q^\prime$, $t=1,2$.
\item if there exist arrows $pre_t\rightarrow in_t\rightarrow mut$ such that the composition is zero, the arrow $pre_t\rightarrow in_t$ factors over $mut\rightarrow in_t$ in $Q^\prime$, $t=1,2$.
\item relations in $\Lambda^\prime$ around the mutated vertex are described by the following diagrams. If arrows and relations exist around the mutation vertex in $\Lambda$ as shown in the first diagram then the second diagram shows the new arrows and relations between the corresponding vertices of $\Lambda^\prime$.
\begin{center}
\includegraphics[scale=.5]{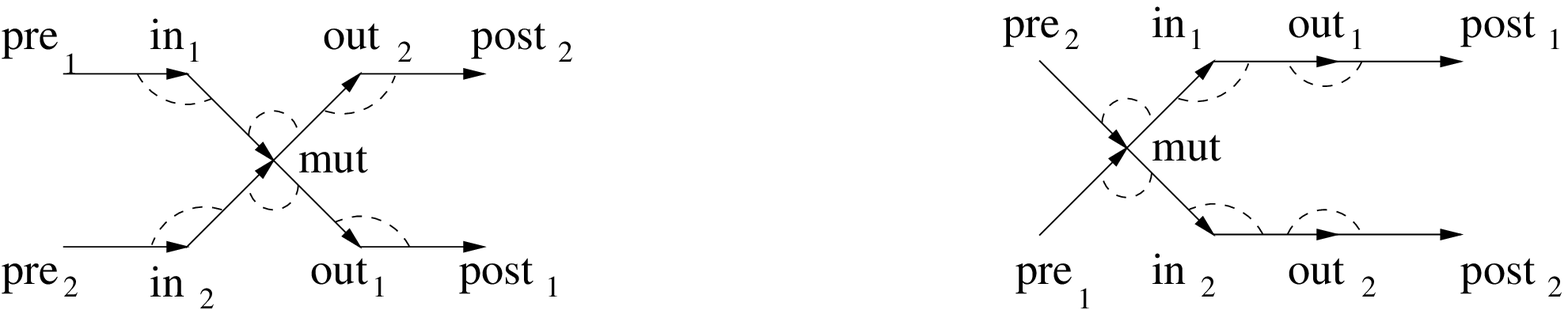}
\end{center}
\end{enumerate}
\end{thm}

\begin{rem} In the last part of the algorithm to determine the endomorphism algebra of the tilting complex in theorem \ref{tilting_reg} it should be noted that we do not suggest that the local configurations shown must exist at the mutation vertex, rather we show the arrows and relation which, \textit{if} they exist, are altered by the mutation. The second diagram describes the new configuration of arrows and relations. It is clear that the mutations in figures (i)-(v),(viii)-(xiii),(ixx) and (xxi) satisfy the conditions of theorem \ref{tilting_reg}.
\end{rem}

\begin{proof}First we prove that $T$ is indeed a tilting complex. In other words we prove that $\Hom_{D^b(\Lambda)}(T,T[k])=0$, for all $k\in\mathbb{Z}, k\neq 0$ (recall $\Hom_{K^b(P(\Lambda))}(T,T[k])=0$ implies $\Hom_{D^b(\Lambda)}(T,T[k])=0$, for all $k\in\mathbb{Z}$ \cite[2.5]{KZ}) and that the subcategory $\add (T)$ generates $K^b(P(\Lambda))$ as a triangulated category (see \cite[3.2]{KZ} for definition of tilting complex).
\\[10pt]
\\The indecomposable direct summands of $T$ are:
$$T_{mut}:\,\,\cdots0\rightarrow P_{in_1}\oplus P_{in_2}\stackrel{[i_1,i_2]}{\rightarrow} P_{mut}\rightarrow 0\rightarrow\cdots$$
$$T_{i}:\,\,\cdots0\rightarrow P_{i}\rightarrow 0\rightarrow\cdots\quad i\neq mut$$
where the $P_i$'s are in degree zero for $i\neq mut$ and $P_{mut}$ is in degree 1.
\\Obviously, $\Hom_{K^b(P(\Lambda))}(T_i,T_j[k])=0$ for some $i,j\neq mut$ and $k\in\mathbb{Z}\, ,k\neq 0$. It is also clear that $\Hom_{K^b(P(\Lambda))}(T,T[k])=0$ for all $k\in\mathbb{Z},\,|k|\geqslant 2$.
\\Note that for $m$-cluster categories of Dynkin type the dimension of the non-zero $\Hom$ spaces in the $m$-cluster category is at most 1 by remark \ref{dimrem}. When proving that certain $\Hom$ spaces are zero in $K^b(P(\Lambda))$ we consider basis elements and show they are homotopic to zero. For two vertices $i$ and $j$ of $Q$ such that there is a non-zero path from $i$ to $j$ in $Q$ we always choose the map $P_i\rightarrow P_j$ induced by the path as our basis of $\Hom_{\mc{C}_m}(P_i,P_j)$.
\\Consider $\Hom_{K^b(P(\Lambda))}(T_{mut},T_i[-1])$ for all $i\neq mut$ and let $\gamma:P_{mut}\rightarrow P_i$ be a basis. Since $\gamma$ must start with either $P_{mut}\stackrel{o_1}{\rightarrow}P_{out_1}$ or $P_{mut}\stackrel{o_2}{\rightarrow}P_{out_2}$ we have that the diagram,
$$
\xymatrix{
\cdots 0\ar[r] & P_{in_1}\oplus P_{in_2}\ar[d] \ar[r]^-{[i_1,i_2]} & P_{mut}\ar[r]\ar[d]^{\gamma}\ar[r] &0\cdots\\
\cdots 0\ar[r] & 0\ar[r] & P_{i}\ar[r] & 0\cdots
}
$$
commutes if and only if $\gamma$ is the zero map. Therefore, $\Hom_{K^b(\Lambda)}(T_{mut},T_i[-1])=0$ for all $i\neq mut$.
\\Next consider $\Hom_{K^b(\Lambda)}(T_i,T_{mut}[1])$ for all $i\neq mut$ and let $\delta:P_i\rightarrow P_{mut}$ be a basis. Since $\delta$ must end with either
$P_{in_1}\stackrel{i_1}{\rightarrow}P_{mut}$ or $P_{in_2}\stackrel{i_2}{\rightarrow}P_{mut}$ we have that chain map,

$$
\xymatrix{
\cdots 0\ar[r] & 0\ar[r]\ar[dd] & P_{i}\ar[r]\ar^{\delta}[dd]\ar^{h}[ddl] & 0\cdots\\
&&&\\
\cdots 0\ar[r] & P_{in_1}\oplus P_{in_2} \ar[r]^-{[i_1,i_2]} & P_{mut}\ar[r] &0\cdots
}
$$

is homotopic to zero via the homotopy $h=\nu_s\delta^\prime$ where $i_s\delta^\prime=\delta$ and $\nu_s$ is the inclusion $P_{in_s}\rightarrow P_{in_1}\oplus P_{in_2}$, $s=1,2$. Thus, $\Hom_{K^b(\Lambda)}(T_i,T_{mut}[1])=0$ for all $i\neq mut$.
\\Finally, we examine $\Hom_{K^b(\Lambda)}(T_{mut},T_{mut}[1])$. Since the maps $[i_1,0]$ and $[0,i_2]$ are a basis for $\Hom_\Lambda(P_{in_1}\oplus P_{in_2},P_{mut})$ the homotopies,
$$
\qquad\qquad\qquad\xymatrix{
\cdots 0 \ar[r] & P_{in_1}\oplus P_{in_2}\ar_(.7){\left[ \begin{array}{cc} \SS 1 & \SS 0 \\ \SS 0 & \SS 0\end{array}\right]}[ddl] \ar[dd]^-{[i_1,0]}\ar[r]^-{[i_1,i_2]} & P_{mut}\cdots\ar^{0}[ddl]\\
& & & & & & & & & & &\\
\cdots P_{in_1}\oplus P_{in_2} \ar[r]^(.6){[i_1,i_2]} & P_{mut}\ar[r] &0\cdots\\
}
$$
$$
\qquad\qquad\qquad\xymatrix{
\cdots 0 \ar[r] & P_{in_1}\oplus P_{in_2}\ar_(.7){\left[ \begin{array}{cc} \SS 0 & \SS 0 \\ \SS 0 & \SS 1\end{array}\right]}[ddl] \ar[dd]^-{[0,i_2]}\ar[r]^-{[i_1,i_2]} & P_{mut}\cdots\ar^{0}[ddl]\\
& & & & & & & & & & &\\
\cdots P_{in_1}\oplus P_{in_2} \ar[r]^(.6){[i_1,i_2]} & P_{mut}\ar[r] &0\cdots
}
$$
show that $\Hom_{K^b(\Lambda)}(T_{mut},T_{mut}[1])=0$. Hence, $\Hom_{D^b(\Lambda)}(T,T[k])=0$, for all $k\in\mathbb{Z}$, $k\neq0$.
\\Now we check that $\add(T)$ generates $K^b(P(\Lambda))$ as a triangulated category.
It is clear that the stalk complexes,
$$P^{\cdot}_i[0]:\quad\cdots\rightarrow 0\rightarrow P_i\rightarrow 0\rightarrow\cdots $$
(where $P^{\cdot}_i[k]$ denotes the complex concentrated in degree $k$ with $P_i$ the only non-zero term)
are in $\add(T)$, $i\neq mut$, and therefore $P^{\cdot}_i[k],\, k\in\mathbb{Z}$ are in the triangulated category generated by $\add(T)$. We require that $P^{\cdot}_{mut}[k],\, k\in\mathbb{Z}$ are also in the triangulated category generated by $\add (T)$.
\\Consider the morphism $T_{mut}\stackrel{f}{\rightarrow} T_{in_1}\oplus T_{in_2}$,
$$\begin{array}{c}
\xymatrix{\cdots\ar[r] & 0\ar[r] & P_{in_1}\oplus P_{in_2}\ar^-{\left[ \begin{array}{cc} \SS 1 & \SS 0 \\ \SS 0 & \SS 1\end{array}\right]}[dd]\ar^(.6){[i_1,i_2]}[r] & P_{mut}\ar[r]\ar[dd] & 0\ar[r] &\cdots \\
& & & & &\\
\cdots\ar[r] & 0\ar[r] & P_{in_1}\oplus P_{in_2}\ar[r] & 0\ar[r] & 0\ar[r] &\cdots
}
\end{array}
$$
The mapping cone, $M(f)$, of $f$ (see \cite[section 2.3]{KZ} for definition) is,
$$
\cdots\rightarrow 0\rightarrow P_{in_1}\oplus P_{in_2}\stackrel{\left[\begin{array}{cc} \SS -i_1 &\SS -i_2 \\\SS 1&\SS 0\\\SS 0&\SS 1
\end{array}\right]}{\longrightarrow} P_{mut}\oplus P_{in_1}\oplus P_{in_2}\rightarrow 0\rightarrow\cdots
$$
where $P_{in_1}\oplus P_{in_2}$ is in degree -1. So that,
$$
d_{M(f)}=\left[
\begin{array}{cc}\SS -i_1&\SS -i_2\\\SS 1&\SS 0\\\SS 0&\SS 1
\end{array}\right]
$$
Define the following maps $M(f)\stackrel{\phi}{\rightarrow} P^{\cdot}_{mut}[0]$ and $P^{\cdot}_{mut}[0]\stackrel{\psi}{\rightarrow} M(f)$:
$$
\begin{array}{c}
\xymatrix{\cdots\ar[r]\ar@{}[d]|{\phi} & 0\ar[r] & P_{in_1}\oplus P_{in_2}\ar^-{d_{M(f)}}[r]\ar[d] & P_{mut}\oplus P_{in_1}\oplus P_{in_2}\ar[r]\ar^-{[1,i_1,i_2]}[d] & 0\ar[r] & \cdots\\
\cdots\ar[r] & 0\ar[r] & 0\ar[r] & P_{mut}\ar[r] & 0\ar[r] & \cdots
}
\end{array}
$$
$$
\begin{array}{c}
\xymatrix{
\cdots\ar[r]\ar@{}[dd]|{\psi} & 0\ar[r]& 0\ar[r]\ar[dd] & P_{mut}\ar[r]\ar^{\left[\begin{array}{c}\SS 1\\\SS 0\\\SS 0
\end{array}\right]}[dd] & 0\ar[r] & \cdots\\
& & & & & &\\
\cdots\ar[r] & 0\ar[r] & P_{in_1}\oplus P_{in_2}\ar^-{d_{M(f)}}[r] & P_{mut}\oplus P_{in_1}\oplus P_{in_2}\ar[r] & 0\ar[r] & \cdots
}
\end{array}
$$
Then $\psi\circ\phi$ is:
$$
\begin{array}{c}
\xymatrix{
\cdots\ar[r]\ar@{}[dd]|{\psi\circ\phi} & 0\ar[r] & P_{in_1}\oplus P_{in_2}\ar^-{d_{M(f)}}[r]\ar^{0}[dd] & P_{mut}\oplus P_{in_1}\oplus P_{in_2}\ar[r]\ar^-{\left[\begin{array}{ccc}\SS 1 &\SS i_1 &\SS i_2\\ \SS 0&\SS 0&\SS 0\\  \SS 0&\SS 0&\SS 0
\end{array}\right]
}[dd] & 0\ar[r] & \cdots\\
& & & & &\\
\cdots\ar[r] & 0\ar[r] & P_{in_1}\oplus P_{in_2}\ar^-{d_{M(f)}}[r] & P_{mut}\oplus P_{in_1}\oplus P_{in_2}\ar[r] & 0\ar[r] & \cdots
}
\end{array}
$$
so that $1-\psi\circ\phi$ is homotopic to zero by the following map,
$$
\begin{array}{c}
\xymatrix {
\cdots\ar[r]\ar@{}[dd]|{1-\psi\circ\phi} & 0\ar[r] & P_{in_1}\oplus
P_{in_2}\ar^-{d_{M(f)}}[r]\ar^<<<<<<{\left[\begin{array}{cc}\SS 1 & \SS
      0 \\ \SS 0 & \SS 1
\end{array}\right]
}[dd]\ar_{0}[ddl] & P_{mut}\oplus P_{in_1}\oplus
P_{in_2}\ar[r]\ar^{\left[\begin{array}{ccc}\SS 0 & \SS -i_1 & \SS
      -i_2\\ \SS 0 & \SS 1 & \SS 0\\ \SS 0 & \SS 0 & \SS 1
\end{array}\right]
}[dd]\ar^{h}[ddl] & 0\ar[r] & \cdots\\
& & & & & \\
\cdots\ar[r] & 0\ar[r] & P_{in_1}\oplus P_{in_2}\ar_-{d_{M(f)}}[r] &
P_{mut}\oplus P_{in_1}\oplus P_{in_2}\ar[r] & 0\ar[r] & \cdots
}
\end{array}
$$where $h=\left[\begin{array}{ccc}\SS 0 & \SS 1 & \SS 0\\ \SS 0&\SS 0&\SS 1
\end{array}\right]$.
\\It is clear that $\phi\circ\psi=1_{P^{\cdot}_{mut}[0]}$. Hence, in $K^{b}(P(\Lambda))\,,M(f)\cong P^{\cdot}_{mut}[0]$.
\\It now follows that $\add (T)$ generates  $K^{b}(P(\Lambda))$ as a triangulated category and so we have proven that $T$ is a tilting complex.
\\[10pt]
\\Since $T$ has finitely many indecomposable summands the endomorphism algebra, $\End_{\mathcal{D}^b(\Lambda)}(T)$, of $T$ can be expressed as a finite quiver with relations denoted by $\Lambda^\prime=kQ^\prime/I^\prime\cong \End_{\mathcal{D}^b(\Lambda)}(T)$. The vertices of the quiver are in one-to-one correspondence with the indecomposable summands, $T_i$, of $T$. To see that the endomorphism algebra $\End_{\mathcal{D}^b(\Lambda)}(T)$ is described by the algorithm in the statement of theorem \ref{tilting_reg} we examine each step separately. Let $C_{\Lambda}=[c_{ij}]_{1\leqslant i,j\leqslant n}$ denote the Cartan matrix of $\Lambda$ and let $\tilde{C}_{\Lambda}=[\tilde{c}_{ij}]_{1\leqslant i,j\leqslant n}$, where $\tilde{c}_{ij}=dim_k \Hom_{K^b(P(\Lambda))}(T_i,T_j)$, be the Cartan matrix of $\Lambda^\prime=\End_{\mc{D}^b(\Lambda)}(T)$.
\\
\\[10pt]$\bullet$ \textit{reverse all arrows going into the mutation vertex}
\\
\\[10pt]Given an arrow $in_1\rightarrow mut$ there is a morphism of complexes $T_{mut}\rightarrow T_{in_1}$ given by:
$$
\begin{array}{c}
\xymatrix{
T_{mut}\ar[d]&\cdots\ar[r] & 0\ar[r] & P_{in_1}\oplus P_{in_2}\ar^-{[i_1,i_2]}[r]\ar^{\pi_1}[d]& P_{mut}\ar[r]\ar[d] & 0\ar[r] & \cdots\\
T_{in_1}&\cdots\ar[r] & 0\ar[r] & P_{in_1}\ar[r] & 0\ar[r] & 0\ar[r] & \cdots
}
\end{array}
$$
where $\pi_1: P_{in_1}\oplus P_{in_2}\rightarrow P_{in_1}$ is the projection map. There is an arrow in the quiver $Q^\prime$ of $\End_{\mathcal{D}^b(\Lambda)}(T)$ corresponding to the morphisms $\pi_1:T_{mut}\rightarrow T_{in_1}$ since it is not possible for this morphism to factor non-trivially through another summand of $T$.
\\Indeed if the morphism $T_{mut}\rightarrow T_{in_1}$ above factors non-trivially through another summand of $T$ then there will exist some vertex $k\neq in_1,mut$ in $Q^\prime$ such that $\tilde{c}_{mut,k}\neq 0$ and $\tilde{c}_{k,in_1}\neq 0$.
\\Theorem \ref{altsum} states that $\tilde{c}_{mut,k}=c_{in_1,k}+c_{in_2,k}-c_{mut,k}$, where the $c_{ij}$'s denote entries of the Cartan matrix of $\Lambda$, and that $\tilde{c}_{k,in_1}=c_{k,in_1}$.
\\Suppose that $\tilde{c}_{k\,in_1}\neq 0$. Now there are two possible cases. If $m=1$ then it is possible that $c_{in_2\,k}=c_{mut\,k}=1$, to see this consider two 3-cycles which share a vertex. Mutating at this shared vertex gives $c_{in_1,out_1}=0$ and $c_{in_2,out_1}=c_{mut,out_1}=1$. If $m\geqslant 2$ then we must have that $c_{in_1,k}=0$ and $c_{in_2\,k}=c_{mut\,k}=0$.  Therefore, $\tilde{c}_{mut\,k}=0$ and so the morphism $\pi_1:T_{mut}\rightarrow T_{in_1}$ corresponds to an arrow in the quiver of $\Lambda^\prime$.
\\A similar argument shows that the morphism $\pi_2:T_{mut}\rightarrow T_{in_2}$ corresponds to an arrow in the quiver $Q^\prime$ of $\End_{\mc{D}^b(\Lambda)}(T)$.
\\
\\[10pt]$\bullet$ \textit{if there exist arrows $in_1\rightarrow mut\rightarrow out_1$, there must exist an arrow $in_1\rightarrow out_1$ in $Q^\prime$.}
\\
\\[10pt]Suppose that there is a non-zero path $in_1\rightarrow mut\rightarrow out_1$ in $\Lambda$. Then we have the following morphism:
$$
\begin{array}{c}
\xymatrix{
T_{in_1}\ar[d]&\cdots\ar[r] & 0\ar[r] & P_{in_1}\ar^{o_1\circ i_1}[d]\ar[r]& 0\ar[r]\ar[d] & 0\ar[r] & \cdots\\
T_{out_1}&\cdots\ar[r] & 0\ar[r] & P_{out_1}\ar[r] & 0\ar[r] & 0\ar[r] & \cdots
}
\end{array}
$$
We claim this morphism corresponds to an arrow $in_1\rightarrow out_1$ in the quiver $Q^\prime$. Again we must be sure that the morphism does not factor through another summand of $T$. We look for a $k\neq in_1,out_1$ such that $\tilde{c}_{in_1,k}\neq 0$ and $\tilde{c}_{k, out_1}\neq 0$. There are two cases.
\\First let $k=mut$. Then by theorem \ref{altsum} $\tilde{c}_{in_1,k}=c_{in_1,in_1}+c_{in_1,in_2}-c_{in_1,mut}=1+0-1=0$. Hence the morphism $o_1\circ i_1:T_{in_1}\rightarrow T_{out_1}$ does not factor through $mut$.
\\Now suppose that $k\neq in_1,out_1,mut$ and that $\tilde{c}_{in_1\,k}\neq 0$. There must be a non-zero path $in_1$ to $k$ in $Q$. Such a path must be through $mut$ and $out_1$ but cannot end at either of these vertices (since $k\neq in_1,out_1,mut$). Therefore no non-zero path from this $k$ back to $out_1$ exists so $\tilde{c}_{k\, out_1}=0$.
\\Hence we have shown that the morphism $o_1\circ i_1:T_{in_1}\rightarrow T_{out_1}$ does not factor. Again, the same argument can be applied to non-zero paths $in_2\rightarrow mut\rightarrow out_2$ in $\Lambda$.
\\
\\[10pt]$\bullet$ \textit{if there exist arrows $pre_1\rightarrow in_1\rightarrow mut$, the arrow $pre_1\rightarrow in_1$ factors over $mut\rightarrow in_1$ in $Q^\prime$}
\\
\\[10pt] If $\eta_1:P_{pre_1}\rightarrow P_{in_1}$ is the map induced by the arrow $pre_1\rightarrow in_1$ in $Q$ then the morphism:
$$
\xymatrix{\cdots\ar[r]&0\ar[r]& P_{pre_1}\ar^{\eta_1}[d]\ar[r]&0\ar[r]&\cdots\\
\cdots\ar[r]&0\ar[r]& P_{in_1}\ar[r]&0\ar[r]&\cdots\\
}
$$
can be factored,
$$
\xymatrix{\cdots\ar[r]&0\ar[r]& P_{pre_1}\ar^{\left[\begin{smallmatrix}\eta_1\\0\end{smallmatrix}\right]}[d]\ar[r]&0\ar[r]\ar[d]&0\ar[r]&\cdots\\
\cdots\ar[r]&0\ar[r]& P_{in_1}\oplus P_{in_2}\ar^-{[i_1,i_2]}[r]\ar^{\pi_1}[d]&P_{mut}\ar[d]\ar[r]&0\ar[r]&\cdots\\
\cdots\ar[r]&0\ar[r]& P_{in_1}\ar[r]&0\ar[r]&0\ar[r]&\cdots\\
}
$$
We must show that the morphism
$$
\xymatrix{\cdots\ar[r]&0\ar[r]& P_{pre_1}\ar^{\left[\begin{smallmatrix}\eta_1\\0\end{smallmatrix}\right]}[d]\ar[r]&0\ar[r]\ar[d]&0\ar[r]&\cdots\\
\cdots\ar[r]&0\ar[r]& P_{in_1}\oplus P_{in_2}\ar^-{[i_1,i_2]}[r]&P_{mut}\ar[r]&0\ar[r]&\cdots
}
$$
corresponds to an arrow in $Q^\prime$. If this morphism were to factor non-trivially over a summand of $T$ then there exists some $k\neq pre_1,mut,in_1$ such that $\tilde{c}_{pre_1\,k}\neq 0$ and $\tilde{c}_{k\,mut}\neq 0$.
Suppose that $\tilde{c}_{pre_1\,k}\neq 0$ then we must have that $c_{k\,in_1}=c_{k\,in_2}=c_{k\,mut}=0$ (otherwise the quiver $Q$ of $\Lambda$ would have a cycle which is not oriented) so that $\tilde{c}_{k\,mut}=0$.
Similarly for arrows $pre_2\rightarrow in_2\rightarrow mut$ in $Q$ such that the composition is zero.
\\
\\[10pt]$\bullet$ \textit{other arrows}
\\
\\[10pt]All other arrows in $Q$ are transferred directly to $Q^\prime$ since theorem \ref{altsum} gives that $\tilde{c}_{ij}=c_{ij}$ for $i,j\neq mut$.
\\
\\[10pt]$\bullet$ \textit{relations}
\\
\\[10pt]To see that the relations described in theorem \ref{tilting_reg} occur in $\End_{\mathcal{D}^b(\Lambda)}(T)$ we check them individually. We demonstrate the relations in $\Lambda^\prime$ between arrows between the vertices $pre_1,in_1,mut,out_1,post_1$. The same arguments apply to the description of the relations between arrows between the vertices $pre_2,in_2,mut,out_2,post_2$.
\\Let $\eta_2$ denote the arrow in the quiver of $\Lambda$ such that the path $i_2\eta_2=0$. Then the arrow $\eta_2$ induces a morphism,
$$
\xymatrix{
\cdots\ar[r] &0\ar[r]&P_{pre_2}\ar[r]\ar[d]_{\left[\begin{smallmatrix}0\\\eta_2\end{smallmatrix}\right]}&0\ar[d]\ar[r] &0\ar[r]& \cdots\\
\cdots\ar[r] &0\ar[r]& P_{in_1}\oplus P_{in_2}\ar[r]^-{[i_1,i_2]} & P_{mut}\ar[r]& 0\ar[r] &\cdots
}
$$
such that the composition,
$$
\xymatrix{\cdots\ar[r] &0\ar[r]&P_{pre_2}\ar[r]\ar[d]_{\left[\begin{smallmatrix}0\\\eta_2\end{smallmatrix}\right]}&0\ar[d]\ar[r] &0\ar[r]& \cdots\\
\cdots\ar[r] &0\ar[r]& P_{in_1}\oplus P_{in_2}\ar[r]^-{[i_1,i_2]}\ar[d]_{\pi_1} & P_{mut}\ar[r]\ar[d]& 0\ar[r] &\cdots\\
\cdots\ar[r] &0\ar[r] & P_{in_1}\ar[r] & 0\ar[r] &0\ar[r]&\cdots
}
$$
is zero.
\\Next let $\zeta_1$ denote the arrow in the quiver of $\Lambda$ such that the path $\zeta_1 o_1=0$. Then the arrow $\zeta_1$ induces a morphism $\zeta_1:T_{out_1}\rightarrow T_{post_1}$,
$$
\xymatrix{
\cdots\ar[r] & 0\ar[r]\ar[d] & P_{out_1}\ar[r]\ar[d]_{\zeta_1} & 0\ar[r]\ar[d] &\cdots\\
\cdots\ar[r] & 0\ar[r] & P_{post_1}\ar[r] & 0\ar[r] &\cdots
}
$$
so that the composition $T_{in_1}\stackrel{o_1\circ i_1}{\rightarrow} T_{out_1}\stackrel{\zeta_1}{\rightarrow} T_{post_1}$,
$$
\xymatrix{\cdots\ar[r] & 0\ar[r]\ar[d] & P_{in_1}\ar[r]\ar[d]_{o_1\circ i_1} & 0\ar[r]\ar[d] &\cdots\\
\cdots\ar[r] & 0\ar[r]\ar[d] & P_{out_1}\ar[r]\ar[d]_{\zeta_1} & 0\ar[r]\ar[d] &\cdots\\
\cdots\ar[r] & 0\ar[r] & P_{post_1}\ar[r] & 0\ar[r] &\cdots
}
$$
is zero.
\\Consider the composition of $\pi_1\circ (o_1\circ i_1):T_{mut}\rightarrow T_{out_1}$,
$$
\xymatrix{
\cdots\ar[r] &0\ar[r]& P_{in_1}\oplus P_{in_2}\ar[r]^-{[i_1,i_2]}\ar[d]_{\pi_1} & P_{mut}\ar[r]\ar[d]& 0\ar[r] &\cdots\\
\cdots\ar[r] &0\ar[r]& P_{in_1}\ar^{o_1\circ i_1}[d]\ar[r]& 0\ar[r]\ar[d]  & 0\ar[r] &\cdots\\
\cdots\ar[r] &0\ar[r]& P_{out_1}\ar[r]& 0\ar[r] & 0\ar[r] &\cdots
}
$$
We have the following homotopy,
$$
\begin{array}{c}
\xymatrix{
\cdots\ar[r] & 0\ar[r] & P_{in_1}\oplus P_{in_2}\ar^-{[i_1,i_2]}[r]\ar^(0.4){[o_1\circ i_1,0]}[d]\ar[dl]& P_{mut}\ar[r]\ar[d]\ar^{o_1}[dl] & 0\ar[r]& \cdots\\
\cdots\ar[r] & 0\ar[r] & P_{out_1}\ar[r] & 0\ar[r] & 0\ar[r] & \cdots
}
\end{array}
$$
Hence the composition is zero in $K^{b}(P(\Lambda))$ and, hence, in $\mathcal{D}^b(\Lambda)$.
\\
\\[10pt]We can now describe all relations in the endomorphism algebra of the complex $T$. Since the class of gentle algebras is closed under derived equivalence, \cite{SZi}, $\End_{\mc{D}^b(\Lambda)}(T)\cong kQ^\prime/I^\prime$ is gentle. Therefore $I^\prime$ is generated by zero relations of length two and we have that for each arrow $\beta$ in the quiver of $\End_{\mc{D}^b(\Lambda)}(T)$ there exists at most one arrow $\alpha$ such that $\beta\alpha=0$ and at most one arrow $\gamma$ such that $\gamma\beta=0$.
\\Relations between the vertices of $\End_{\mc{D}^b(\Lambda)}(T)$ not labeled by $pre_t,in_t,mut,$
\\$out_t,post_t$, $t=1,2$ are the same as those for $\Lambda$. Indeed, let $i,j$ be two such vertices. The indecomposable summands of $T$ corresponding to these vertices are stalk complexes of the form,
$$
\xymatrix{
\cdots\ar[r] &0\ar[r]&P_i\ar[r]&0\ar[r]&\cdots\\
\cdots\ar[r] &0\ar[r]&P_j\ar[r]&0\ar[r]&\cdots
}
$$
where $P_i,P_j$ are the indecomposable projectives corresponding to vertices $i,j$ in $\Lambda$.
\\For such $i,j$, theorem \ref{altsum} gives the entries of the Cartan matrix $\tilde{C}_{\End_{\mc{D}^b(\Lambda)}(T)}$ are $\tilde{c}_{ij}=c_{ij}$. Hence, zero relations between arrows connecting vertices not labeled $pre_t,in_t,mut,out_t,post_t$, $t=1,2$ in the quiver $Q^\prime$ of $\End_{\mc{D}^b(\Lambda)}(T)$ are exactly the same as those between arrows connecting vertices not labeled $pre_t,in_t,mut,out_t,post_t$, $t=1,2$ in $\Lambda$.
\\Since we have already described the arrows and relations between the vertices of $Q^\prime$ labeled $pre_t,in_t,mut,out_t,post_t$, $t=1,2$ it only remains to check relations which are compositions of arrows coming into $pre_t$, $t=1,2$ and the arrow $pre_t\rightarrow mut$ and relations which are compositions of the arrow $out_t\rightarrow post_t$ and arrows out of $post_t$.
\\Suppose that the path $x\rightarrow pre_1\rightarrow in_1$ is zero in $\Lambda$ then the following composition of morphisms,
$$
\xymatrix{
\cdots\ar[r] &0\ar[r]\ar[d]&P_x\ar[r]\ar[d]&0\ar[d]\ar[r]&\cdots\\
\cdots\ar[r]&0\ar[r]\ar[d]&P_{pre_1}\ar[r]\ar^{\left[\begin{smallmatrix}\eta_1\\0\end{smallmatrix}\right]}[d]&0\ar[d]\ar[r]&\cdots\\
\cdots\ar[r] &0\ar[r] &P_{in_1}\oplus P_{in_2}\ar^-{[i_1,i_2]}[r]&P_{mut}\ar[r]&0\ar[r] &\cdots
}
$$
(which corresponds to arrows) is zero in $\End_{\mc{D}^b(\Lambda)}(T)$. Also if there is a zero path $out_1\rightarrow post_1\rightarrow y$ then the composition of morphisms,
$$
\xymatrix{
\cdots\ar[r] &0\ar[r]\ar[d]&P_{out_1}\ar[r]\ar^{\left[\begin{smallmatrix}\zeta_1\\0\end{smallmatrix}\right]}[d]&0\ar[d]\ar[r]&\cdots\\
\cdots\ar[r] &0\ar[r]\ar[d]&P_{post_1}\ar[r]\ar[d]&0\ar[d]\ar[r]&\cdots\\
\cdots\ar[r] &0\ar[r]& P_{y}\ar[r]& 0\ar[r] & 0\ar[r] &\cdots
}
$$
(which corresponds to arrows) is zero in $\End_{\mc{D}^b(\Lambda)}(T)$. Similarly for $t=2$.
\\Note that since $\End_{\mc{D}^b(\Lambda)}(T)$ is gentle any relations involving arrows into, or out of, vertices labeled $in_t,mut,out_t$ must be those we have already described. Hence we have described all possible relations in $\End_{\mc{D}^b(\Lambda)}(T)$.
\\This concludes the proof of \ref{tilting_reg}.
\end{proof}
The next result relates the elementary moves and the associated algebra mutations to the tilting complex just defined.

\begin{prop}\label{mm}
Let $\Lambda$ be an $m$-cluster-tilted algebra which has a local configuration of type (i)-(vi), (viii)-(x) or (xv). Then the mutated algebra $\Lambda^\prime$ is isomorphic to $\End_{\mathcal{D}^b(\Lambda)}(T)$, where $T$ is the tilting complex in theorem \ref{tilting_reg}. Hence $\Lambda$ and $\Lambda^{\prime}$ are derived equivalent.
\end{prop}

\begin{proof}We are required to prove that the elementary polygonal move $\mu_m$ performed at the mutation vertex of $\Lambda$ produces an algebra $\Lambda^\prime$ which is isomorphic to  $\End_{\mathcal{D}^b(\Lambda)}(T)$.
\\We will show that the algorithm in theorem \ref{tilting_reg} describes exactly how to produce $\Lambda^\prime$ from $\Lambda$.
\\To see that the arrows going into $mut$ are reversed consider the following local situation,
\begin{center}
\includegraphics[scale=0.5]{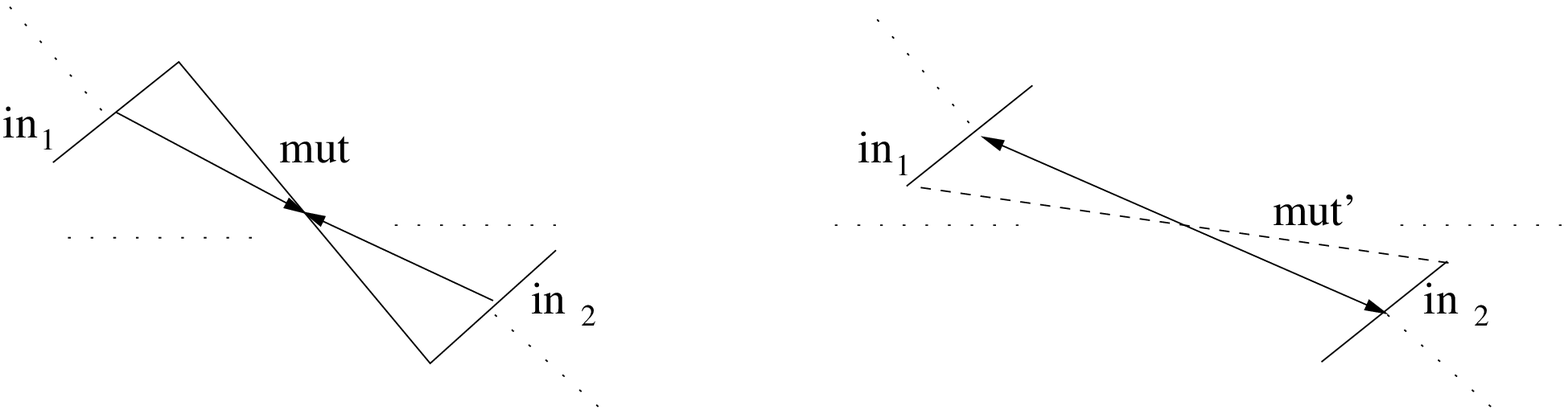}
\end{center}
Here we have labeled the $m$-allowable diagonals with the notations of theorem \ref{tilting_reg}. We make no assumptions on the division of $P$ beyond what is shown, and as such we only show the three possible $m$-allowable diagonals of relevance to this step. The same is true for the local configurations which correspond to the second and third parts of the algorithm in theorem \ref{tilting_reg}.
\\Next we examine the second step of the algorithm which states that if the path $in_1\rightarrow mut\rightarrow out_1$ exists in $\Lambda$ then we must have an arrow $in_1\rightarrow out_1$ in the quiver of $\Lambda^\prime$, the mutated algebra and the arrow $mut\rightarrow out_1$ does not exist in $\Lambda^\prime$. Again we show the local configuration of relevance in the division of $P$, the regular $(n+1)m+2$-gon.
\begin{center}
\includegraphics[scale=0.5]{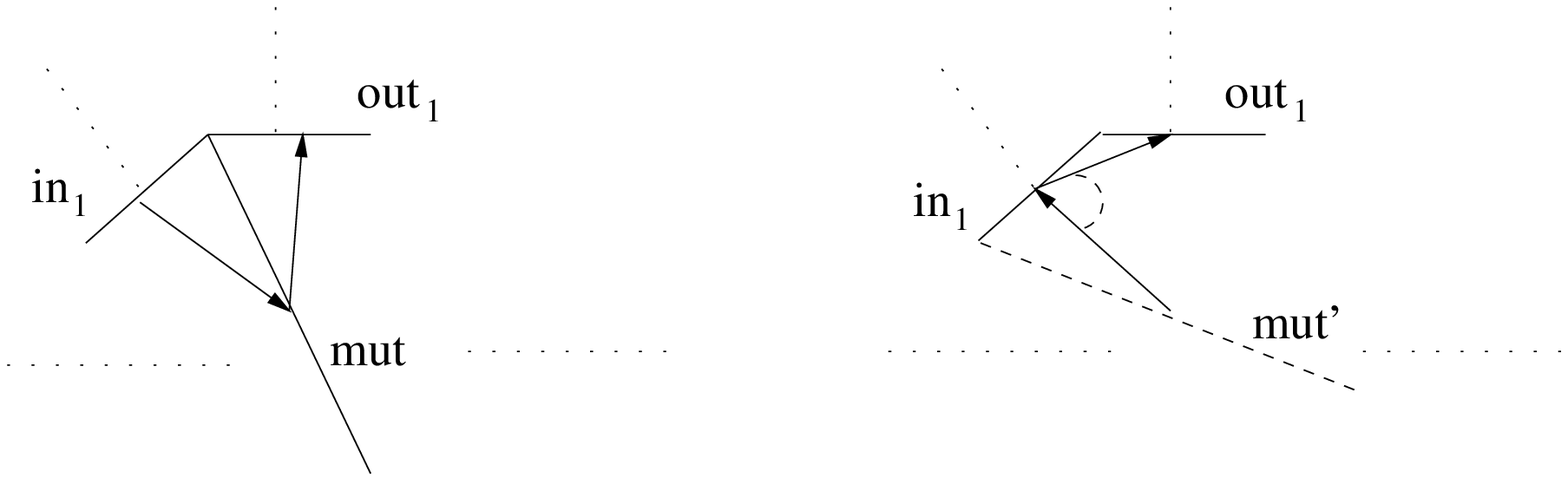}
\end{center}

Finally, if the following vertices and arrows exist in $\Lambda$,
$$pre_1\stackrel{\eta_1}{\rightarrow}in_1\stackrel{i_1}{\rightarrow}mut\stackrel{o_1}{\rightarrow}out_1\stackrel{\zeta_1}{\rightarrow}post_1
$$
$$pre_2\stackrel{\eta_2}{\rightarrow}in_2\stackrel{i_2}{\rightarrow}mut\stackrel{o_2}{\rightarrow}out_2\stackrel{\zeta_2}{\rightarrow}post_2
$$
where $i_t\eta_t=0$, $\zeta_t o_t=0$ and $o_t i_t\neq 0$, $1\leqslant t\leqslant 2$.
Then we have the following local configurations in $P$ before and after the application of $\mu_m$,
\begin{center}
\includegraphics[scale=0.5]{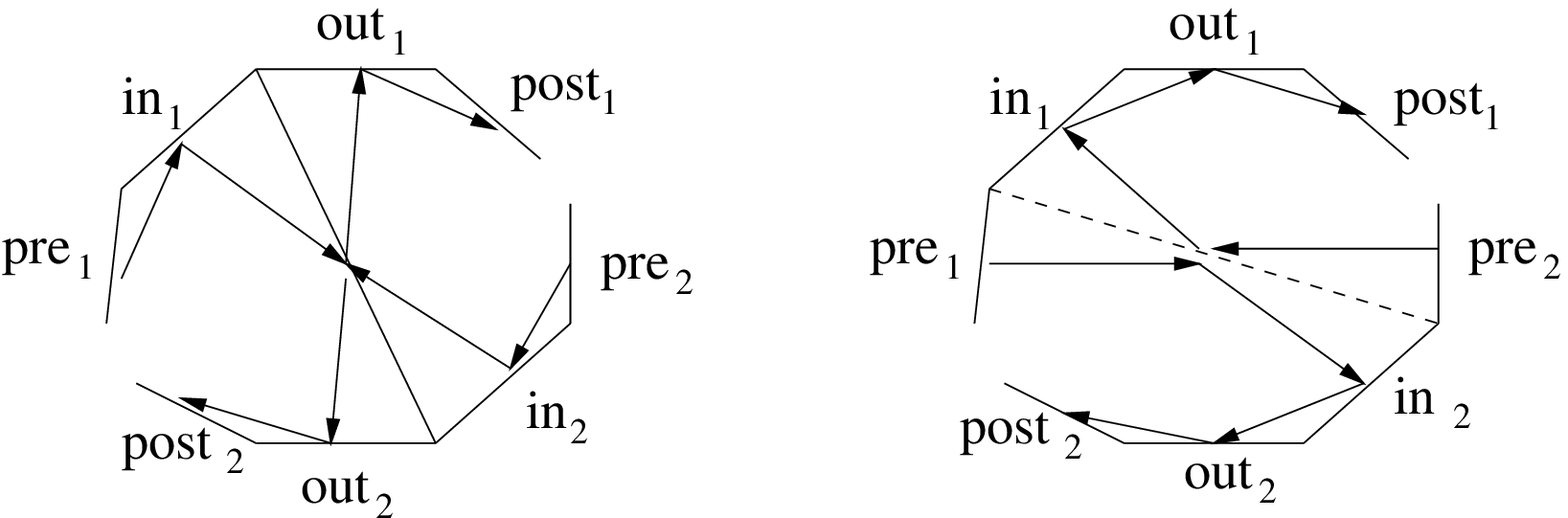}
\end{center}
so that the relations in the last part of the algorithm in the statement of theorem \ref{tilting_reg} are satisfied.
\\Hence we have shown that $\Lambda^\prime$ and $\End_{\mc{D}^b(\Lambda)}(T)$ are isomorphic algebras.
\end{proof}
Thus we have demonstrated that in the local configurations labeled (i)-(vi),(viii)-(x),(xv) we have a derived equivalence between the algebra $\Lambda$ and the mutated algebra $\Lambda^\prime$.
\\It remains to prove the same result for moves (vi),(vii),(xiv)-(xviii),(xx), and (xxii)-(xxvi) which are applications of the elementary polygonal move $\mu^{-1}_m$ defined in definition \ref{elmove}. The following two results achieve this. We will keep the notations used in theorem \ref{tilting_reg}.

\begin{thm}\label{tilting_dual}Suppose $\Lambda$ is a connected $m$-cluster-tilted algebra of type $A_n$ in which we can perform $\mu^{-1}_m$ at a given vertex, $mut$, and preserve the number of $m+2$-cycles and connectedness. Then the following complex, $T^\prime$ in $K^b(P(\Lambda))$, where $P_{mut}$ is in degree 0, is a tilting complex.
$$T^\prime:\,\,\cdots 0\rightarrow P_{mut}\stackrel{[0,o_1,o_2]}{\longrightarrow}\left(\bigoplus_{i\neq mut}P_i\right)\oplus  P_{out_1}\oplus P_{out_2}\rightarrow 0\rightarrow\cdots$$
Moreover, the endomorphism algebra, $\End_{\mathcal{D}^b(\Lambda)}(T^\prime)$, of $T^\prime$ can be obtained from $\Lambda$ using the following algorithm:
\begin{enumerate}
\item reverse all arrows outgoing from the mutation vertex.
\item if there exist arrows $in_t\rightarrow mut\rightarrow out_t$, there must exist an arrow $in_t\rightarrow out_t$ in $Q^\prime$ and the arrow $mut\rightarrow out_t$ does not exist in $q^\prime$, $t=1,2$.
\item if there exist arrows $mut\rightarrow out_t\rightarrow post_t$, the arrow
    $mut\rightarrow out_t$ factors over $out_t\rightarrow mut$ in $Q^\prime$, $t=1,2$.
\item relations in $\Lambda^\prime$ around the mutated vertex are described by the following diagrams. If arrows and relations exist around the mutation vertex in $\Lambda$ as shown in the first diagram then the second diagram shows the new arrows and relations between the corresponding vertices of $\Lambda^\prime$.
\begin{center}
\includegraphics[scale=.5]{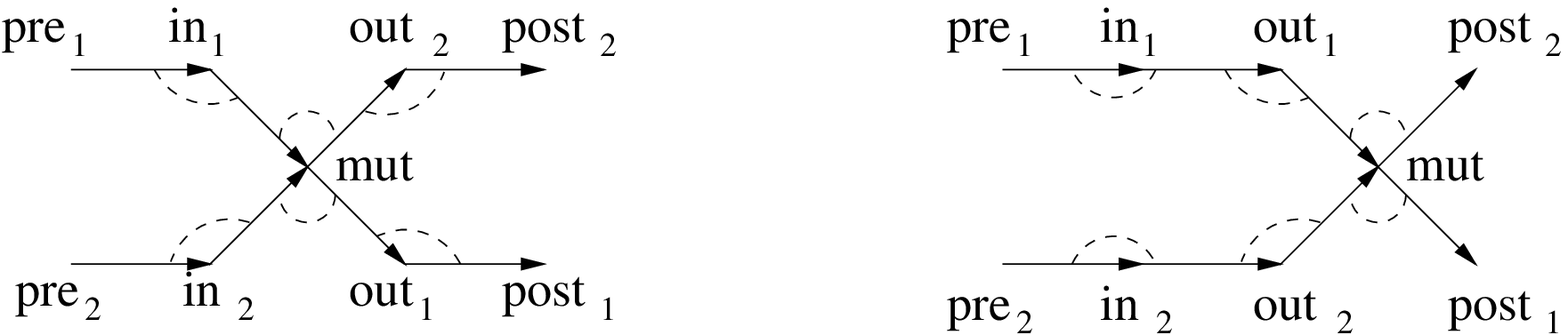}
\end{center}
\end{enumerate}
\end{thm}
\begin{rem} The comments in the remark after the statement of theorem \ref{tilting_reg} are relevant to the statement of theorem \ref{tilting_dual}, except that theorem \ref{tilting_dual} is applicable to local configurations (vi),(vii),(xiv)-(xviii),(xx), and (xxii)-(xxvi).
\end{rem}
\begin{proof}The proof is exactly analogous to the proof of theorem \ref{tilting_reg}.
\end{proof}
The next result is the analogue of proposition \ref{mm}.
\begin{prop}\label{mm1}Let $\Lambda$ be an $m$-cluster-tilted algebra which has a local configuration of type (vi),(vii),(xiv)-(xviii),(xx), and (xxii)-(xxvi). Then, if $\Lambda^\prime$ denotes the mutated algebra, $\Lambda^\prime$ is isomorphic to $\End_{\mathcal{D}^b(\Lambda)}(T^\prime)$, where $T^\prime$ is the tilting complex in \ref{tilting_dual}. Hence $\Lambda$ and $\Lambda^{\prime}$ are derived equivalent.
\end{prop}

We have now proven that in the algebra mutations described in the local configurations (i)-(xviii) $\Lambda$ and $\Lambda^\prime$ are derived equivalent.
\\[10pt]
\\Now we describe a very useful tilting complex. Suppose the following local configuration in an $m$-cluster-tilted algebra $\Lambda=kQ/I$,
\begin{center}
\includegraphics[scale=0.6]{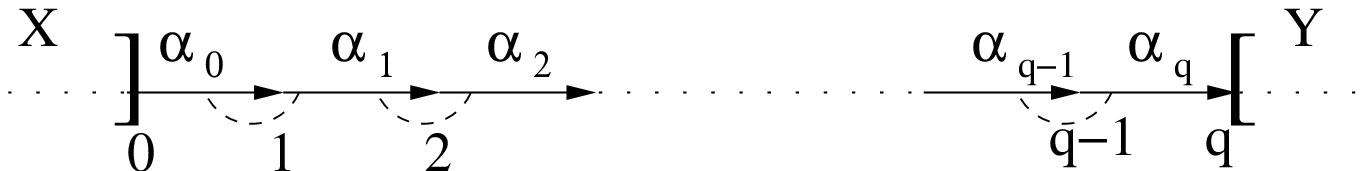}
\end{center}
where there are no paths between the regions of $Q$ labeled $X$ and $Y$; $\alpha_{i+1}\alpha_i$ is a relation for $0\leqslant i\leqslant q-2$; and there are no arrows (other than those shown) incident at vertices $1,2,\ldots,q-1$. We assume also that in the configuration shown the string of consecutive relations is maximal, that is if there is an arrow $x$ such that the composition,
$$
q-1\stackrel{\alpha_{q-1}}{\rightarrow}q\stackrel{x}{\rightarrow}j\cdots
$$
is a relation then there must be some other arrow incident at vertex $q$. Also, if there is an arrow $w$ such that the composition,
$$
i\stackrel{w}{\rightarrow}0\stackrel{\alpha_0}{\rightarrow}1\cdots
$$
is a relation there must exist another arrow incident at vertex $0$. Notice that vertices 0 and $q$ lie in regions $X$ and $Y$ respectively.
\\For a vertex $i\in X$ we define $T_i^X$ to be the stalk complex,
$$
\cdots0\rightarrow P_i\rightarrow 0\cdots
$$
where $P_i$ is in degree 0. For a vertex $j\in Y$ we define $T_j^Y$ to be the stalk complex,
$$
\cdots0\rightarrow P_j\rightarrow 0\cdots
$$
where $P_j$ is in degree q. Next for $1\leqslant r\leqslant q-1$ define $T_r$ to be the complex,
$$
\cdots 0\rightarrow P_0\stackrel{\alpha_0}{\rightarrow}P_1\stackrel{\alpha_1}{\rightarrow}\cdots\rightarrow P_{r-1}\stackrel{\alpha_{r-1}}{\rightarrow}\cdots\rightarrow P_{r}\stackrel{\alpha_{r}}{\rightarrow}0\rightarrow \cdots.
$$
Then we have the following proposition.
\begin{prop}\label{rel_rem}
Let $\Lambda$ be an $m$-cluster-tilted algebra where the local configuration described above exists. Then we have a tilting complex $T$ whose indecomposable summands are $T_i^X$, $T_j^Y$ and $T_r$, $1\leqslant r\leqslant q-1$. The endomorphism algebra, $\End_{\mc{D}^b(\Lambda)}(T)=kQ^\prime/I$, is given by preserving all arrows and relations between the vertices of regions $X$ and $Y$ and making the changes shown in the following figure to the arrows and relations between vertices $0,1,\ldots,q$.
\begin{center}
\includegraphics[scale=0.6]{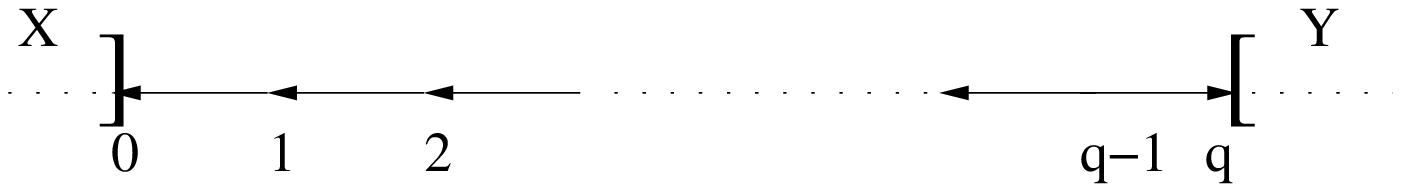}
\end{center}
\end{prop}
\begin{rem}This tilting complex will be used in section \ref{sec:algo} to help reduce the connected components of $m$-cluster-tilted algebras to the normal form. Its key feature is that it removes the relations $\alpha_{i+1}\alpha_i$, $0\leqslant i\leqslant q-2$.
\end{rem}
\begin{proof}
First we prove that $\Hom_{\mc{D}^b(\Lambda)}(T,T[k])=0$, $k\neq 0$. It suffices to show that $\Hom_{K^bP((\Lambda))}(T,T[k])=0$, $k\neq 0$.
\\Let $i$ be a vertex in the region $X$. Then we claim that for $k\neq 0$,
$$\Hom_{K^bP((\Lambda))}(T^X_i, T^Y_j[k])=0$$
for all vertices $j\in Y$. If $k<0$ then it is clear that $\Hom_{K^bP((\Lambda))}(T^X_i, T^Y_j[k])=0$. If $k>0$ then $\Hom_{K^bP((\Lambda))}(T^X_i, T^Y_j[k])=0$ since there are no non-zero paths from $X$ to $Y$ in $Q$.
\\Also for a vertex $j\in Y$ we claim that for $k\neq 0$,
$$\Hom_{K^b(P(\Lambda))}(T^Y_j,T^X_i[k])=0$$
for all vertices $i\in X$. If $k>0$ then clearly $\Hom_{K^b(P(\Lambda))}(T^Y_j,T^X_i[k])=0$. If $k<0$ then we have $\Hom_{K^b(P(\Lambda))}(T^Y_j,T^X_i[k])=0$ since there are no non-zero paths from region $Y$ to region $X$.
\\Now consider $\Hom_{K^b(P(\Lambda))}(T_i^X,T_r[k])$, $1\leqslant r\leqslant q-1$. If $k\geqslant 2$ then $\Hom_{K^b(P(\Lambda))}(T_i^X,T_r[k])=0$ since there are no non-zero paths from an vertex $i$ in $X$ to vertices $2,3,\ldots,q-1$. If $k<0$ then clearly $\Hom_{K^b(P(\Lambda))}(T_i^X,T_r[k])=0$. It remains to check that $\Hom_{K^b(P(\Lambda))}(T_i^X,T_r[1])=0$.
\\First, for any two vertices $u$ and $v$ of $Q$ such that there is a non-zero path in $Q$ from $u$ to $v$ we choose the map between $P_u$ and $P_v$ induced by the path to be a basis of $\Hom_{\mc{C}_m}(P_u,P_v)$. We can always do this since the $\Hom$ spaces in the $m$-cluster categories of type $A_n$ have dimension at most 1.
\\So, if there exists a non-zero path from some vertex $i\in X$ to vertex 1 as follows,
$$ i\stackrel{\epsilon_1}{\rightarrow}\cdots\stackrel{\epsilon_s}{\rightarrow}P_0\stackrel{\alpha_0}{\rightarrow}P_1
$$
then we have the following homotopy,
$$
\xymatrix{
\cdots\ar[r]&0\ar[r]&0\ar[r]&P_i\ar^{\epsilon}[d]\ar_{h}[dl]\ar[r]&0\ar[r]\ar_{0}[dl]&\cdots\ar[r]&0\ar[r]&0\ar[r]\ar_{0}[dl]&0\ar[r]&\cdots\\
\cdots\ar[r]&0\ar[r]&P_0\ar^{\alpha_0}[r]&P_1\ar^{\alpha_1}[r]&P_2\ar^{\alpha_2}[r]&\cdots\ar[r]&P_{r-1}\ar^-{\alpha_{r-1}}[r]&P_r\ar[r]
&0\ar[r]&\cdots
}
$$
where $\epsilon=\epsilon_s\circ\epsilon_{s-1}\circ\cdots\circ\epsilon_1\alpha_0$ is a basis of and $\Hom_{\mc{C}_m}(P_i,P_1)$ and $h=\epsilon_s\circ\epsilon_{s-1}\circ\cdots\circ\epsilon_1$. Hence,
$\Hom_{K^b(P(\Lambda))}(T_i^X,T_r[1])=0$ and we have proven that $\Hom_{K^b(P(\Lambda))}(T_i^X,T_r[k])=0$, $k\neq 0$.
\\Next we show that $\Hom_{K^b(P(\Lambda))}(T_r,T_i^X[k])=0$, $k\neq 0$, $1\leqslant r\leqslant q-1$ and $i\in X$. If $k>0$ then the claim is clearly true. If $k<0$ then $\Hom_{K^b(P(\Lambda))}(T_r,T_i^X[k])=0$ since there are no non-zero paths in $Q$ from vertices $1,2,\ldots,q-1$ to any vertex in region $X$.
\\We also have that $\Hom_{K^b(P(\Lambda))}(T_r,T_j^Y[k])=0$, $k\neq 0$. If $k<0$ then it is clear that $\Hom_{K^b(P(\Lambda))}(T_r,T_j^Y[k])=0$. If $k>0$ then $\Hom_{K^b(P(\Lambda))}(T_r,T_j^Y[k])=0$ since there are no non-zero paths in $Q$ from vertices $1,2,\ldots,q-2$ to any vertex in region $Y$. Next consider, $\Hom_{K^b(P(\Lambda))}(T_j^Y,T_r[k])$ for any vertex $j\in Y$. If $k>0$ then it is clear that $\Hom_{K^b(P(\Lambda))}(T_j^Y,T_r[k])=0$. If $k<0$ then we have that $\Hom_{K^b(P(\Lambda))}(T_j^Y,T_r[k])=0$ since there are no non-zero paths in $Q$ from the vertices in region $Y$ to vertices
\\$0,1,2,\ldots,q-2$.
\\
\\[10pt]It remains to prove that $\Hom_{K^b(P(\Lambda))}(T_{r_1},T_{r_2}[k])=0$, $1\leqslant r_1,r_2\leqslant q-1$. For $
k\geqslant 2$ or $k<0$ there are no possible non-zero morphisms $T_{r_1}\rightarrow T_{r_2}[k]$ so in these cases $\Hom_{K^b(P(\Lambda))}(T_{r_1},T_{r_2}[k])=0$.
\\Next we consider $\Hom_{K^b(P(\Lambda))}(T_{r_1},T_{r_2}[1])$. There are two separate cases.
First let $r_2\leqslant r_1$ and let $\phi\in\Hom_{K^b(P(\Lambda))}(T_{r_1},T_{r_2}[1])$ be a non-zero morphism. We now prove that $\phi$ is homotopic to zero. Let $\phi_i=\alpha_i$, $0\leqslant i\leqslant r_2-1$ be the first non-zero component of the map $\phi$. We can assume that this is the case by scaling the morphism $\phi$ appropriately.
\\The following figure illustrates the construction of the homotopy. The component $h_i:P_i\rightarrow P_i$ of the homotopy must be the zero map, therefore to satisfy the definition of a homotopy $h_{i+1}$ must be the identity. There are now two possibilities for $\phi_{i+1}=\lambda_{i+1}\alpha_{i+1}$ either it is zero, in which case $h_{i+2}=-id$, or it is not, in which case $h_{i+2}=(\lambda_{i+1}-1)id$. In both cases we have that $h_{i+2}$ is a scalar multiple of the identity so we can write $h_{i+2}=\mu_{i+2}id$. We can continue to construct the components of the homotopy in this manner, at each stage we can deduce that $h_{i+p}=\mu_{i+p}id$ where $p\geqslant 1$. The next figure shows the homotopy we have described.
$$
\xymatrix{
\cdots\ar[r]&P_{i-1}\ar^{\alpha_i}[r]\ar^(0.4){0}[dd]&P_i\ar^{\alpha_i}[r]\ar^(0.4){\alpha_i}[dd]\ar^{0}[ddl]&P_{i+1}\ar^{\alpha_{i+1}}[rr]
\ar^(0.4){\lambda_{i+1}\alpha_{i+1}}[dd]\ar^{1}[ddl]&&P_{i+2}\ar^{\alpha_{i+2}}[rr]\ar^(0.4){\lambda_{i+2}\alpha_{i+2}}[dd]\ar^{\mu_{i+2}1}[ddll]
&&P_{i+3}\ar^{\alpha_{i+3}}[r]\ar^(0.4){\lambda_{i+3}\alpha_{i+3}}[dd]\ar^{\mu_{i+3}1}[ddll]&\cdots\\
&&&&&&&&&&&\\
\cdots\ar[r]&P_{i}\ar^{\alpha_i}[r]& P_{i+1}\ar^{\alpha_{i+1}}[r] & P_{i+2}\ar^{\alpha_{i+2}}[rr] && P_{i+3}\ar^{\alpha_{i+3}}[rr]&& P_{i+4}\ar^{\alpha_{i+4}}[r] & \cdots
}
$$
Eventually we encounter the following situation.
$$
\xymatrix{
\cdots\ar[r]&P_{r_2-2}\ar^{\alpha_{r_2-2}}[rr]\ar^(0.4){\lambda_{r_2-2}\alpha_{r_2-2}}[dd]&& P_{r_2-1}\ar^{\alpha_{r_2-1}}[rr]\ar^(0.4){\lambda_{r_2-1}\alpha_{r_2-1}}[dd]\ar^{\mu_{r_2-1}1}[ddll]&&
P_{r_2}\ar^{\alpha_{r_2}}[r]\ar^(0.4){0}[dd]\ar^{\mu_{r_2}1}[ddll]&
P_{r_2+1}\ar[r]\ar[dd]\ar^{0}[ddl]&
\cdots\ar[r]&P_{r_1}\ar[r]\ar[dd]&0\ar[r]&\cdots\\
&&&&&&&&&&&&&\\
\cdots\ar[r]&P_{r_2-1}\ar^{\alpha_{r_2-1}}[rr]&& P_{r_2}\ar[rr] && 0\ar[r] & 0\ar[r]&\cdots\ar[r]& 0\ar[r] &0\ar[r]& \cdots
}
$$
So we have that setting $h_{s}=0$, $s\geqslant r_2$ defines a homotopy and therefore since $\phi$ was arbitrary $\Hom_{K^b(P(\Lambda))}(T_{r_1},T_{r_2}[1])=0$.
\\Now we consider the second case, that is assume now that $r_2>r_1$. By theorem \ref{altsum} we have that,
$$
(-1)dim_k\Hom_{\mc{D}^b(\Lambda)}(T_{r_1},T_{r_2}[1])=\sum_{a,b}(-1)^{a-b}dim_k\Hom_\Lambda(T_{r_1}^a,T_{r_2}^b)
$$
where $T_{r_1}^a$ is the degree $a$ term of $T_{r_1}$ and $T_{r_2}^b$ is the degree $b$ term of $T_{r_2}$.
\\It follows that,
\begin{eqnarray*}
(-1)dim_k\Hom_{\mc{D}^b(\Lambda)}(T_{r_1},T_{r_2}[1])&=&-c_{0,0}+c_{0,1}-c_{1,1}+c_{1,2}-c_{2,2}+\cdots{}\\
&&{}\cdots-c_{i,i}+c_{i,i+1}-\cdots{}\\
&&{}\cdots-c_{r_1-1,r_1-1}+c_{r_1-1,r_1}-c_{r_1,r_1}+c_{r_1,r_1+1}\\
&=&0
\end{eqnarray*}
where the $c_{y,z}$ are entries of the Cartan matrix $C_\Lambda$ of $\Lambda$. Note that for $0\leqslant y,z\leqslant q-1$ we have $c_{y,z}=0$ if $y>z$. Also, $c_{y,z}=0$ if $z\geqslant y+2$.
Therefore, if $r_2>r_2$ we have that $dim_k\Hom_{\mc{D}^b(\Lambda)}(T_{r_1},T_{r_2}[1])=0$.
\\We have now proven that $\Hom_{\mc{D}^b(\Lambda)}(T,T[k])=0$, $k\neq 0$. Arguments similar to those used in theorems \ref{tilting_reg} and \ref{tilting_dual} give that $add(T)$ generates $K^b(P(\Lambda))$ as a triangulated category and so $T$ is a tilting complex.
\\
\\[10pt]Next we must describe $\End_{\mc{D}^b}(T)$. First we note that arrows and relations between vertices of the quiver $Q^\prime$ of $\End_{\mc{D}^b}(T)$ corresponding to the $T_i^X$s are the same as the arrows and relations between vertices of type $X$ in the quiver $Q$ of $\Lambda$. Similarly, the vertices of the quiver $Q^\prime$ corresponding to the $T_j^Y$s are the same as the arrows and relations between vertices of type $Y$ in $Q$. To see this let $\tilde{c}_{y,z}$ denote the entries of the Cartan matrix $\tilde{C}$ of $\End_{\mc{D}^b}(T)$ and consider $i_1,i_2\in X$. By theorem \ref{altsum} we have that $\tilde{c}_{i_1,i_2}=c_{i_1,i_2}$. Also for $j_1,j_2\in Y$ theorem \ref{altsum} gives $\tilde{c}_{j_1,j_2}=c_{j_1,j_2}$.
\\Next consider the following morphism $T_r\stackrel{\theta}{\rightarrow} T_{r-1}$, $1\leqslant r\leqslant q-1$,
$$
\xymatrix{
\cdots\ar[r] &0\ar[r]\ar[d] & P_{0}\ar^{\alpha_0}[r]\ar^{\theta_0}[d] & P_{1}\ar^{\alpha_1}[r]\ar^{\theta_1}[d] & P_{2}\ar^{\alpha_2}[r]\ar^{\theta_2}[d] &\cdots\ar[r] & P_{r-1}\ar^{\alpha_{r-1}}[r]\ar^{\theta_{r-1}}[d] &P_{r}\ar[r]\ar[d] &0\ar[r] &\cdots\\
\cdots\ar[r] &0\ar[r] & P_{0}\ar^{\alpha_0}[r] & P_{1}\ar^{\alpha_1}[r] & P_{2}\ar^{\alpha_2}[r] &\cdots\ar[r] & P_{r-1}\ar[r] &0\ar[r] &\cdots
}
$$
By theorem $\ref{altsum}$ we can easily deduce that $dim_k\Hom_{\mc{D}^b(\Lambda)}(T_r,T_{r-1})=1$ (the argument is similar to the argument above which shows that $dim_k\Hom_{\mc{D}^b(\Lambda)}(T_{r_1},T_{r_2}[1])=1$). Therefore we can assume that $\theta$ is a basis for $\Hom_{\mc{D}^b(\Lambda)}(T_r,T_{r-1})$.
\\If $\theta$ is to correspond to an arrow in $Q^\prime$we must show that it does not factor over another summand of $T$.
\\If $\theta$ were to factor over a summand of $T$ then there would exist some $k\neq r,r-1$ such that $\tilde{c}_{r,k}\neq 0$ and $\tilde{c}_{k,r-1}\neq 0$. Note that for any $T_j^Y$ there are no non-zero morphisms $T_j^Y\rightarrow T_{r-1}$, hence its is not possible for $T_r\stackrel{\theta}{\rightarrow} T_{r-1}$ to factor over any $T_j^Y$.
\\Now let $i\in X$ and consider $T_i^X$. Then by theorem \ref{altsum} $\tilde{c}_{i,r-1}=c_{i,0}-c_{i,1}$. There are two possible cases. The first is that $c_{i,0}=c_{i,1}=1$ and the second is that $c_{i,0}=c_{i,1}=0$. In either case $\tilde{c}_{i,r-1}=0$ so that $\theta$ does not factor over any $T_i^X$.
\\It remains to show that $\theta$ does not factor over any $T_s$, $1\leqslant s\leqslant q-1$, $s\neq r,r-1$. Let $s>r$ and consider a morphism $f:T_r\rightarrow T_s$. We claim that this morphism is zero. Suppose that $f$ is non-zero and that $P_i\stackrel{f_i}{\rightarrow} P_i$ is the last non-zero component of the morphism $f$, then by scaling $\theta$ we can assume that $f_i$ is the identity. The following diagram shows that $f$ cannot be a morphism and hence we have a contradiction.
$$
\xymatrix{
\cdots\ar^{\alpha_{i-1}}[r]&P_i \ar^{1}[d]\ar^{\alpha_{i}}[r]&P_{i+1}\ar[r]\ar^{0}[d]&\cdots\ar[r]&P_r\ar[r]\ar^{0}[d]&0\ar[r]\ar[d]&\cdots\\
\cdots\ar^{\alpha_{i-1}}[r]&P_i \ar^{\alpha_{i}}[r]&P_{i+1}\ar[r]&\cdots\ar[r]&P_r\ar[r]&P_{r+1}\ar[r]&\cdots\ar[r]&P_s\ar[r]&0\ar[r]&\cdots
}
$$
It follows that
the map $\theta$ does not factor over $T_s$, $s>r$. Similarly, the map $\theta$ cannot factor over $T_s$, $s<r-1$ and so we have proven that $\theta:T_r\rightarrow T_r-1$ corresponds to an arrow in the quiver $Q^\prime$ of $\End_{\mc{D}^b(\Lambda)}(T)$, $1\leqslant r\leqslant q-1$.
\\Next consider the non-zero morphism $T_{q-1}\rightarrow T_q$,
$$
\xymatrix{
\cdots\ar[r]&P_{q-2}\ar^{\alpha_{q-2}}[r]\ar[d]&P_{q-1}\ar[r]\ar^{\alpha_{q-1}}[d]&0\ar[r]&\cdots\\
\cdots\ar[r]&0\ar[r]&P_q\ar[r]&0\ar[r]&\cdots
}
$$
If this morphism were to factor over a summand of $T$ then there would exist some $k\neq q-1,q$ such that $\tilde{c}_{q-1,k}\neq 0$ and $\tilde{c}_{k,q}\neq 0$. It is clear that if $k$ is a vertex in region $X$ or if $k=q-2,q-3,\ldots,1$ then $\tilde{c}_{k,q}=0$. Now let $k\neq q$ be a vertex in region $Y$. Then by theorem \ref{altsum} we have that $\tilde{c}_{q-1,k}=c_{q-1,k}$. Hence it is not possible that $\tilde{c}_{q-1,k}\neq 0$ and $\tilde{c}_{k,q}\neq 0$, since this would create an configuration in $Q$, the quiver of $\Lambda$, which would contradict the description of $m$-cluster-tilted algebras in section \ref{sec:m_clust}.
\\It remains to state that there are no relations in $\Lambda$ of the form $1\rightarrow 0\rightarrow i$ for some $i\in X$ and that any relation of the form $q-1\rightarrow q\rightarrow j$ in $\Lambda$ is preserved in $\Lambda^\prime$. Both of these statements are clear when we consider the following compositions of morphisms. We have scaled the morphisms where we can.
$$
\xymatrix{
\cdots\ar[r]&0\ar[r]&P_0\ar^{\alpha_0}[r]\ar^{1}[d]&P_1\ar[r]\ar[d]&0\ar[r]&\cdots\\
\cdots\ar[r]&0\ar[r]&P_0\ar[r]\ar^{\nu}[d]&0\ar[r]\ar[d]&0\ar[r]&\cdots\\
\cdots\ar[r]&0\ar[r]&P_i\ar[r]&0\ar[r]&0\ar[r]&\cdots
}
$$
$$
\xymatrix{
\cdots\ar[r]&P_{q-3}\ar^{\alpha_{q-3}}[r]&P_{q-2}\ar^{\alpha_{q-2}}[r]\ar[d]&P_{q-1}\ar[r]\ar^{\alpha_{q-1}}[d]&0\ar[r]&\cdots\\
\cdots\ar[r]&0\ar[r]&0\ar[r]\ar[d]&P_q\ar[r]\ar^{\nu^\prime}[d]&0\ar[r]&\cdots\\
\cdots\ar[r]&0\ar[r]&0\ar[r]&P_j\ar[r]&0\ar[r]&\cdots
}
$$

We have now completed the proof of the proposition since by \cite{SZi} we know that $\End_{\mc{D}^b(\Lambda)}(T)$ is gentle, so there are no other possible relations to describe.
\end{proof}

This complex and the derived equivalence it induces will be important in our reduction of a connected component to the normal form.
\begin{ex} We include explicit examples of one of each of the three types of tilting complexes we have outlined in this section. The figures are followed by the relevant tilting complex.
\\
\\[10pt]An example of the tilting complex in theorem \ref{tilting_reg}.
\begin{center}
\includegraphics[scale=0.5]{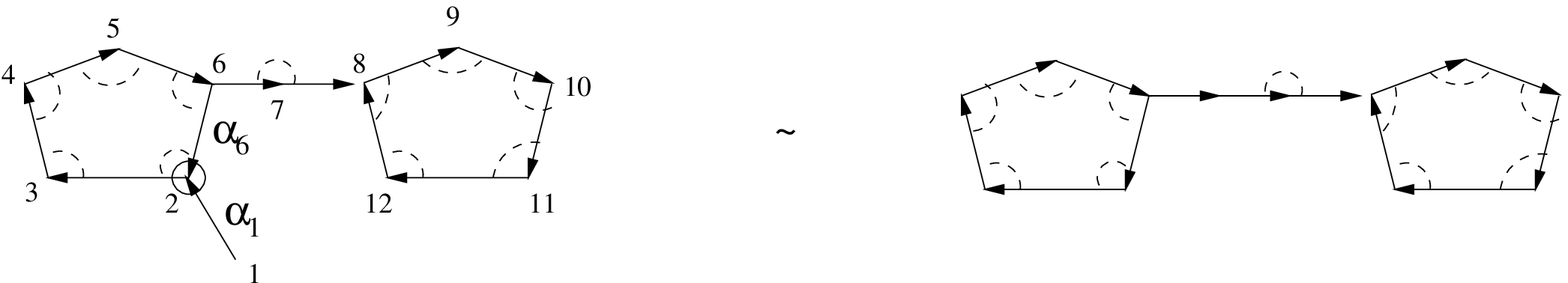}
\end{center}
$$
\cdots\rightarrow 0\rightarrow P_1\oplus P_6\oplus(\bigoplus_{i\neq 2} P_i)\stackrel{[\alpha_1,\alpha_6,0]}{\longrightarrow}
P_2\rightarrow 0\rightarrow\cdots
$$
An example of the tilting complex in theorem \ref{tilting_dual}.
\begin{center}
\includegraphics[scale=0.5]{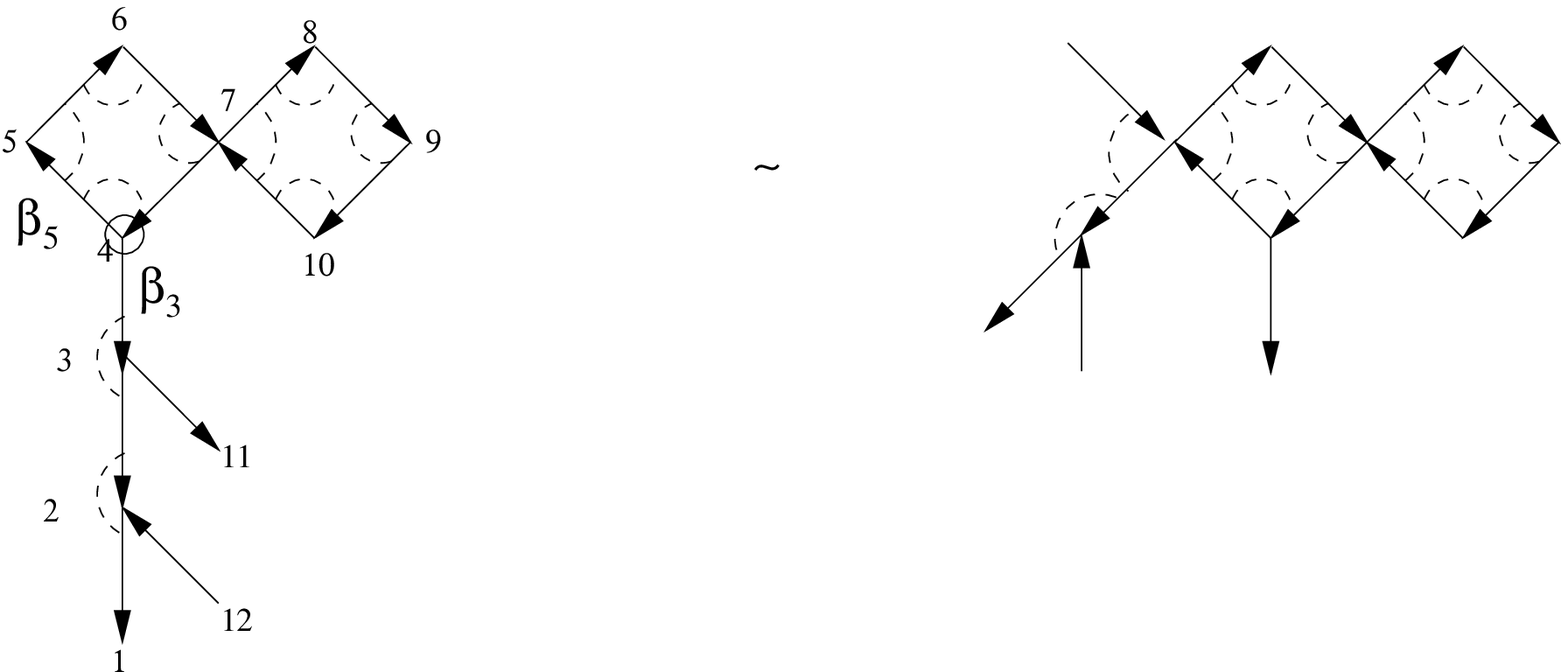}
\end{center}
$$
\cdots\rightarrow 0\rightarrow P_4\stackrel{[0,\beta_5,\beta_3]}{\longrightarrow}
(\bigoplus_{i\neq 4} P_i)\oplus P_5\oplus P_3\rightarrow 0\rightarrow\cdots
$$
An example of the tilting complex in theorem \ref{rel_rem}.
\begin{center}
\includegraphics[scale=0.5]{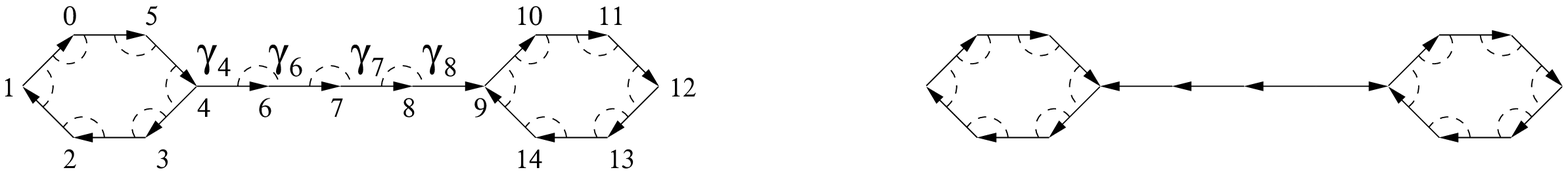}
\end{center}
We give the indecomposable summands for the tilting complex relating to the above figure.
$$
\cdots\rightarrow 0\rightarrow P_i\rightarrow 0\rightarrow\cdots
$$
$0\leqslant i\leqslant 5$ and the $P_i$ terms are in degree 0.
$$
\cdots\rightarrow 0\rightarrow P_j\rightarrow 0\rightarrow\cdots
$$
$9\leqslant j\leqslant 14$ and the $P_j$ terms are in degree 4.
$$
\cdots\rightarrow 0\rightarrow P_4\stackrel{\gamma_4}{\rightarrow}P_6\rightarrow 0\rightarrow \cdots
$$
$$
\cdots\rightarrow 0\rightarrow P_4\stackrel{\gamma_4}{\rightarrow}P_6\stackrel{\gamma_6}{\rightarrow}P_7\rightarrow 0\rightarrow \cdots
$$
$$
\cdots\rightarrow 0\rightarrow P_4\stackrel{\gamma_4}{\rightarrow}P_6\stackrel{\gamma_6}{\rightarrow}P_7\stackrel{\gamma_7}{\rightarrow}P_8\rightarrow 0\rightarrow \cdots
$$
$$
\cdots\rightarrow 0\rightarrow P_4\stackrel{\gamma_4}{\rightarrow}P_6\stackrel{\gamma_6}{\rightarrow}P_7\stackrel{\gamma_7}{\rightarrow}P_8\stackrel{\gamma_8}{\rightarrow}P_9\rightarrow 0\rightarrow \cdots
$$
where the $P_4$ term always lies in degree zero.
\end{ex}
\section{Reduction to Normal Form}\label{sec:algo}
In this section we provide the procedure to reduce any connected component of an $m$-cluster-tilted algebra, $\Lambda$, to the normal form, defined in definition \ref{nformdef}, using the local mutations (i)-(xxvi) from section \ref{sec:tilt}. Once this procedure has been demonstrated, we will have proven theorem \ref{main}.
\begin{thm}\label{norm_form}Every connected component $\Lambda$ with $r\in\mathbb{N}$ $m+2$-cycles and $s\in\mathbb{N}$ simple modules is derived equivalent to the normal form,
\begin{center}
\includegraphics[scale=0.5]{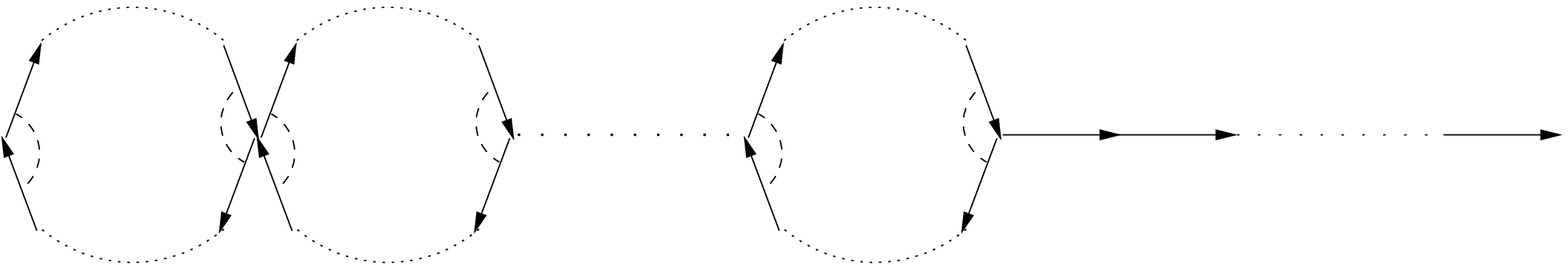}
\end{center}
where there are  $r$ cycles and a total of $s$ vertices.
\end{thm}
\begin{proof}Let $\Lambda$ be as in theorem \ref{norm_form}. The following algorithm proves theorem \ref{main} when combined with theorems \ref{tilting_reg}, \ref{tilting_dual} and proposition \ref{rel_rem}.
\\We assume that the quiver of $\Lambda$ contains at least one $m+2$-cycle. We will see later how to deal with a connected component which has no $m+2$-cycles.
\\Since every $m$-cluster-tilted algebra, $\Lambda$ can be described by a maximal division of a regular convex $(n+1)m+2$-gon, $P$, by $m$-allowable diagonals (see section \ref{sec:m_clust}) it follows that in the quiver, $Q$, of $\Lambda$ there exists a $m+2$-cycle which is connected by paths in $Q$ to other $m+2$-cycles through at most one vertex. Such a cycle will be called an \textit{initial cycle}. Notice that initial cycles need not be unique and if $Q$ has only one cycle this cycle is clearly initial.
\\It is important to note that in the following that whenever we apply the elementary polygonal move $\mu_m$ in the polygon $P$ and consider the corresponding algebra mutation we achieve a derived equivalence via theorem \ref{tilting_reg}. Similarly, whenever we apply the elementary polygonal move $\mu_m^{-1}$ in $P$ and consider the corresponding algebra mutation we achieve a derived equivalence via theorem \ref{tilting_dual}.
\\
\\[10pt]Now choose an initial cycle and label the vertices as follows. The vertex which is (potentially) connected to other cycles will be labeled 0 and the other vertices are labeled $\{1,2,\ldots,m+1\}$ in an anti-clockwise direction around the cycle (that is, in the opposite direction to the orientation which we take to be clockwise).
\\Starting at vertex $i=1$ initiate the following algorithm:
\begin{enumerate}
\item if there exists an arrow into vertex $i$ which is not part of the initial cycle then perform $\mu_m$ at this vertex. That is, apply the elementary move $\mu_m$, defined in definition \ref{elmove}, to the $m$-allowable     diagonal corresponding to vertex $i$ in the division of $P$ which is associated with $\Lambda$. The possible local configurations and resultant algebras relevant to this step are shown in figures (ii) and (x)-(xiii).
    \begin{rem} Note that figure (xiii) shows two $m+2$-cycles which share a vertex. The mutation vertex is this shared vertex and the mutation shown corresponds to $\mu_m$. This situation is actually a special case of figure (xi), namely the case where the quivers labeled $Q_2$ and $Q_5$ join together. Figure (xii) is also a special case of figure (xi). The list of figures is extended for clarity, since it may not be immediately apparent that some figures are special cases of the others. For example, figures(i), (iii) and (iv) could also be united into a single diagram.
    \end{rem}
\item if there exists no such arrow, but there exists an arrow out of $i$ which is again not part of the initial cycle perform $\mu_m$ at the target of this outgoing arrow, \textit{except} where we encounter the configurations shown in (vi) and (vii). In the these cases we must apply $\mu_m^{-1}$. The possible local configurations and resultant algebras relevant to this step are shown in figures (i), (iii), (iv), (v), (vi), (vii), (viii) and (ix).
    There is one other situation which can occur and must be treated slightly differently. Suppose that the following local situation occurs:
    \begin{center}
    \includegraphics[scale=0.5]{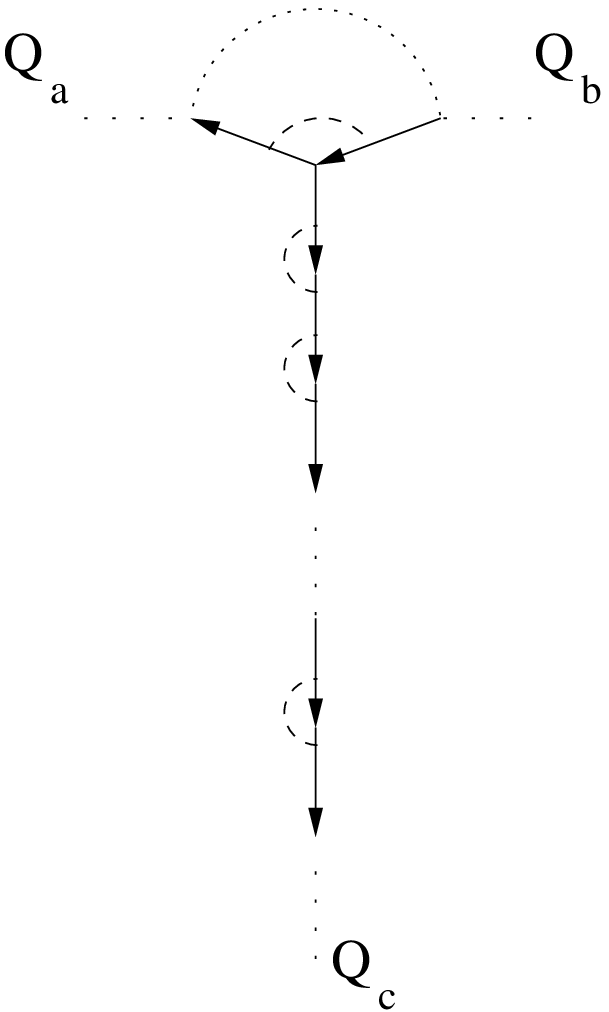}
    \end{center}
    Namely, from vertex $i$ we have an oriented path of length $s$ where $2\leqslant s\leqslant m+1$ in which each composition of two consecutive arrows is a zero relation.
    \\We can apply proposition \ref{rel_rem} to achieve a derived equivalence with an algebra with the following local configuration:
    \begin{center}
    \includegraphics[scale=0.5]{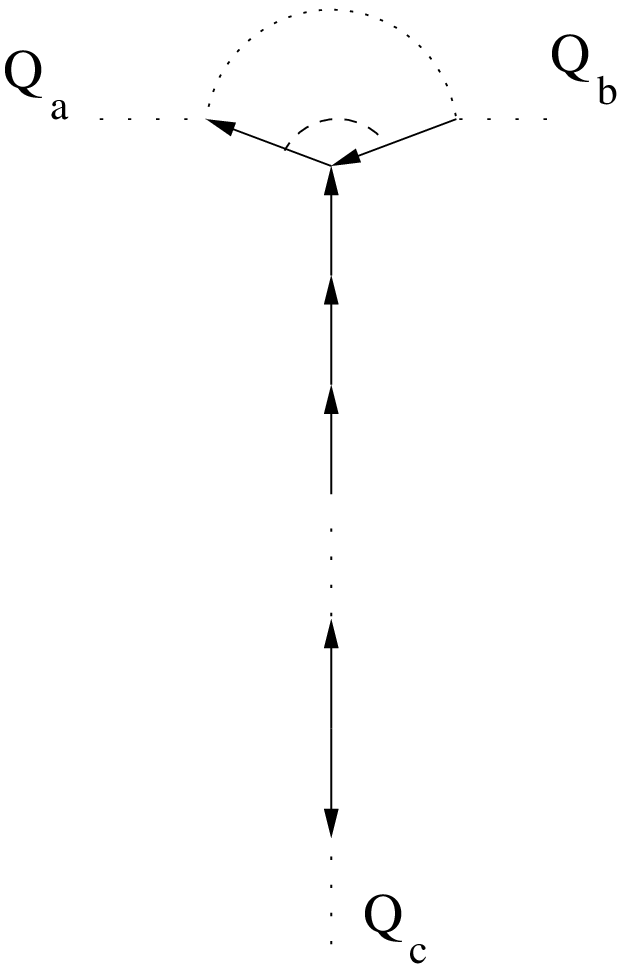}
    \end{center}
    and elsewhere has the same quiver with relations as our original algebra.
\item repeat steps one and two until there are no arrows into or out of $i$ which are not part of the initial cycle.
\item repeat the first three steps for vertex $i+1$, unless $i=m+1$ in which case stop.
    \begin{rem}Suppose at a given vertex $i$ of the initial cycle we have an arrow out of $i$ which is not part of the initial cycle then the application of the second step always serves to create a situation where there is an arrow \textit{into} $i$ which is not part of the initial cycle and so the first step can then be applied.
    \end{rem}
\end{enumerate}
The algorithm halts when there are no arrows which are not part of the chosen initial cycle into or out of vertex $m+1$.
To complete the procedure we only use local moves (i)-(xiii), as listed in section \ref{sec:tilt}, and proposition \ref{rel_rem}. After the application of the above steps only vertex 0 has any arrows incident or outgoing which are not part of the initial cycle. We say that we have \textit{cleared} the initial cycle.
\\
\\[10pt]Now assume that further $m+2$-cycles exist in $\Lambda$. We must choose and be able to clear the second $m+2$-cycle.
\\It is convenient to choose an $m+2$-cycle which is closest to the initial cycle (that is, this second cycle is linked to the initial cycle by a shortest un-oriented path in the quiver, but in general the next procedure works for an arbitrary choice).
\\Once the second $m+2$-cycle has been chosen label its vertices as follows. Label the vertex which is connected to the initial cycle with 0. Next label the remaining vertices $\{1,2,\ldots,m+1\}$ following the clockwise orientation of the cycle. Define vertices of type $B$ to be those which are numbered $\{1,2,\ldots,(\frac{m-1}{2})\}$ if $m$ is odd and those which are numbered $\{1,2,\ldots,\frac{m}{2}\}$ if $m$ is even. Define vertices of type $A$ to be those labeled  $\{(\frac{m+4}{2}),\ldots,m+1\}$ if $m$ is even and to be those labeled $\{(\frac{m+3}{2}),\ldots,m+1\}$ if $m$ is odd.
\begin{rem}
We must ensure that our algorithm produces a quiver with the same arrangement of $m+2$-cycles as in the normal form defined in definition \ref{nformdef}. The above definition of vertices of types $A$ and $B$ will coincide with the definition of vertices of types $A$ and $B$ in definition \ref{nformdef}. This is not clear now, but will become apparent in due course.
\end{rem}
The algorithm used to clear the initial cycle is, in general, no longer sufficient to clear our chosen second cycle. However, it can still be applied to vertices in the second cycle of type $A$. We must describe what to do at vertex 0 and at vertices of type $B$ in the second cycle.
\\At vertex 0 in the second cycle it is our aim to describe a procedure which results in the initial cycle and second cycle sharing a vertex, but which leaves the initial cycle otherwise unaffected. This procedure is made slightly more complicated by this requirement.
\\At vertex 0 in the second cycle initiate the following steps:
\begin{enumerate}
\item if there exists an arrow into the vertex 0 in the second cycle which is not part of the initial or second cycle then apply $\mu_m^k$ at 0 in the second cycle. If the arrow is the last arrow in a chain of $t$ consecutive relations, $1\leqslant t\leqslant m-1$,  then the index $k=t+1$. If the arrow is not the last arrow in a chain of relations then $k=1$. The possible figures relevant to this step are (ii), (x), (xi), (xii) and (xiii).
    \begin{rem}We saw in section \ref{sec:m_clust} that outwith an $m+2$-cycle there can be at most $m-1$ consecutive relations.
    \end{rem}
\item if there exists an arrow out of the vertex 0 in the second cycle which is not part of the initial or second cycle then apply $\mu_m^{-k}$ at 0 in the second cycle. If the arrow is the first arrow in a chain of $t$ consecutive relations, $1\leqslant t\leqslant m-1$,  then the index $k=t+1$. If the arrow is not the last arrow in a relation then $k=1$. The possible figures relevant to this step are (xvi), (xxiii)-(xxvi).
\item repeat steps one and two until there are no arrows into or out of 0 which are not part of the initial cycle or second cycle.
\end{enumerate}
The result is that now the quiver of $\Lambda$ consists of the initial $m+2$-cycle which is connected at one vertex to the second  $m+2$-cycle. In the initial cycle there are no arrows into or out of vertices which are not part of the initial cycle (with the sole exception of the vertex connected to the second cycle). Notice that, unlike the initial cycle, the second $m+2$-cycle can be connected to many other $m+2$-cycles.
\\
\\[10pt]Next, to clear the vertices of type $A$ (that is vertex $i$, $(\frac{m+3}{2})\leqslant i\leqslant m+1$ or $(\frac{m+4}{2})\leqslant i\leqslant m+1$, if $m$ is odd or even respectively) we apply exactly the same procedure as when clearing the initial cycle. We always begin with $i=m+1$ and proceed to $m, m-1,\ldots,x$ where $x=(\frac{m+3}{2})$ or $x=(\frac{m+4}{2})$ if $m$ is odd or even respectively.
\\Once the process has been applied at all vertices of type $A$ there will be no arrows into or out of these vertices which are not part of the second cycle.
\\
\\[10pt]We require a different set of instructions to clear the vertices of type $B$ without affecting the already cleared initial cycle.
\\Starting at vertex $j=1$ initiate the following procedure.
\begin{enumerate}
\item If at vertex $j$ there is an arrow out of $j$ which is not part of the second cycle then apply $\mu_m^{-1}$ at $j$. The possible local configurations and resultant algebras relevant to this step are shown in figures (xiv), (xv) and (xxiii)-(xxvi).
\item If there exists no such arrow, but there is an arrow into vertex $j$ which is not part of the second cycle then apply $\mu_{m}^{-1}$ at the tail of this incoming arrow, except where we encounter the configurations shown in (ixx) and (xxi) where we must apply $\mu_m$ at the tail of the arrow. The possible local configurations and resultant algebras relevant to this step are shown in figures (xvi)-(xxii).
    \\Note that not all configurations can be cleared from vertex $j$ using the above steps. In particular if we have a local configuration of the following type:
    \begin{center}
    \includegraphics[scale=0.5]{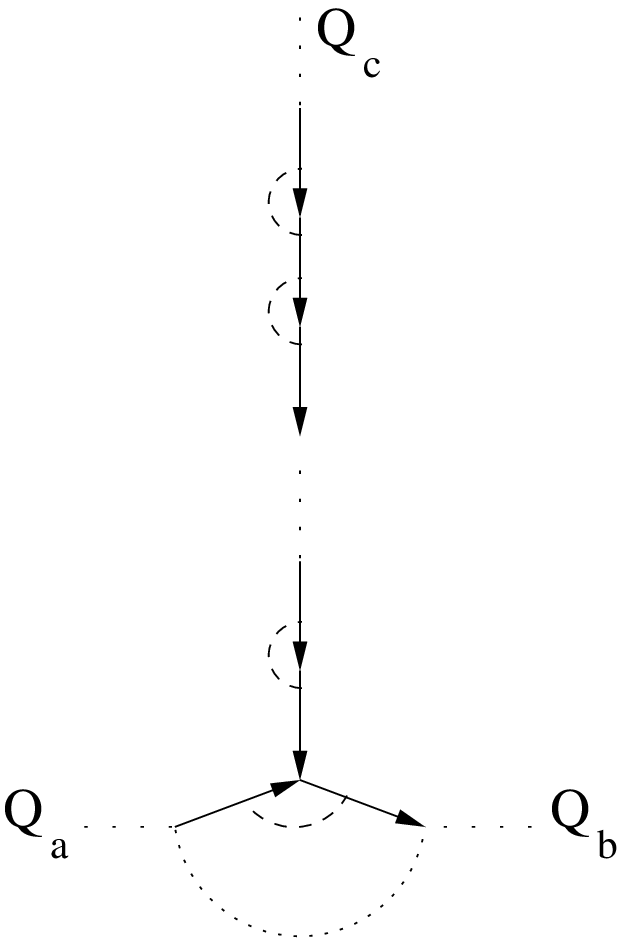}
    \end{center}
    Namely, into vertex $j$ we have an oriented path of length $s$ where $2\leqslant s\leqslant m+1$ in which each composition of two consecutive arrows is a zero relation.
    \\We can apply proposition \ref{rel_rem} to achieve a derived equivalence with an algebra with the following local configuration:
    \begin{center}
    \includegraphics[scale=0.5]{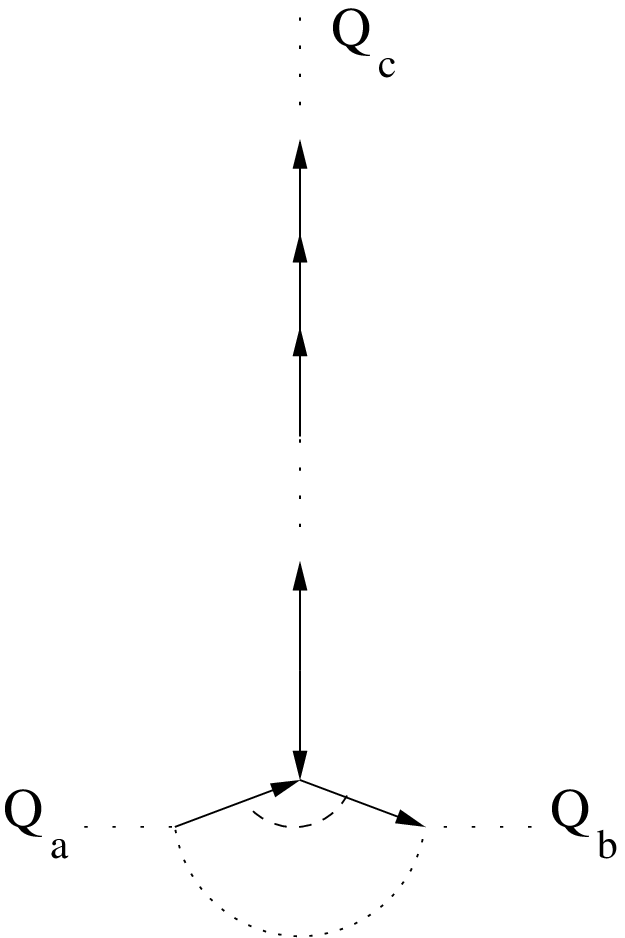}
    \end{center}
    and elsewhere has the same quiver with relations as our original algebra.
\item repeat steps (1) and (2) for vertex $j+1$, $1\leqslant j\leqslant x$ where $x=(\frac{m-3}{2})$ if $m$ is odd and $x=(\frac{m-2}{2})$ if $m$ is even.
\end{enumerate}
Once we have applied the above algorithm for vertices of type $B$ to the second cycle we have shown that $\Lambda$ is derived equivalent to an algebra whose quiver has two $m+2$-cycles which share a vertex. These two cycles are connected to other parts of the quiver only through one vertex of the second cycle in such a way that vertices of type $A$ and $B$ coincide with the definitions in \ref{nformdef}.
\\Now we simply choose a third cycle (if it exists) in exactly the same way as we chose the second cycle. The same reduction algorithm can be applied to the third cycle. Indeed, we may continue to choose and clear cycles as described above until there are no more $m+2$-cycles. Therefore, we have been able to apply successive mutations to $\Lambda$ to achieve,
\begin{center}
\includegraphics[scale=0.5]{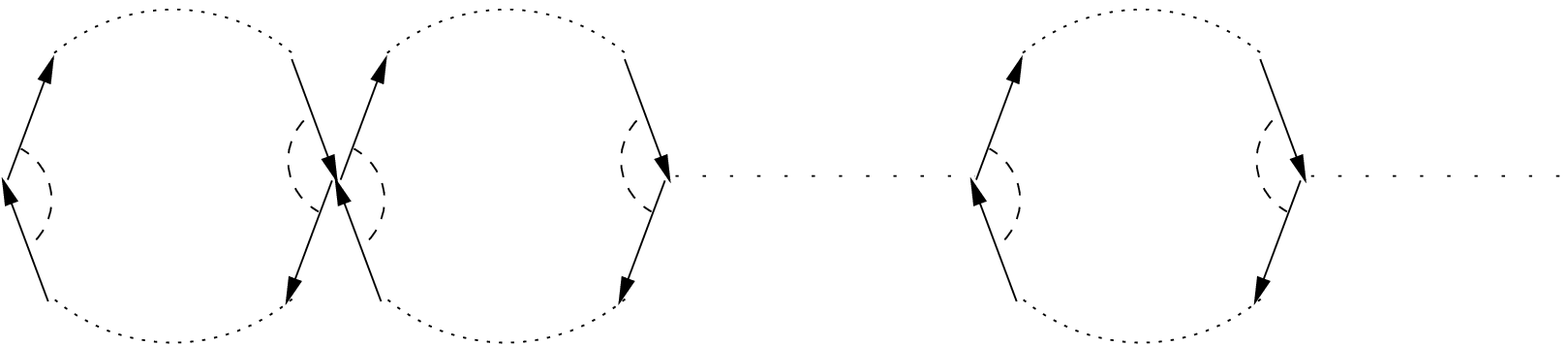}
\end{center}
so that $\Lambda$ is derived equivalent to this algebra. Let $\Lambda_2$ denote the above algebra.
Finally, it remains to prove that we can reduce the portion of the quiver which contains no $m+2$-cycles to the linear orientation with no relations, whilst maintaining derived equivalence at each step. We restrict to using the mutations induced by the elementary polygonal moves $\mu_m$ and $\mu_m^{-1}$ and so ensuring a derived equivalence at each stage by theorems \ref{tilting_reg} and \ref{tilting_dual} respectively. Notice that this part of the algorithm is the only part needed if $\Lambda$ has no $m+2$-cycles, a possibility we have ignored till now.
\\
\\[10pt]Since $\Lambda$ is gentle (see results in section \ref{sec:m_clust}) and, by our reduction algorithm, is derived equivalent to $\Lambda_2$ we have, by \cite{SZi}, that $\Lambda_2$ is gentle. Therefore, the quiver of $\Lambda_2$ which contains no $m+2$-cycles is a tree. Since it is a finite quiver, we must have that there are a finite number of end points of the underlying tree. In order to show that we can reduce the cycle-free part of the quiver of $\Lambda_2$ to the linear orientation on $A_s$ (where $s=m+2+(r-2)(m+1)+m$ and $r$ in the number of cycles in the quiver of $\Lambda_2$) we will need the following argument.
\\Suppose that we have some orientation on the graph $A_s$. The path algebra of such a quiver (where we assume that there are no relations) certainly describes an $m$-cluster-tilted algebra which has some corresponding division of the $(n+1)m+2$-gon, $P$. We claim that by using $\mu_m$ and $\mu_m^{-1}$ we can mutate the given orientation on $A_s$ to an orientation where all arrows are oriented in the same direction.
\begin{rem}For our purposes we will need to consider both possible orientations where all arrows are in the same direction. This will become clear in due course.
\end{rem}
First assume we have the following configuration occurring in some orientation of $A_s$,
\begin{center}
\includegraphics[scale=0.6]{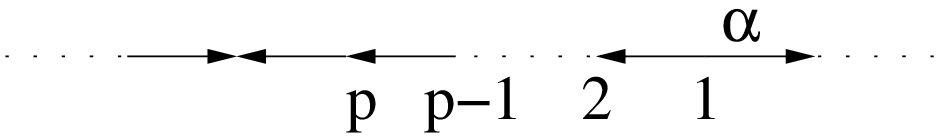}
\end{center}
Then we apply $\mu_m^{-1}$ successively to the vertices $1,2,\ldots,p-2,p-1$ to get,
\begin{center}
\includegraphics[scale=0.6]{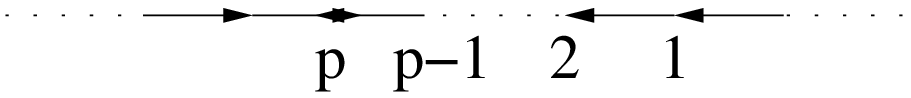}
\end{center}
\begin{rem}
Notice that it is always possible to ensure the existence of the arrow $\alpha$ in the left to right direction via a source/sink mutation if required. The effect of this series of mutations is to interchange the position of one arrow oriented from left to right with a collection of consecutive arrows oriented from right to left.
\end{rem}
Repetition of this process will result in a configuration of the form,
\begin{center}
\includegraphics[scale=0.5]{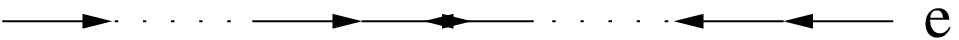}
\end{center}
Now apply $\mu_m^{-1}$ to the vertex $e$ to get,
\begin{center}
\includegraphics[scale=0.5]{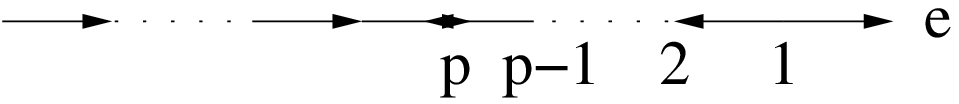}
\end{center}
Next repeat the first step, that is, apply $\mu_m^{-1}$ in turn to $1,2,\ldots,p-2,p-1$. The result is,
\begin{center}
\includegraphics[scale=0.5]{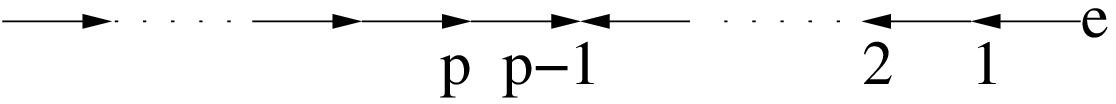}
\end{center}
but now we have effectively reduced the number of arrows oriented from right to left by 1. Hence, we can eventually achieve the orientation of $A_s$ where all arrows are oriented from left to right.
\\The proof that we mutate to the orientation where all arrows are oriented from right to left is similar, simply reverse the orientations of the arrows in the relevant figures above and apply $\mu_m$ in place of $\mu_m^{-1}$. Notice that we never mutate at the left endpoint of $A_s$. This means that we can use this process to reduce the cycle-free part of $\Lambda_2$ to the linear orientation if it has no zero relations.
\\
\\[10pt]Now we are ready to discuss relations in the cycle-free part of the quiver of $\Lambda_1$ and how to reduce the number of them. First, notice that in the cycle free part of the quiver we have that,
$$
\textnormal{number of relations }\geqslant\textnormal{ number of endpoints -1}\qquad\qquad(\ast)
$$
The equality holds where each relation occurs at a vertex of valency greater than 2. If relations occur around vertices of valency two then we will have an inequality.
\\Take any endpoint, $e$, in the cycle-free part of the quiver of $\Lambda_1$ and consider the following local situations which could occur,
\begin{center}
\includegraphics[scale=0.5]{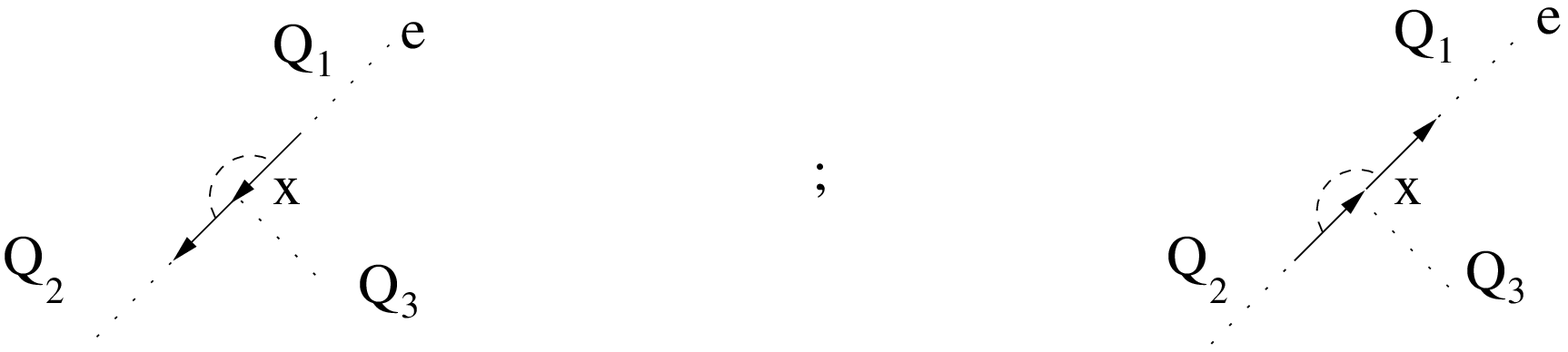}
\end{center}
where $Q_1$ consists of some orientation of $A_t$ for some $t\in\mathbb{N}$ which has no relations and $Q_2$ and $Q_3$ are arbitrary. Note that $Q_1$, $Q_2$ and $Q_3$ might all be empty, we do not discount this possibility. In effect we are considering the first relation encountered as we proceed back up the levels of the tree from an arbitrary endpoint.
\\We have the following algorithm for situations on the left in the above figure.
\\It has been demonstrated above that we can choose to orient $Q_1$ as follows and maintain a derived equivalence. Therefore, we can assume that $Q_1$ has the form shown in the following figure.
\begin{center}
\includegraphics[scale=0.5]{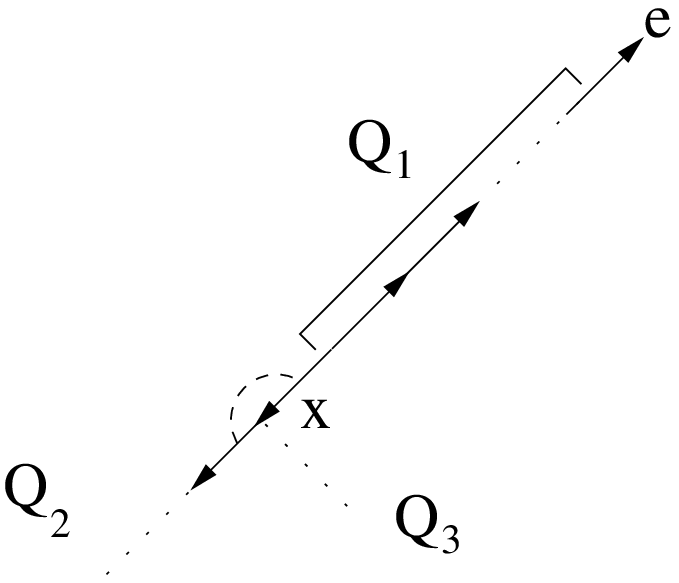}
\end{center}
Now for each arrow in $Q_1$ apply $\mu_m$ once at vertex $x$. The result is that we have mutated to,
\begin{center}
\includegraphics[scale=0.5]{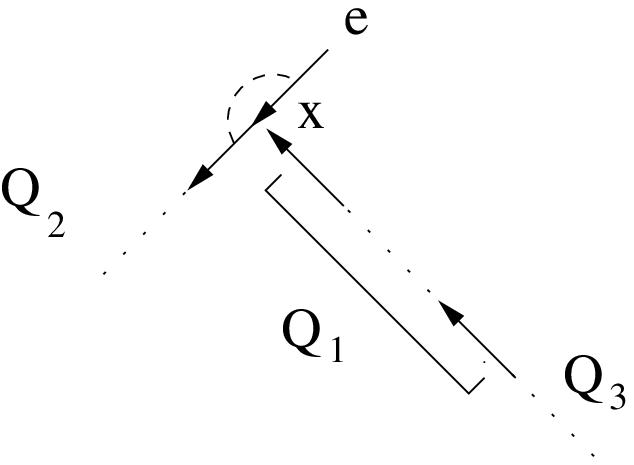}
\end{center}
and maintained the derived equivalence.
\\Now perform $\mu_m^{-1}$ at vertex $y$ to get,
\begin{center}
\includegraphics[scale=0.5]{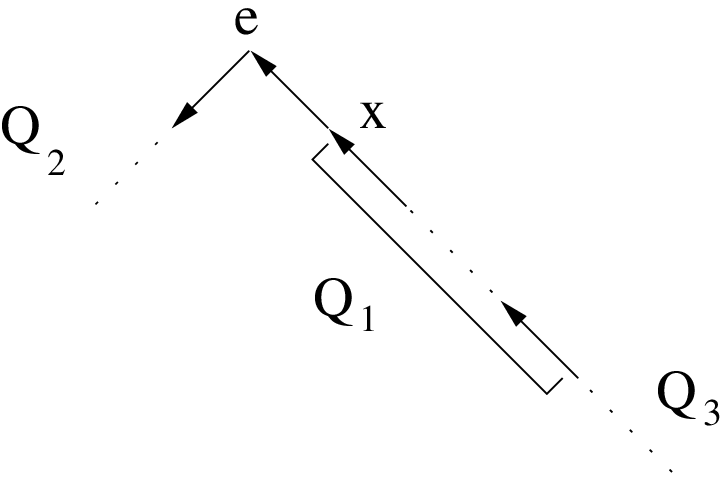}
\end{center}
It is important to note that at no stage during these mutations have we created any relations in the cycle-free part of the quiver and at each stage we have ensured a derived equivalence.
\\Now, if in the initial configuration the vertex $x$ had valency greater than 2 then we have removed a relation and reduced the number of endpoints, if the vertex $x$ had valency 2 then we have removed a relation, but the number of endpoints is preserved. Therefore, we preserve the inequality ($\ast$).
\\For the second situation where in the initial configuration we reverse the orientation of all arrows we can achieve the same result by applying $\mu_m^{-1}$ in place of $\mu_m$.
\\Finally, one can now continue to choose end points and eliminate relations. Notice that eventually the number of end points will be 1 and the number of relations will be zero. In order to achieve an orientation of the cycle-free part of the quiver in which all arrows are oriented in the same direction we apply the arguments above.
\end{proof}
\begin{rem}Some remarks should be made following the proof of the theorem.
\\It should be noted that at no stage in the algorithm do we mutate to an algebra which is not an $m$-cluster-tilted algebra, however in general it is certainly possible to do so. Using the inverse of the derived equivalence of proposition \ref{rel_rem} it is possible for an $m$-cluster-tilted algebra to be derived equivalent to an algebra which is not $m$-cluster-tilted. For example, the first algebra shown here is a 4-cluster-tilted algebra and is derived equivalent to the second algebra. The second algebra is, however, not a 4-cluster-tilted algebra since it has 4 consecutive relations not in a cycle, contrary to remark \ref{difference}.
\begin{center}
\includegraphics[scale=0.5]{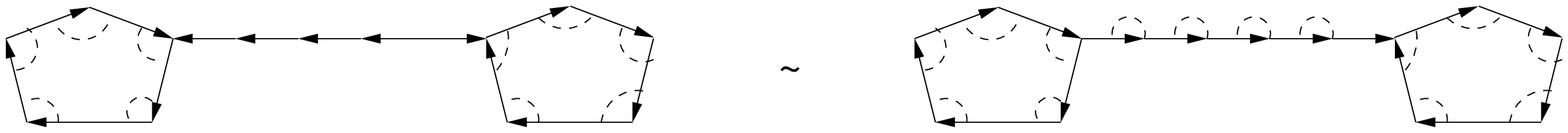}
\end{center}
So while it is known that any algebra derived equivalent to an $m$-cluster-tilted algebra of type $A_n$ is gentle, by \cite{SZi} it may not necessarily be $m$-cluster-tilted. Hence the class of $m$-cluster-tilted algebras of type $A_n$ is not closed under derived equivalence.
\\Also note that, whilst developed independently, the final part of the algorithm which removes relations from the cycle-free part of the quiver of $\Lambda_1$ features an argument similar to that used in \cite{AH} where generalized tilted algebras of type $A_n$ are considered. Our method was developed from a derived equivalence viewpoint as opposed to the classical tilting theory approach in \cite{AH} and is motivated by the polygonal combinatorics which does not feature in \cite{AH}.
\end{rem}
\begin{center}
\textsc{Acknowledgements}
\end{center}
I would like to thank my Ph.D. supervisor Dr. Thorsten Holm for many helpful discussions while preparing this paper.

\end{document}